%% file: wavelet_0_main.tex
\documentclass[review]{siamart190516}

\newif\ifarxiv

\arxivtrue

\input{header}

\input{wavelet_author}

\ifarxiv
\usepackage{filecontents}
\makeatletter\@input{wavelet_0_supp_aux.tex}\makeatother
\nolinenumbers
\fi
\definecolor{airforceblue}{rgb}{0.36, 0.54, 0.66}
\definecolor{carolinablue}{rgb}{0.6, 0.73, 0.89}
\definecolor{moonstoneblue}{rgb}{0.45, 0.66, 0.76}
\definecolor{amber}{rgb}{1.0, 0.75, 0.0}
\definecolor{cadmiumorange}{rgb}{0.93, 0.53, 0.18}
\definecolor{persianorange}{rgb}{0.85, 0.56, 0.35}
\definecolor{cadetgrey}{rgb}{0.57, 0.64, 0.69}
\definecolor{darkgray}{rgb}{0.66, 0.66, 0.66}
\definecolor{asparagus}{rgb}{0.53, 0.66, 0.42}
\definecolor{cambridgeblue}{rgb}{0.64, 0.76, 0.68}
\definecolor{olivine}{rgb}{0.6, 0.73, 0.45}
\definecolor{carnelian}{rgb}{0.7, 0.11, 0.11}
\definecolor{cardinal}{rgb}{0.77, 0.12, 0.23}
\definecolor{firebrick}{rgb}{0.7, 0.13, 0.13}

\def\figgrey{darkgray}

\begin{document}
\maketitle

\begin{abstract}
We present the derivation, implementation, and analysis of a multiresolution adaptive grid framework for numerical simulations on \revtwo{octree-based} 3D block-structured collocated grids with distributed computational architectures. 
Our approach provides a consistent handling of non-lifted and lifted interpolating wavelets of arbitrary order demonstrated using second, fourth, and sixth order wavelets, \revtwo{combined} with standard finite-difference based discretization operators. 
We first validate that the wavelet family used provides strict and explicit error control when coarsening the grid, \revtwo{and show} that lifting wavelets increase the grid compression rate while conserving discrete moments across levels. \revtwo{Further, we demonstrate that} high-order PDE discretization schemes \revtwo{combined with sufficiently high order wavelets retain the expected} convergence order even at resolution jumps.
\revtwo{We then simulate the advection of a scalar} to analyze convergence for the temporal evolution of a PDE. \revtwo{The results} shows that our wavelet-based refinement criterion is successful at controlling the overall error while the coarsening criterion is effective at retaining the relevant information on a compressed grid. 
Our software exploits \revtwo{a} block-structured grid data structure for efficient multi-level operations, \revtwo{combined with a} parallelization strategy \revtwo{that} relies on a one-sided MPI-RMA communication approach with active PSCW synchronization. \revtwo{Using performance tests up to $16,384$ cores, we demonstrate that this leads to a highly scalable performance.} \typo{The associated code is available under a BSD-3 license at \url{https://github.com/vanreeslab/murphy}.}
\end{abstract}



\input{wavelet_1_introduction}

\input{wavelet_2_method_short}

\input{wavelet_3_implementation}

\input{wavelet_4_validation}

\input{wavelet_5_results}

\input{wavelet_6_scalability}

\input{wavelet_7_conclusion}

\section*{Acknowledgments}
We would like to acknowledge the insightful discussions with James Gabbard, Matthieu Duponcheel, Pierre Balty, and Philipe Chatelain, as well as the help received from Howard Pritchard and Greg Hernandez to run the PSCW MPI-RMA on different infrastructures.
Further we wish to acknowledge the financial support from an Early Career Award from the Department of Energy, Program Manager Dr. Steven Lee, award number DE-SC0020998 (TG, WVR), the Belgian American Educational Fundation (B.A.E.F.) (TG) and from Wallonie-Bruxelles International (WBI) excellence fellowship (TG).
Computational resources have been provided by the Consortium des Équipements de Calcul Intensif (CÉCI) funded by the Fonds de la Recherche Scientifique de Belgique (F.R.S.-FNRS) under Grant No. 2.5020.11 and by the Walloon Region. \typo{Additional resources include the Tier-1 supercomputer of the Fédération Wallonie-Bruxelles, infrastructure funded by the Walloon Region under the grant agreement n°1117545, the Tier-0 LUMI accessed throught the pilot-phase and in collaboration with the CÉCI, and MELUXINA an infrastructure accessed through EuroHPC under the grant EHPC-DEV-2022D05-149. Finally this research used resources of the National Energy Research Scientific Computing Center (NERSC), a U.S. Department of Energy Office of Science User Facility located at Lawrence Berkeley National Laboratory, operated under Contract No. DE-AC02-05CH11231 using NERSC award ASCR-ERCAP20232.}

{\footnotesize

\bibliography{AllPublications.bib}
}

\end{document}


\maketitle
\input{wavelet_s2_implementation}

\input{wavelet_s4_algo}

\input{wavelet_s3_rkandweno}

\input{wavelet_s1_scalability}

{\footnotesize

\bibliography{AllPublications.bib}
}

%% file: header.tex
\usepackage[utf8]{inputenc}
\usepackage[english]{babel}
\usepackage{amsmath}
\usepackage{amsfonts}
\usepackage{amssymb}
\usepackage{amsopn}
\usepackage{bm}

\usepackage{fourier}
\usepackage[left=1.4cm,right=1.4cm,top=1.5cm,bottom=1.5cm]{geometry}

\usepackage{bm}
\usepackage{lineno}
\usepackage{hyperref}
\usepackage[boxruled,vlined,linesnumbered,algo2e]{algorithm2e} 
\usepackage{setspace}
\usepackage{xcolor}
\usepackage{wrapfig}
\usepackage{graphicx}
\usepackage{pstool}
\usepackage{multirow}
\usepackage{setspace}
\usepackage{enumitem}
\usepackage{multirow}
\usepackage{subcaption}
\usepackage{comment}

\graphicspath{{./figures/}}

\usepackage{pgfplots}
\usepgfplotslibrary{external} 
\pgfplotsset{compat=newest}

\usepackage{tikz}%
\usetikzlibrary{shapes,snakes}
\usetikzlibrary{patterns}
\usetikzlibrary{shapes.misc}
\usetikzlibrary{arrows}
\usepgfplotslibrary{groupplots}
\usepgfplotslibrary{dateplot}

\usepackage[numbers,sort&compress]{natbib}

\input{symbols_and_cmd}

\newsiamremark{remark}{Remark}
\newsiamremark{hypothesis}{Hypothesis}
\crefname{hypothesis}{Hypothesis}{Hypotheses}
\newsiamthm{claim}{Claim}

%% file: symbols_and_cmd.tex
\newcommand{\abs}[1]{\left\vert #1 \right\vert} %

\newcommand{\cfl}{\operatorname{CFL}}

\newcommand{\code}[1]{\texttt{\small{#1}}}
\newcommand{\cons}[1]{\code{CONS-#1}}

\newcommand{\fdot}[2]{\left \langle #1 \;,\; #2 \right \rangle }

\newcommand{\diff}[1]{\code{DIFF-#1}}

\newcommand{\expo}[1]{\operatorname{exp} \left[ #1 \right]}

\newcommand{\norm}[1]{\big\vert\!\big\vert #1 \big\vert\!\big\vert}
\newcommand{\normi}[1]{\vert\!\vert #1 \vert\!\vert_{\infty}}

\newcommand{\order}[1]{\mathcal{O}\left( #1 \right)}

\newcommand{\proj}[2]{P^{#2}\!\left[{#1} \right]}
\newcommand{\projval}[2]{\overline{P}^{#2}\!\left[{#1} \right]}

\newcommand{\bu}{\mathbf{u}}

\newcommand{\wave}[2]{wavelet~\code{#1.#2}}
\newcommand{\waveshort}[2]{\code{#1.#2}}

\newcommand{\x}{\mathbf{x}}

\newcommand{\dbx}{\operatorname{d}\!\mathbf{x}}
\newcommand{\dx}{\;\operatorname{d}\!{x}}

\SetKwFunction{FnGet}{MPI\_Get}
\SetKwFunction{FnPut}{MPI\_Put}
\SetKwFunction{FnPost}{MPI\_Win\_Post}
\SetKwFunction{FnStart}{MPI\_Win\_Start}
\SetKwFunction{FnComplete}{MPI\_Win\_Complete}
\SetKwFunction{FnWait}{MPI\_Win\_Wait}
\SetKwFunction{FnSyncStatus}{SyncStatus}

\SetKwFunction{FnCopy}{copy}
\SetKwFunction{FnSample}{sample}
\SetKwFunction{FnBounday}{BoundaryCondition}

\SetKwFunction{FnRefine}{WaveletRefine}
\SetKwFunction{FnCoarsen}{WaveletCoarsen}
\SetKwFunction{FnCriterion}{WaveletCriterion}
\SetKwFunction{FnOverwrite}{WaveletSubstitution}
\SetKwFunction{FnInterpolate}{WaveletInterpolate}
\SetKwFunction{FnSmooth}{WaveletUpdateAfterCoarsening}

\SetKwFunction{FnStatus}{StatusUpdate}
\SetKwFunction{FnPolicy}{EnforceStatusPolicy}

\SetKwFunction{FnTreeRefine}{p4estRefine}
\SetKwFunction{FnTreeCoarsen}{p4estCoarsen}
\SetKwFunction{FnTreeBalance}{p4estBalance}
\SetKwFunction{FnTreePartition}{p4estPartition}

\SetKwFunction{FnGhostPost}{GhostPull\_Post}
\SetKwFunction{FnGhostWait}{GhostPull\_Wait}
\SetKwFunction{FnGhost}{Ghost}

\SetKw{KwOr}{or}

\newcommand{\be}{\begin{equation}}
\newcommand{\eec}{\;\;,\end{equation}}
\newcommand{\ee}{ \end{equation}}
\newcommand{\eed}{\;\;.\end{equation}}
\newcommand{\bse}{\begin{subequations}}
\newcommand{\ese}{\end{subequations}}

\newenvironment{eqd}{\begin{equation}}{\;\;.\end{equation}}
\newenvironment{eqc}{\begin{equation}}{\;\;,\end{equation}}

\crefname{algocf}{Algorithm}{Algorithms}
\Crefname{algocf}{Algorithm}{Algorithms}

\newcommand{\Fref}[1]{\Cref{#1}}

\newcommand{\Sect}[1]{\Cref{#1}}

\newcommand{\tref}[1]{\cref{#1}}
\newcommand{\fref}[1]{\cref{#1}}
\newcommand{\eqqref}[1]{\cref{#1}}
\newcommand{\algref}[1]{\cref{#1}}
\newcommand{\frefs}[2]{\cref{#1} and \cref{#2}}
\newcommand{\sect}[1]{\cref{#1}}

\definecolor{darkGreen}{RGB}{0,160,0}
\definecolor{darkBlue}{RGB}{51,51,204}
\definecolor{darkRed}{RGB}{160,0,0}
\definecolor{orange}{RGB}{255,120,0}
\definecolor{myGreen}{rgb}{0.0706, 0.5333, 0.1921}
\definecolor{asparagus}{rgb}{0.53, 0.66, 0.42}
\definecolor{bostonuniversityred}{rgb}{0.8, 0.0, 0.0}
\definecolor{airforceblue}{rgb}{0.36, 0.54, 0.66}
\definecolor{cadmiumorange}{rgb}{0.93, 0.53, 0.18}

\ifarxiv
\newcommand{\revone}[1]{{#1}}
\newcommand{\revtwo}[1]{{#1}}

\newcommand{\typo}[1]{{#1}}
\else
\newcommand{\revone}[1]{\textcolor{airforceblue}{#1}}
\newcommand{\revtwo}[1]{\textcolor{bostonuniversityred}{#1}}

\newcommand{\typo}[1]{\textcolor{cadmiumorange}{#1}}
\fi

\definecolor{matlabBlue}{rgb}{0, 0.4470, 0.7410}
\definecolor{matlabOrange}{rgb}{0.8500, 0.3250, 0.0980}
\definecolor{matlabYellow}{rgb}{0.9290, 0.6940, 0.1250}
\definecolor{matlabPurple}{rgb}{0.4940, 0.1840, 0.5560}
\definecolor{matlabGreen}{rgb}{0.4660, 0.6740, 0.1880}
\definecolor{matlabRed}{rgb} {0.7765, 0.1569, 0.0902}
\definecolor{matlabBlack}{rgb}{0, 0, 0}

\definecolor{keyBlue}{rgb}{0.0392, 0.2274, 0.4235}
\definecolor{keyBlueTrans}{rgb}{0.4274, 0.4980, 0.5764}
\definecolor{keyGreen}{rgb}{0.05490, 0.38039, 0.003921}
\definecolor{keyGreenTrans}{rgb}{0.7019, 0.8313, 0.7176}
\definecolor{keyOrange}{rgb}{0.98039, 0.4980, 0.03529}
\definecolor{keyGray}{rgb}{0.50196, 0.50196, 0.50196}

\definecolor{pyBlue}{rgb}{0.11372,0.42352,0.67059}
\definecolor{pyBlueLight}{rgb}{0.6196,0.79215,0.88235}
\definecolor{pyOrange}{rgb}{1.0,0.4549,0.0627}
\definecolor{pyOrangeLight}{rgb}{0.9921568,0.643137,0.376470}
\definecolor{pyGreen}{rgb}{0.1529,0.5882,0.1529}
\definecolor{pyRed}{rgb}{0.81569,0.13725,0.14117}
\definecolor{pyPurple}{rgb}{0.53725,0.36078,0.7098}
\definecolor{myCacaDoie}{rgb}{0.580392156862745,0.43,0.27}
\definecolor{pyBrown}{rgb}{0.549019607843137,0.337254901960784,0.294117647058824}
\definecolor{pyRose}{rgb}{0.890196078431372,0.466666666666667,0.76078431372549}

\definecolor{pyOrange2}{rgb}{0.9921,0.5098,0.2078}
\definecolor{pyBlue2}{rgb}{0.3765,0.6431,0.8157}

\definecolor{pyBar1}{rgb}{0.1647,0.2745,0.2745}
\definecolor{pyBar2}{rgb}{0.0000,0.4588,0.4588}
\definecolor{pyBar3}{rgb}{0.4863,0.7804,0.9998}
\definecolor{pyBar4}{rgb}{0.4549,0.4549,0.4549}

\definecolor{keyGreen2}{rgb}{0.0627 , 0.4196 , 0.3882}%
\definecolor{keyOrange2}{rgb}{0.9843 , 0.2941 , 0.2431}%

\definecolor{coluni}{rgb}{0.58823,0.58823,0.58823}
\definecolor{col40}{rgb}{0.192156862745098,0.509803921568627,0.741176470588235}
\definecolor{col42}{rgb}{0.901960784313726,0.333333333333333,0.0509803921568627}
\definecolor{col20}{rgb}{0.388235294117647,0.474509803921569,0.223529411764706}
\definecolor{col22}{rgb}{0.549019607843137,0.427450980392157,0.192156862745098}
\definecolor{col60}{rgb}{0.458823529411765,0.419607843137255,0.694117647058824}
\definecolor{col62}{rgb}{0.482352941176471,0.254901960784314,0.450980392156863}

\definecolor{colweakstencil}{rgb}{0.419607843137255,0.682352941176471,0.83921568627451}
\definecolor{colweakadapt}{rgb}{0.992156862745098,0.552941176470588,0.235294117647059}
\definecolor{colweakghost}{rgb}{0.549019607843137,0.635294117647059,0.32156862745098}

\newcommand{\myLine}[3]{\raisebox{2pt}{\tikz{\draw[-,#1,#2,line width = #3 pt](0,0) -- (10pt,0);}}}
\newcommand{\ThickLine}[2]{\myLine{#1}{#2}{1}}

\newcommand{\lineSquare}[3]{
        \tikz{
            \draw[-,#1,solid,line width = #3 pt](0,0) -- (12.5pt,0);
            \filldraw[draw=#1,fill=#2, line width=#3 pt] (4.5pt,-2pt) rectangle (8.5pt,2pt);
        }
}

\newcommand{\ThickLineSquare}[2]{\lineSquare{#1}{#2}{1.0}}

\newcommand{\lineCircle}[4]{
        \tikz{
            \draw[-,#1,#4,line width = #3 pt](0,0) -- (13pt,0);
            \filldraw[draw=#1,fill=#2, line width=#3 pt] (6.5pt,0pt) circle (2pt);
        }
}
            
\newcommand{\ThickLineCircle}[2]{\lineCircle{#1}{#2}{1.0}{solid}}

\tikzset{cross/.style={cross out, draw, minimum size=2*(#1-\pgflinewidth), inner sep=0pt, outer sep=0pt}}

\newcommand{\CaptionTwoZero}{{\protect\ThickLineCircle{col20}{col20}}}
\newcommand{\CaptionTwoTwo}{{\protect\ThickLineCircle{col22}{col22}}}
\newcommand{\CaptionFourZero}{{\protect\ThickLineCircle{col40}{col40}}}
\newcommand{\CaptionFourTwo}{{\protect\ThickLineCircle{col42}{col42}}}
\newcommand{\CaptionSixZero}{{\protect\ThickLineCircle{col60}{col60}}}
\newcommand{\CaptionSixTwo}{{\protect\ThickLineCircle{col62}{col62}}}
\newcommand{\CaptionUni}{{\protect\ThickLineCircle{coluni}{coluni}}}

\newcommand{\CaptionWeakGhost}{ {\protect\ThickLineCircle{colweakghost}{colweakghost}} }
\newcommand{\CaptionWeakStencil}{ {\protect\ThickLineCircle{colweakstencil}{colweakstencil}} }
\newcommand{\CaptionWeakAdapt}{ {\protect\ThickLineCircle{colweakadapt}{colweakadapt}} }

%% file: wavelet_author.tex
\title{MURPHY - A scalable multiresolution framework for scientific computing on 3D block-structured collocated grids}

\headers{MURPHY}{Thomas Gillis and Wim M.~van Rees}
\author{Thomas Gillis \thanks{Department of Mechanical Engineering, Massachusetts Institute of Technology, 77 Massachusetts Avenue, Cambridge, MA 02139, USA} 
\and Wim M.~van Rees \footnotemark[1] \thanks{Corresponding author: \email{wvanrees@mit.edu}}}

%% file: wavelet_1_introduction.tex
\section{Introduction}

Solutions to partial differential equations (PDEs) are typically characterized by unsteady spatial scale separations. In incompressible fluid dynamics, for example, flows within a boundary layer include important structures on the smallest viscous length scales, whereas the wake is characterized by much larger, inertial structures. Moreover, these flows are intrinsically unsteady and often coupled with the motion or deformation of immersed boundaries, making the resolution requirements difficult to predict \textit{a priori}. Especially in three-dimensional \revtwo{problems}, where performance and memory constraints of computing resources \revtwo{constrain} the range of applications, it is desirable to construct methods that can adapt the local resolution of the grid to the physical requirements of the PDE solution evolved on that same grid. To achieve accurate solutions, such methods need to be able to discretize and evolve the PDE consistently across different resolutions; detect the need to refine or opportunity to coarsen; and perform the actual coarsening or refinement of the field. Further, their (parallel) implementation needs to be sufficiently efficient so that any increase in computational overhead due to the required memory access patterns, load balancing, and synchronizations does not cancel the gains from the reduction in computational elements, compared  to uniform resolution grids. Consequently, the algorithms and implementations of such adaptive grid refinement methods are significantly more complex than uniform grid methods, and have seen significant development over the last three to four decades of research in scientific computing. %

An established family of methods is formed by the patch-based adaptive mesh refinement (pAMR) approach, which considers nested overlapping grids of increasing resolution \citep{Berger:1984, Berger:1989}. This approach is prevalent across many application domains and has mature implementations inside several software frameworks such as \code{Chombo} \citep{Adams:2019}, \code{SAMRAI} \citep{Hornung:2002}, and \code{AMReX} \citep{Zhang:2019}.
Overlapping grids simplifies the use of coarsening and refinement operators, and provide straightforward integration with multigrid-based elliptic solvers. Excluding more expensive adjoint approaches \citep{Fidkowski:2011}, grid adaptation decisions in pAMR are typically made on the basis of either a heuristic measure of the field values and/or their derivatives, or an estimate of the truncation error through a Richardson extrapolation technique \citep{Berger:1989}, or some combination of both \citep{Roy:2009}. If the chosen measure is larger than a user-defined threshold, the grid is refined, and if it is smaller than a second user-defined threshold, the grid is coarsened. \revtwo{Such} heuristic measure\revtwo{s} \revtwo{are} easy to implement but \revtwo{generally} provide no \textit{a priori} sense of the error made. The Richardson extrapolation \revtwo{method} is more rigorous but also more challenging to implement \citep{Adams:2019}, and requires \revtwo{the simultaneous} evolution of the discrete equations on multiple levels. After the decision to adapt a grid has been made, the actual coarsening/refinement approaches \revtwo{and the evaluation of differential operators across resolution boundaries} are most often based on polynomial interpolation. Most approaches achieve second-order accuracy in space throughout these operations, though extensions to fourth order have been demonstrated  as well \citep{Zhang:2012, Hooft:2022}. 

\revtwo{Our work falls under a} different family of methods which does not consider overlapping grids but instead uses an octree-based approach. This was demonstrated first in \code{Gerris} \citep{Popinet:2003} and has a more general distributed implementation in \code{p4est} \citep{Burstedde:2011}, which itself is successfully used across different application clients such as the finite-volume solver \code{ForestClaw} \citep{Calhoun:2017}. 
The dyadic recursive structure associated with these approaches can be combined with a wavelet-based multiresolution analysis of any signal on the grid  \citep{Mallat:1989,Meyer:1993}. Combining wavelet-based grid adaptation with node-based collocated PDE solution methods such as finite-difference techniques leads to the wavelet-collocation method \citep{Harten:1995, Cohen:2002, Schneider:2009}. Wavelet collocation methods typically rely on interpolating wavelets, which distinguish the information between two levels through the deviation of fine-level values from an interpolating polynomial constructed from coarse-level values \citep{Donoho:1992a}. The interpolating wavelets have been cast in a formal basis through the addition of the biorthogonality \citep{Cohen:1992} and extended to include moment-preservation and reduce aliasing through the lifting scheme \cite{Sweldens:1996,Daubechies:1998,Sweldens:1998}. The lifted interpolating wavelets form the basis of so-called second generation wavelet collocation methods \citep{Vasilyev:2000}. Existing codes using wavelet-based grid adaptation with finite-difference based PDE evolutions are \code{MRAG} \citep{Rees:2014, Rossinelli:2015}, which uses non-lifted interpolating wavelets for incompressible flow simulations with shared-memory parallelism, and \code{wabbit} \citep{Engels:2021}, using non-lifted and lifted interpolating wavelets on distributed memory architectures for weakly compressible flow simulations. Overall, despite the potential advantages of high-order grid adaptation, formal multiresolution analysis, explicit error control, and parallel performance, there is a significantly smaller body of work on the implementation details and performance analysis of non-lifting and lifting wavelets on block-structured grids, as compared to pAMR methods. 

In this work, we describe the derivation, implementation, and analysis of a parallel, scalable implementation of a 3D multiresolution adaptive grid solver for partial differential equations on collocated grids, supported by non-lifted or lifted wavelets for scale detection, grid adaptation, and ghost reconstruction. In \sect{sect_multires_global} we provide a brief background of the wavelet theory and show how this translates to block-structured multi-level grids. We emphasize here the consistent treatment of resolution jumps to preserve polynomial order and lifting properties \revtwo{in multiple dimensions}. \Sect{sect_implementation_global} details our implementation, using one-sided MPI-RMA communication strategies to handle parallel communication.
In \sect{sect_validation_global} \revtwo{we validate our approach on static grid adaptation tests} by demonstrating error control and convergence of high-order finite-difference schemes for all wavelets, and moment conservation for lifted wavelets.
\Sect{sec:results_fd} applies \revtwo{the resulting software} to solve \revtwo{partial differential equations}, where we provide detailed analysis of the error as a function of the wavelet-based thresholding parameters. \revtwo{Here we limit ourselves to linear and nonlinear scalar advection equations as examples of challenging problems involving dynamically changing scales, though our framework as presented can already handle a wider set of problems including advection-diffusion and reaction-diffusion equations, and will be further extended in future work.}
In \sect{sec:scalability} we demonstrate that our code retains parallel efficiency across \typo{more than sixteen thousand compute cores} and conclude our work with a perspective and future work in \sect{sec:conclusion}.

%% file: wavelet_2_method_short.tex
\section{Wavelet-based multiresolution}
\label{sect_multires_global}
Our work relies on a few key contributions that have been made to the field of wavelets theory that include the multiresolution analysis, biorthogonal interpolating wavelets, and the lifting scheme. Though a complete overview of wavelet theory is beyond the scope of this \revtwo{manuscript, we provide} a concise review of the concepts required to detail the mathematical framework for our multiresolution grid adaptation below. 

\subsection{Interpolating wavelets}

Throughout our work we use interpolating wavelets. These were first introduced by \citep{Donoho:1992a}, and generalize the polynomial interpolation procedure on nested dyadic grids presented in \citep{Deslauriers:1987,Deslauriers:1989} using wavelet theory, as detailed in \citep{Donoho:1996}.
\revtwo{Interpolating wavelets can be constructed } through \revtwo{polynomial} interpolation, thus avoiding the Fourier transform.
Through the introduction of the orthogonal multiresolution analysis first, and the biorthogonal multiresolution analysis second, they have been formalized within the framework of second-generation wavelets, and through the lifting scheme they can be generalized to achieve moment-conservation properties and reduce aliasing. Here we will briefly touch upon these concepts and their mathematical background.

\subsubsection{The orthogonal multiresolution analysis}
The multiresolution analysis \citep{Mallat:1989a, Cohen:1992,Sweldens:1996} defines nested orthogonal subspaces decomposition of $L_2(\mathbb{R})$ indexed by $L$. Mathematically the nested spaces can be written as $V^{L} \subset V^{L+1}$ for $L \in \mathbb{Z}$, where the union $\cup_{L \in \mathbb{Z}} V^L$ is dense and orthogonality implies that the intersection $\cap_{L \in \mathbb{Z}} V^L$ is empty. The subspaces are further defined with dilation and translation characteristics that guarantee the existence 
of a unique function $ \varphi (x) $, such that for any $L \in \mathbb{Z}$ the translated and dilated family of functions $\varphi_k^L(x) = \sqrt{2^L} \; \varphi (2^L \; x - k )$ for $k \in \mathbb{Z}$ is an orthonormal basis of $V^L$ \citep{Mallat:1989a}. 
Orthonormality here means that $\fdot{\varphi_k^L(x)}{\varphi_i^L(x)} = \delta_{i,k} $ $\forall \left\{k,i\right\} \in \mathbb{Z}$, where $\fdot{f(x)}{g(x)} = \int_{-\infty}^{\infty} f(x) g(x) \dx$.
The difference between two spaces $V^L$ and $V^{L+1}$ is characterized by a new subspace $W^L$ as the orthogonal complement of $V^L$ to $V^{L+1}$, so that $W^{L} = \left( V^{L} \cap V^{L+1} \right)^{\perp}$, hence $V^{L} \oplus W^{L} = V^{L+1}$.
Similarly to $V^L$, a basis for $W^L$ is obtained through the dilatation and translation of the wavelet function $\psi(x)$, such that $\psi^{L}_m(x) = \sqrt{2^L} \; \psi(2^{L} x - m)$ with $m \in \mathbb{Z}$.

With these definitions, a given function $f(x) \in L_2(\mathbb{R})$ can be projected onto either basis to define the scaling coefficients $\lambda_k^{L}$ and detail coefficients $\gamma_k^{L}$
\be
\lambda_k^{L} \triangleq \fdot{f(x)}{\varphi^{L}_k(x)} \qquad \text{and} \qquad \gamma_m^{L} \triangleq \fdot{f(x)}{\psi^{L}_m(x)}
\eed
\revtwo{We can then} build a hierarchy of projections of $f(x)$ into the wavelet subspaces. We start with the projection of $f(x)$ onto level $L$ denoted as 
\be
\proj{f}{L}(x) \triangleq \sum_{k} \lambda_k^{L} \; \varphi^{L}_k(x)
\eec
Further, given that $V^{L} \oplus W^{L} = V^{L+1}$, we can relate the projection of $f(x)$ onto level $L+1$ to lower levels through the refinement relation:
\be
\proj{f}{L+1} (x) \triangleq \sum_{j} \lambda_j^{L+1} \; \varphi^{L+1}_j(x) = \sum_{k} \lambda_k^{L} \; \varphi^{L}_k(x) \;+\; \sum_{m} \gamma^{L}_m \; \psi^{L}_m(x) 
\label{eq_refinement_relation}
\eed
Applied recursively, \eqqref{eq_refinement_relation} can be used to create a hierarchy of nested decompositions from level $L_0$ to level $L$:

\begin{eqd}
\proj{f}{L+1} (x) = \sum_j \lambda^{L+1}_j \varphi^{L+1}_j(x) =  \sum_{k} \lambda^{L_0}_k \; \varphi^{L_0}_k(x) \;+\; \sum_{L_0 \leq l \leq L} \sum_{m} \gamma^{l}_m \; \psi^{l}_m(x)
\label{eq_multires}
\end{eqd}

\subsubsection{Biorthogonality and linear filters}
To generalize the multiresolution analysis to a broader class of wavelet functions such as the symmetric or interpolating ones, one can relax the criteria for finding scaling and corresponding wavelet functions using the dual multiresolution analysis based on biorthogonality \cite{Cohen:1992, Sweldens:1998}.
With biorthogonal wavelets, the basis of $V^L$, $\varphi_k^L(x)$, does no longer need to be orthonormal.
Instead, one uses another subspace $\tilde{V}^L$ and the associated basis functions $\tilde{\varphi}_k^L(x)$ such that:
\begin{eqd}
\fdot{\varphi_i^L(x)}{\tilde{\varphi}_k^L(x)} = \delta_{i,k}  \quad \text{with} \quad \left\{i,k\right\} \in \mathbb{Z}
\end{eqd}
The spaces $V^L$ and $\tilde{V}^L$ have non-orthogonal complements $W^L$ and $\tilde{W}^L$, respectively, such that $V^L \perp \tilde{W}^L$ and $\tilde{V}^L \perp W^L$.
This leads to the definition of the primal (dual) scaling functions $\varphi$ ($\tilde{\varphi}$), and the primal (dual) wavelet functions $\psi$ ($\tilde{\psi}$), which form bases of their respective subspaces and satisfy
\begin{eqd}
\fdot{\tilde{\varphi}_i^L(x)}{\psi_k^L(x)} = \fdot{\tilde{\psi}_i^L(x)}{\varphi_k^L(x)} = 0 \qquad \text{and}\qquad \fdot{\tilde{\varphi}_i^L(x)}{\varphi_k^L(x)} = \fdot{ \tilde{\psi}_i^L(x)}{\psi_k^L(x)} = \delta_{i,k}
\end{eqd}
The definitions of the scaling and detail coefficients become
\begin{eqc}
\lambda_k^{L} \triangleq \fdot{ f(x) }{ \tilde{\varphi}^{L}_k(x) } \qquad \text{and} \qquad \gamma_m^{L} \triangleq \fdot{ f(x)}{\tilde{\psi}^{L}_m(x)}
\end{eqc}
and the refinement relation remains unchanged
\begin{eqd}
\proj{f}{L+1} (x) = \sum_{j} \lambda_j^{L+1} \; \varphi^{L+1}_j(x) = \sum_{k} \lambda_k^{L} \; \varphi^{L}_k(x) \;+\; \sum_{m} \gamma^{L}_m \; \psi^{L}_m(x)
\end{eqd}

Following the nested subspace decomposition a linear filter can be associated to each primal/dual basis function \citep{Sweldens:1998}, binding two levels together
\be
\varphi^{L}_k(x) = \sum_j h_{k,j} \; \varphi^{L+1}_{j} (x) \text{,}\qquad
 \psi^{L}_m (x) = \sum_n g_{m,n}  \; \varphi^{L+1}_{n} (x) \text{,}\qquad
 \tilde{\varphi}^{L}_k (x) = \sum_j \tilde{h}_{k,j}   \; \tilde{\varphi}^{L+1}_{j} (x)  \qquad \text{and} \qquad
 \tilde{\psi}^{L}_m (x) = \sum_n \tilde{g}_{m,n}  \; \tilde{\varphi}^{L+1}_{n} (x)
\eed
Combined with the biorthogonal refinement relation, the filters provide the relations for the forward wavelet decomposition (also known as the analysis operation):
\be
\lambda^{L}_k = \sum_j \tilde{h}_{k,j} \lambda^{L+1}_j \;\triangleq\; \tilde{H}_{k,j} \; \lambda_j^{L+1} \qquad \text{and} \qquad \gamma^{L}_m = \sum_j \tilde{g}_{m,j} \lambda^{L+1}_j \;\triangleq\; \tilde{G}_{m,j} \; \lambda^{L+1}_j
\eec
where we used the Einstein summation convention to simply the notation. The inverse wavelet decomposition (also known as synthesis operation) is obtained \revtwo{as}
\be
\lambda^{L+1}_j = \sum_k h_{j,k} \; \lambda^{L}_k \quad+\quad \sum_m g_{j,m} \; \gamma^{L}_m \;\triangleq\; H_{j,k}\; \lambda^{L}_{k} + G_{j,m} \; \gamma^{L}_m
\eed

Both the forward and inverse transforms have linear computational complexity in the number of degrees-of-freedom and are easily represented as a block diagram, as illustrated in \fref{block_fwt_iwt_lifting_b}.

 \subsubsection{Interpolating wavelets}
 \label{sec:donoho_interp_wavelets}
Our work \revtwo{relies on} the interpolating wavelets first proposed by \citep{Donoho:1992a} based on the Deslauries-Dubuc interpolation filters \citep{Deslauriers:1987,Deslauriers:1989}.
This wavelet family provides a non-orthogonal basis \cite{Donoho:1992a} but their construction can be framed within the context of a biorthogonal multiresolution analysis as detailed in \citep{Sweldens:1996,Sweldens:1998}.%

Like the orthogonal wavelets, interpolating wavelets are characterized by the dilatation and translation of a (non-orthonormal) scaling function $\varphi_k^L = \varphi (2^{L} x - k )$ \cite{Donoho:1992a}. The interpolating property of the scaling function is given by $\varphi (i - k ) = \delta_{i,k}$ for $i,k \in \mathbb{Z}$. It is convenient at this point to define $x_{L,k} \triangleq k 2^{-L}$ as the coordinate associated to an index $k \in \mathbb{Z}$ at level $L$ such that $\varphi_k^L(x_{L,i}) = \varphi(i - k)$. The interpolating nature then implies that the evaluation of the projection of the function at $x_{L,i}$ is equal to the associated scaling coefficient at that coordinate, $\lambda_i^L$:
\begin{eqd}
\proj{f}{L} (x_{L,i})  = \sum_k \lambda^L_k \; \varphi^L_k(x_{L,i}) = \lambda_i^L
\label{eq:scaling_interp}
\end{eqd}
Interpolating wavelets can be classified by their degree of interpolation $N$, which corresponds to the number of moments of the scaling function,
\begin{eqd}
\label{eq:moments_scaling_function}
\int_{-\infty}^{\infty} x^p \; \varphi(x) \dx = \delta_p \qquad 0 \leq p < N
\end{eqd}
This relation guarantees the ability of the scaling functions to exactly reproduce polynomials of order $N-1$.

Interpolating wavelets are \revtwo{well suited to} wavelet collocation methods because it is convenient to use function evaluations at $x_{L,k}$  interchangeably with scaling coefficients $\lambda_k^L$. However, in general this does incur an error associated with the truncation of detail coefficients at any given level. This error can be captured by comparing the exact function with the function $\projval{f}{L}(x)$ defined as $\projval{f}{L}(x) \triangleq \sum_k f(x_{L,k}) \; \varphi_k^L(x)$. \revtwo{This corresponds to a similar projection as \eqref{eq:scaling_interp},} but replacing $\lambda^L_k$ with $f(x_{L,k})$. \revtwo{A bound on the difference} can be found as \citep{Sweldens:1994}:
\begin{eqd}
\abs{\projval{f}{L}(x) - f(x)}  \leq \order{2^{-L \; N}} %
\label{eq:interp_error_bound}
\end{eqd}
Specifically at location $x_{L,k}$ we then find $\abs{f(x_{L,k})- \lambda^L_k} \leq \order{2^{-L \; N}}$. Since \revtwo{$2^{-L} = x_{L,k+1} - x_{L,k} \sim h$ relates to the grid spacing on level $L$}, this \revtwo{relation implies } that using function values in an $N$th order interpolating wavelet-based projection incurs a discretization error of $\mathcal{O}(h^N)$.

The simplest family of interpolating wavelets are the Donoho interpolating wavelets \citep{Donoho:1992a,Sweldens:1996,Sweldens:1998}, which are classified here with the code $N.0$, with $N$ the degree of interpolation. %
For the Donoho interpolating wavelets, the dual scaling function is a Dirac impulse located at the origin \citep{Sweldens:1996}, $\tilde{\varphi}_k^L(x) = \delta(x-x_{L,k})$. Hence, it follows that at any level
\begin{eqc}
\lambda_k^L = \fdot{f(x)}{\tilde{\varphi}^L_k(x)} = f(x_{L,k})
\label{eq:donoho_interp}
\end{eqc}
which means the scaling coefficients at level $L$ do not just equate the function \textit{projection} evaluation at level $L$, as in \eqref{eq:scaling_interp}, but they equate the function evaluation itself: $\projval{f}{}(x_{L,k}) = \proj{f}{}(x_{L,k}) = f(x_{L,k}) = \lambda^L_k$. However, these wavelets do not conserve moments when compressing information, and are characterized by considerable aliasing in the wavelet transform as reported in \citep{Sweldens:1996}. One potential avenue to address these issues is by increasing $\tilde{N}$, the number of zero moments of the wavelet function. This can be done through the lifting approach proposed in \citep{Sweldens:1996} as discussed in the next section.

\subsubsection{Lifted interpolating wavelets}
\begin{figure}[h!]
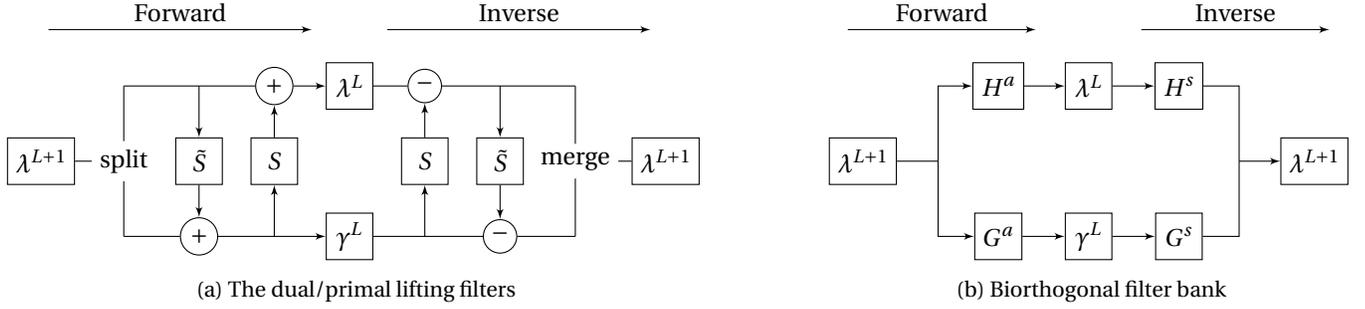

\begin{minipage}{0.58\textwidth}
\centering
\InputIfFileExists{figures/fwt_iwt_lifting_tikz}{}{}
\subcaption{The dual/primal lifting filters}
\label{block_fwt_iwt_lifting_a}
\end{minipage}%
\hfill%
\begin{minipage}{0.38\textwidth}
\centering
\InputIfFileExists{figures/ftw_iwt_tikz}{}{}{}{}
\subcaption{Biorthogonal filter bank}
\label{block_fwt_iwt_lifting_b}
\end{minipage}%
\caption{The forward and inverse wavelet transform as the combination of the dual and primal lifting steps, expressed through the dual/primal filters (a) and through the full biorthogonal filter bank (b).}
\label{block_fwt_iwt_lifting}
\end{figure}

The lifting scheme has been introduced as a general and convenient way to construct biorthogonal and second generation wavelets and their associated linear filters \citep{Sweldens:1996,Sweldens:1996a,Daubechies:1998,Sweldens:1998}.
Although the scheme can be generalized to any wavelet family, we \revtwo{restrict ourselves here to the interpolating wavelet}.
Starting with a set of scaling coefficients on level $L+1$, the lifting scheme uses the following three different steps to obtain the set of scaling and detail coefficients on level $L$:
\begin{enumerate}
    \item The \textit{splitting} step splits the fine scaling coefficients into temporary coarse scaling (even indices) and detail (odd indices) coefficients. This step is also known as the application of the ``Lazy wavelet'', and can be captured by the filters $\tilde{H}_{k,j} = \delta_{2k,j}$ and $\tilde{H}_{m,j}= \delta_{2m+1,j}$. After this first step, we have a set of coarse scaling and detail coefficients
    \begin{eqd}
    \lambda^{L}_k = \lambda^{L+1}_{2k} \quad , \quad     \gamma^{L}_m = \lambda^{L+1}_{2m+1}
    \end{eqd}
    \item The \textit{dual lifting} applies the filter $\tilde{S}$ to the scaling coefficients and uses the result to update the detail coefficients: 
    \begin{eqd}
    \gamma^{L}_m \leftarrow \gamma^{L}_{m} + \tilde{S}_{m,k} \; \lambda^{L}_k
    \end{eqd}
    \item The \textit{primal lifting} applies the filter $S$ to the detail coefficients and uses the result to update the scaling coefficients:
    \begin{eqd}
    \lambda^{L}_k \leftarrow \lambda^{L}_{k} + {S}_{k,m}  \gamma^{L}_m
    \end{eqd}
\end{enumerate}
The successive application of the three steps is illustrated in \fref{block_fwt_iwt_lifting}, and can be expressed through composite filters $H^a$ and $G^a$ that combine all stages into single operators. Reversing the sequence of operations and individual steps leads to the corresponding inverse transform, captured by the $H^s$ and $G^s$ filters.

The Donoho interpolating wavelets discussed in the previous section can be cast in the format of the lifting scheme by setting $\tilde{S}$ to the filter coefficients of \citep{Deslauriers:1987,Deslauriers:1989}, which ensure exact interpolation for polynomials of degree up to $N-1$, and setting all primal lifting filter coefficients $\tilde{S} = 0$. 

To `lift' these wavelets, following \citep{Sweldens:1996}, one can choose the primal filter $\tilde{S}$ in such a way that the first $\tilde{N}$ moments of the primal wavelet function vanish:
\begin{eqd}
\int_{-\infty}^{\infty} x^p \; \psi_k^L(x) \dx = 0 \qquad 0 \le p < \tilde{N}
\label{eq:choice_stilde}
\end{eqd}
If this holds, we ensure the conservation of the first $\tilde{N}$ moments across levels:
\begin{eqc}
\int_{-\infty}^{\infty} x^p  \; \left[ \sum_{k} \lambda^{L+1}_k \; \varphi^{L+1}_k(x) \right] \dx = \int_{-\infty}^{\infty} x^p  \; \left[ \sum_{k} \lambda^{L}_k \; \varphi^{L}_k(x) \;+\; \sum_{m} \gamma^{L}_m \; \psi^{L}_m(x) \right] \dx = \int_{-\infty}^{\infty} x^p  \; \left[ \sum_{k} \lambda^{L}_k \; \varphi^{L}_k(x) \right] \dx
\label{eq:moment_conservation_identity}
\end{eqc}
for $0 \le p < \tilde{N}$. 
Using the moment properties on $\varphi(x)$ from the interpolating wavelet, the $p$th moment on level $L$ can be expressed as: 
\begin{eqc}
\int_{-\infty}^{\infty} x^p \proj{f}{L}(x) \dx = \sum_k \; \lambda_k^L \int_{-\infty}^{\infty}  x^p \; \varphi^{L}_k(x) \dx = 2^{-L} \; \sum_k \; \lambda_k^L \; \left( x_{k,L}\right)^p
\label{eq:moment_interpolating_1d}
\end{eqc}
which allows us to rewrite the moment conservation identity of lifted interpolating wavelets \eqqref{eq:moment_conservation_identity} between two levels $L+1$ and $L$ as
\begin{eqc}
\sum_j \; \lambda_j^{L+1}  \left( x_{L+1,j} \right)^p = 2 \; \sum_k \; \lambda_k^{L} \; \left(x_{L,k} \right)^p
\label{eq:moment_conservation_measure}
\end{eqc}
for $0 \le p < \tilde{N}$. 

To satisfy \eqqref{eq:choice_stilde}, on a uniform grid and assuming $\tilde{N} \leq N$ it can be shown that $S^{\tilde{N}} = - 1/2 \; \tilde{S}^{\tilde{N}}$ \citep[theorem 12]{Sweldens:1996}, where $\tilde{S}$ are the dyadic interpolation coefficients from \citep{Deslauriers:1987,Deslauriers:1989} and given in \tref{table_duallifting_coef}. The lifting scheme thus results in an interpolating filter bank for interpolating wavelets indexed by the corresponding interpolation and moment properties $N.\tilde{N}$, which can be used for refinement ($\tilde{H} = H^a$ and $\tilde{G} = G^a$) and coarsening ($H = H^s$ and $G=G^s$) operations. Non-lifted (Donoho) interpolating wavelets have $\tilde{N} = 0$, whereas lifted wavelets have $\tilde{N} > 0$; numerical values for the filters are given in \sect{sect:a0_implementation}. All our results in this work are restricted to \revtwo{$N \in \{2,4,6\}$ and $\tilde{N}\in\{0,2\}$.} 

\begin{table}[h!]
\centering
\begin{tabular}{r|cccccc}
N  & $S_{-2}$ 	& $S_{-1}$ 	& $S_0$	& $S_1$	& $S_2$ 	& $S_3$ 	\\ \hline
$2$ & 		& 			&  $-1/2$ 	& $-1/2$  	&		&		\\
$4$ &  		& $ 1/16 $ 	& $-9/16$	& $-9/16$ 	& $ 1/16$			\\
$6 $& $-3/256$ &$ 25/256 $	& $-75/128$ & $ -75/128 $ & $25/256 $ &$ -3/256$\\
\end{tabular}
\caption{Dual-lifting coefficients, $\tilde{S}_i$ \protect\citep[Appendix B.3]{Bernard:1999}. For uniform grids the filter is symmetric, \textit{i.e.} $\tilde{S}_{-i} = \tilde{S}_i$.}
\label{table_duallifting_coef}
\end{table}

Lastly, we note that lifting the interpolating wavelets leaves the primal scaling function $\varphi(x)$ unaffected, but does change the dual scaling function $\tilde{\varphi}(x)$ from a delta function to a continuous distribution. This means the identity $ f(x_{L,k}) = \lambda^L_k$ of the Donoho interpolating wavelets is lost, and instead we fall back on the general error bound  $\abs{\lambda^L_k - f(x_{L,k})} \leq \order{2^{-L \; N}}$ provided in \eqqref{eq:interp_error_bound}.

\subsubsection{Compression}
\label{sec:wavelet_compression}
Relying on the multiresolution theory, \revtwo{compression can be achieved by discarding all the detail coefficients whose absolute values are smaller than a} tolerance $\epsilon$. 
This yields a coarser (or compressed) representation of the information,
\begin{eqd}
\proj{f}{L}_{\epsilon} (x)= \sum_{k} \lambda^{L_0}_k \; \varphi^{L_0}_k(x) \;+\; \sum_{L_0 \leq l \leq L \; } \sum_{\abs{\gamma_m} > \epsilon} \gamma^{l}_m \; \psi^{l}_m(x)
\label{eq_multires_epsilon}
\end{eqd}
It can be shown \cite{Donoho:1992a,Vasilyev:2000, Kevlahan:2005} that the error committed by this approximation is of the order of $\epsilon$,
\be
\norm{ \proj{f}{L} (x) - \proj{f}{L}_{\epsilon} (x) }_{\infty} \leq C_1 \epsilon
\label{eq:error_mr}
\eec
where $C_1$ depends on $f(x)$.
In practice, with a reasonably smooth function the value of $C_1 \approx 1$, which means that $\epsilon$ is an accurate estimate of the local error committed.%

\subsection{Extension to block-structured grids}
\label{sec:block_struct}

In this section we describe how to adapt the multiresolution analysis described above to block-structured grids. Throughout this work, we limit the jump of resolution between two adjacent blocks to $1$ (denoted as the 2:1 constraint), which considerably simplifies the operations and the implementation complexity. Below we discuss three fundamental operations that are required in the implementation: \textit{coarsening} describes the compression of data in $2^d$ blocks into a single block at the next lower resolution level, with $d$ the number of spatial dimensions of the grid; \textit{refinement} describes the refinement of data in a single block into $2^d$ new blocks at the next higher resolution level; and \textit{ghost point reconstruction} relates to the construction of ghost points for blocks across a resolution jump, so that finite-difference stencils can be evaluated on each block at its own local uniform resolution. 

We describe each of these operations in $1$D in more detail below, focusing on the \wave{2}{2} wavelet for simplicity. Subsequently, we will discuss how implementation choices lead to the treatment of grid points near resolution jumps, and how we define the criteria for compressing or refining a block. We note beforehand that our sketches use grid `blocks' and associated numbering that do not reflect a practical setting, but rather provide the minimum number of points needed to explain the respective operations for ease of interpretation.

\subsubsection{Coarsening}
\label{sec_mr_coarsening}
\begin{figure}[h!]
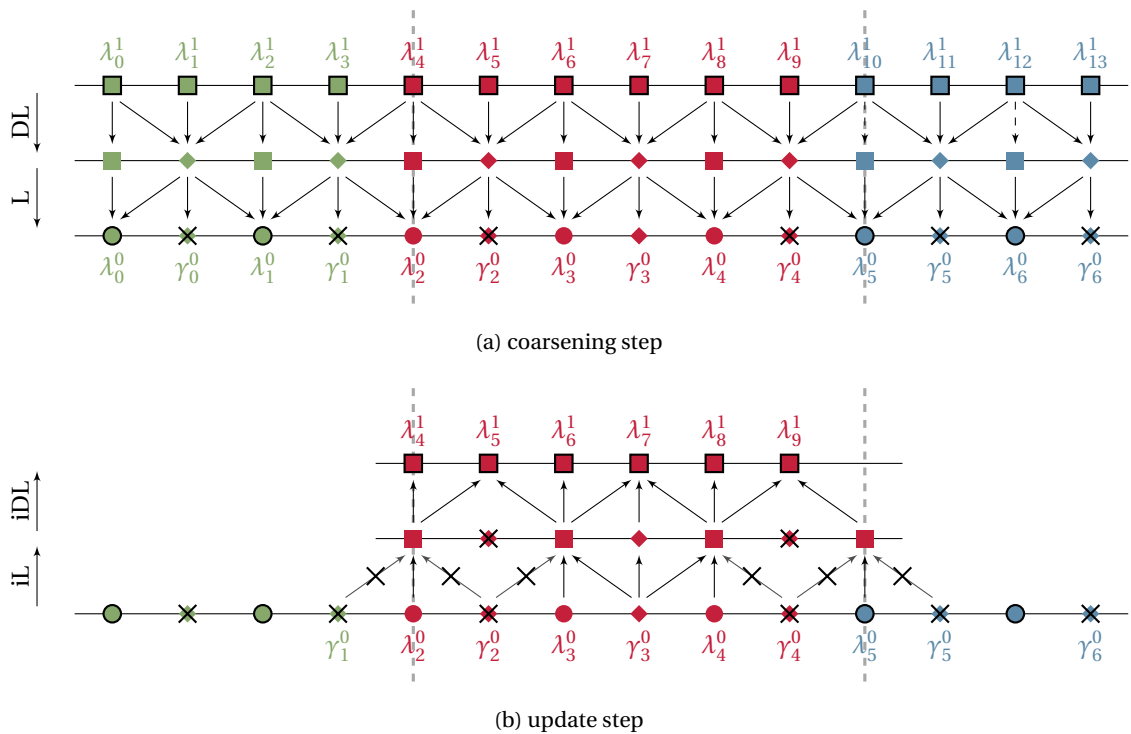

    \centering
    \begin{minipage}{\textwidth}
    \centering
    \InputIfFileExists{figures/w2.2_1d_coarsen}{}{}
    \subcaption{coarsening step}
    \label{fig:2.2_coarsen_2_a}
    \end{minipage}
        
    \begin{minipage}{\textwidth}
    \centering
    \InputIfFileExists{figures/w2.2_1d_smoothing}{}{}
	 \subcaption{update step}
	 \label{fig:2.2_coarsen_2_b}
    \end{minipage}%
\caption{Sketch of the steps required to coarsen the left (green) and right (blue) grid regions using \wave{2}{2}, with the central (red) region remaining at fine resolution. First we apply the dual lifting (DL) and lifting (L) steps successively to compute the coarse-level scaling coefficients, and discard the neglected detail coefficients (a). When $\tilde{N} > 0$ the information encoded in the discarded detail coefficients must also be removed from the fine-level information on the central (red) region, by using the inverse lifting (iL) and inverse dual lifting (iDL) filters during the update step (b). }
    \label{fig:2.2_coarsen_2}
\end{figure}

Starting from a uniform resolution one can coarsen a block \revtwo{using} the filter $H_a$. The coarsening pattern is illustrated in \fref{fig:2.2_coarsen_2} for \wave{2}{2} where the `left' (green) and `right' (blue) fine scaling coefficients at level $L+1 = 1$ (top row) are converted into coarse scaling coefficients at level $L = 0$ (bottom row) through subsequent application of the dual lifting ($\tilde{S}$) and lifting ($S$) filters. After these steps, only the scaling coefficients at level $L$ are retained while the detail coefficients are discarded. Due to the lifting step, ghost points are required for a block to coarsen when using wavelets with $\tilde{N} > 0$, with the precise number specified in the first column of \tref{table_wavelet_npoints}. In the example of \fref{fig:2.2_coarsen_2_a} for \wave{2}{2}, the green region needs one ghost point at the back ($\lambda_4^1$, needed to compute $\lambda_1^0$) and the blue region needs two ghost points at the front ($\lambda_8^1$ and $\lambda_9^1$, required to compute $\lambda_5^0$).

Discarding the detail coefficients on level $L$ does not affect the scaling coefficients on that level, however when $\tilde{N} > 0$ discarding these details does affect the scaling coefficients of adjacent blocks whose resolution has not changed. This can be seen and accounted for by performing an inverse wavelet transform from level $L$ back to level $L+1$ after discarding the details, and \revtwo{updating} the values of the affected scaling coefficients:
\begin{eqd}
\lambda_{j}^{L+1} \leftarrow \lambda_j^{L+1} - G^s_{j,n} \gamma^L_n 
\end{eqd}
where $\gamma^L_n$ are all the detail coefficients that we have discarded. To perform the update step on a fine block whose neighbor has coarsened, the fine block needs to have enough ghosts points to compute the values of $\gamma^L_n$ that are discarded, which increases significantly the ghost point requirements of the update step as shown in \tref{table_wavelet_npoints}. The distance (in index space) to the farthest detail to be discarded depends on the wavelet order and $\tilde{N}$, and is shown in \tref{table_wavelet_npoints} under the column \textit{coarse region extension}. 

In the specific example of \fref{fig:2.2_coarsen_2}, though the middle red region does not change resolution, we must still remove all the information associated with the discarded detail coefficients of the coarsened green and blue regions. For \wave{2}{2}, this corresponds to discarding the information associated with $\gamma_{1}^0$ on the left and $\gamma_{5}^0$ on the right. To achieve this, starting from the original uniform grid on the top line of \fref{fig:2.2_coarsen_2_a}, the red region requires two ghost points in front ($\lambda_2^1$ and $\lambda_3^1$, needed to compute $\gamma_1^0$) and three ghost points in the back ($\lambda_{10}^1$, $\lambda_{11}^1$, and $\lambda_{12}^1$, needed to compute $\gamma_5^0$), in order to perform the update step associated with the coarsening of both its neighboring grid regions. Lastly, as indicated in \fref{fig:2.2_coarsen_2}, for \wave{2}{2} we additionally discard the detail coefficient $\gamma_2^0$ when the left region coarsens, and $\gamma_4^0$ when the right region coarsens, which will be explained further in the next two sections below.

\subsubsection{Ghost point reconstruction}
\label{sec_ghost_reconstruction}
For any grid configuration with blocks at multiple levels of resolution, we have to be able to compute ghost points for each block at their local resolution level. We choose here to rely on the wavelets to do so for all ghosting operations, in order to be consistent with the grid adaptation operations. 

\begin{figure}[h!]
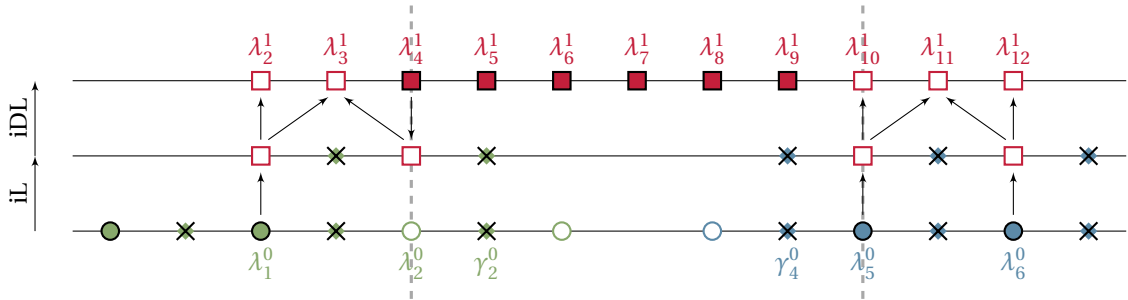

    \centering
    \InputIfFileExists{figures/w2.2_1d_ghost_inverse}{}{}
    \caption{Sketch of the process used to reconstruct ghost points (open symbols) for a 1D fine-resolution grid region (in red) surrounded by coarse grid regions on the front (green) and back (blue), using \wave{2}{2}. The arrows denote the flow of information from grid values to the unknown ghost points. After the ghost points for the fine region have been computed, the ghost points for both coarse regions can be found from a standard coarsening procedure (\fref{fig:2.2_coarsen_2_a})}
    \label{fig:w2.2_1d_ghost_inverse}
 \end{figure}

\Fref{fig:w2.2_1d_ghost_inverse} shows the computation of the ghost points for a fine region (in red) surrounded by neighboring coarse regions (in green and blue) for \wave{2}{2} in 1D. Ghost points to be computed are shown with open circles (on the coarse level) and open squares (on the fine level), whereas known scaling coefficients are shown in colored symbols with black outlines. A naive wavelet transform indicates the immediate problem that the ghost points for the finer region and those for the coarser region are interdependent: for instance, to compute $\lambda_3^1$ we would need to apply inverse lifting on $\lambda_2^0$, but $\lambda_2^0$ is in turn dependent on $\lambda_3^1$ through the dual lifting. This interdependency gets more intricate in higher dimensions and for higher-order wavelets.  Though these systems can be solved as proposed in \cite{Rossinelli:2015} for non-lifting wavelets, the associated implementation requires expansive look-up tables that significantly increase the memory footprint of the solver, especially in 3D \cite{Rees:2014}.
To avoid this complication, we choose to discard the fine-region detail coefficients that cause the interdependency between the wavelet transforms on the two levels, which we denote the `coarse-extension assumption'. In the case illustrated in \fref{fig:w2.2_1d_ghost_inverse} specifically, this means we discard the detail coefficients $\gamma_2^0$ and $\gamma_4^0$.
On a more abstract level, this choice  means that we effectively extend the coarse-level region across the resolution jump into the first few grid points on the neighboring fine resolution block.
To consistently follow-through with this assumption requires additional steps in our implementation that we will discuss more below. 
For now, with this assumption in place the ghost reconstruction across a resolution jump can be done in two steps. First, we can use the wavelet transform to compute ghost values for the fine resolution block through what is essentially a local refinement of the coarse grid, which can now be done explicitly.%

In \fref{fig:w2.2_1d_ghost_inverse} we illustrate this approach for the computation of the fine ghost values with the \wave{2}{2}. Focusing first on the right resolution jump, our assumption of explicitly discarding $\gamma_4^0$ formally sets the ghost point $\lambda^{1}_{10} = \lambda_5^0$ and further allows us to directly evaluate $\lambda^{1}_{11}$ as $\lambda^{1}_{11} = 1/2 (\lambda_5^0 + \lambda_6^0)$. This is essentially a refinement procedure where we take the coarse-level scaling coefficients and refine them under the assumption that all unknown detail coefficients involved in this process are zero, irrespective of the side of the interface where they exist. 

For the left resolution jump, the specific grid layout enables us to proceed in one of two ways. The first way is conceptually similar to how we describe the right interface, where by assumption $\gamma_2^0 = 0$ even though this detail resides in the fine resolution grid. Under this assumption, $\lambda_2^0 = \lambda_4^{1}$ and the ghost point $\lambda_3^{1} = 1/2 \left( \lambda_1^0 + \lambda_2^0 \right)$ is obtained from the inverse wavelet transform of the coarse-level data.
The second way relies on the inverse of the dual lifting step: we can directly write $\lambda_3^{1} = 1/2 \left( \lambda_4^{1} + \lambda_1^0 \right)$ without explicitly considering $\gamma_2^0$.
The two approaches are identical if $\gamma_2^0 = 0$. 
Note that this is only possible on the left interface; on the right interface, the value of $\lambda_{10}^{1}$ is not readily available unless $\gamma_4^0 = 0$, since  $\lambda_{10}^{1}$ belongs to the coarse region. In practice, we use the first approach and explicitly assume that both $\gamma_2^0$ and $\gamma_4^0$ are zero. We do so because the implementation of the inverse dual lifting in multiple dimension is not trivial (see \sect{sect_multiple_dim}), which would complicate the second approach.%

This concludes the computation of the ghost points for the fine level; afterwards, we can treat the region of the fine resolution as a local uniform grid that we coarsen in order to obtain the ghost points for the coarse grid levels. This procedure poses no further difficulties and is identical to the coarsening described above.

\subsubsection{Refinement}
\label{sec:refinement}
\begin{figure}[h!]
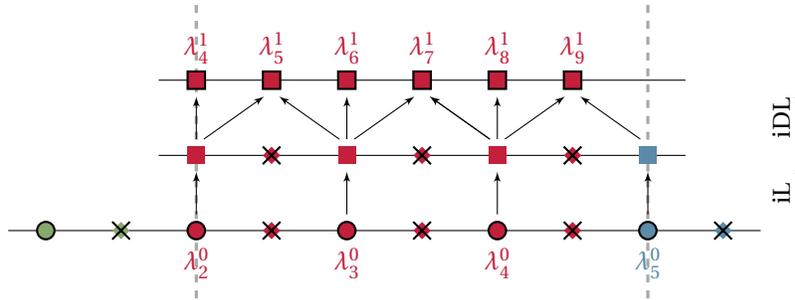

    \centering
    \InputIfFileExists{figures/w2.2_1d_refinement}{}{}
    \caption{Sketch of the refinement procedure of a 1D coarse grid region (in red) using \wave{2}{2}. Starting from the coarse-level scaling coefficients (bottom line), we use inverse lifting (iL) and inverse dual lifting (iDL) to compute the fine-level scaling coefficients (top line). For this wavelet, only one ghost point on the back of the refined region is required, indicated in blue.}
    \label{fig:w2.2_1d_refinement}
 \end{figure}

The refinement operation of a block away from resolution boundaries is trivially done through the subsequent application of the lifting and dual lifting filters. Near resolution boundaries, we retain the coarse-extension assumption introduced for the computation of ghost points described in the previous subsection, which enables the explicit computation of the fine-level scaling coefficients. This process is illustrated in \fref{fig:w2.2_1d_refinement} for the special case of the \wave{2}{2}, in which case only the ghost point $\lambda_{5}^0$ is required on the right resolution jump to compute the new scaling coefficient $\lambda_9^1$. The number of ghost points needed for a block to refine for any other wavelet considered here is shown in the last column of \tref{table_wavelet_npoints}.

\subsubsection{Substitution}
\label{sec_mr_substitution}
In the previous subsections we motivated and detailed the coarse-extension assumption, where we neglect specific fine-level detail coefficients near coarse-fine resolution jumps to facilitate explicit ghost reconstruction and refinement operations. We explained that this is essentially equivalent to an extension of the coarse-level region into a small band of the adjacent fine-level block. In practice however, these specific detail coefficients on the fine-level block may not be zero due to field operations on the associated fine-level scaling coefficient, such as during the evolution of a PDE. Without addressing this, we would inconsistently neglect high-frequency information during the grid adaptation and ghost reconstruction due to our coarse-extension assumption.

\begin{figure}[h!]
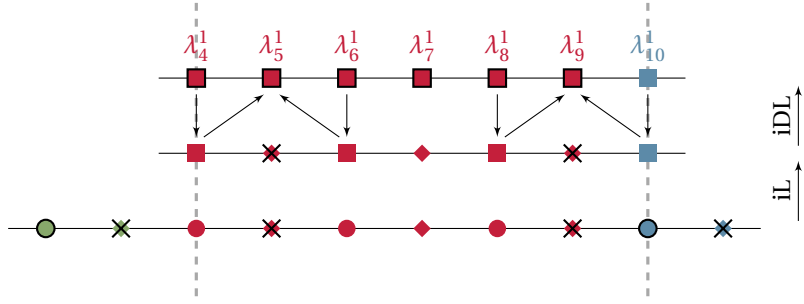

    \centering
    \InputIfFileExists{figures/w2.2_1d_substitution}{}{}
    \caption{Sketch of the substitution process, where spurious information associated with neglected detail coefficients in fine-resolution grid regions adjacent to coarse-resolution regions is explicitly discarded. For \wave{2}{2} in $1$D, as shown here, substitution is equivalent to the inverse dual lifting. }
    \label{fig:w2.2_1d_substitution}
 \end{figure}

To avoid that this spurious information persists and leads to an inconsistent wavelet transform on the two sides of the interface, we perform an additional ``substitution'' step where we overwrite each fine-level scaling coefficient associated with a neglected detail coefficient locally with wavelet-reconstructed values that will enforce a zero detail coefficient. To achieve this step, we use the dual lifting relationship
\begin{eqc}
\gamma_m^L = 0  \quad \Rightarrow \quad G^a_{m,j} \; \lambda_j^{L+1} \;=\; G^a_{m,{j \neq \left(2m+1\right)}} \; \lambda_{j\neq \left(2m+1\right)}^{L+1}  +  \lambda_{2m+1}^{L+1} = 0 \quad \Rightarrow \quad  \lambda_{2m+1}^{L+1} \leftarrow - G^a_{m,{j\neq \left(2m+1\right) }}  \; \lambda_{j\neq \left(2m+1\right)}^{L+1} 
\end{eqc}

where we used the fact that $G^a_{m,0} = 1$ for all wavelets considered in this work. We note that in $1$D, the proposed approach is exactly the inverse of the dual lifting step, as illustrated in \fref{fig:w2.2_1d_substitution}. This substitution step is subtle so we point out that this step does not affect the order of accuracy of the wavelet operations, only removes detail coefficient values that are generated during the PDE evolution starting from values below the coarsening threshold, and can be rigorously understood as the consistent enforcement of the coarse-extension assumption. In practice, we apply the substitution step as part of the ghost point reconstruction, immediately after the computation of the fine ghost points and before computation of the coarser ghost points.

\subsubsection{Block adaptation criteria}
\label{subsubsec:adaptation}
\Sect{sec:wavelet_compression} describes how to compress a signal in 1D given its wavelet transformation. Here we explain how we transform this condition to compress data on a block-structured grid, detect emerging scales and \revtwo{manage} the need to refine during a simulation.

\paragraph{Compression}
Starting with the former, for each block at level $L$ we compute all associated details $\gamma^{L-1}$ through the forward wavelet transform, and take the maximum value for each block $b$ as $\normi{\gamma^{L-1}}^b$.
Consistent with our coarse-extension assumption above, we consider within this local infinite norm also the set of details that we require to be negligible when computing ghost points, the refinement relation, and the coarsening steps, even though these detail coefficients might physically reside in adjacent blocks. The \revtwo{number of} additional details considered beyond the block boundary is given in the column \textit{coarse region extension} in \tref{table_wavelet_npoints}.
Once we have the maximum detail coefficient on each block, we decide on an action to take. Due to the octree nature of our grid, we can only coarsen all leaf blocks within a single tree node simultaneously. Therefore, for each set of $2^d$ leaf blocks in the grid where $d$ is the spatial dimension, we reduce them to a single coarser level block if each of the leaf blocks $b$ satisfies $\normi{\gamma^{L-1}}^b < \epsilon_c$, with $\epsilon_c$ the coarsening threshold. This criterion constitutes a generalization of \eqqref{eq_multires_epsilon} \revtwo{to a block-structured grid. This approach implies that the compression rate of a given signal decreases as the block size increases, due to a reduced granularity in the grid.} \revtwo{We emphasize here} that \revtwo{by} including the details neglected during the coarse-extension assumption \revtwo{within} our definition of $\normi{\gamma^{L-1}}^b$, \revtwo{the compression approach remains} consistent with \eqqref{eq_multires_epsilon}. 

\paragraph{Refinement}
The criterion to refine is necessarily more ad-hoc, as it aims to predict where new scales are expected during the evolution of the equations based on an instantaneous analysis of the field. Different approaches exist \revtwo{within existing wavelet-based adaptive grid methods}, such as increasing the resolution of neighboring blocks to take into account the smallest scales created by the PDE \citep{Kevlahan:2005}. %
Here we follow \citep{Rossinelli:2015} to rely on a user-specified tolerance $\epsilon_r > \epsilon_c$ that determines whether refinement is necessary. Using this approach, a block $b$ is refined if $\normi{\gamma^{L-1}}^b > \epsilon_r$, i.e.\ if the detail coefficients of the current information exceed a user defined threshold. We will analyze this choice in the validation and result sections below.

\revtwo{Under the above compression policy,} we are guaranteed by the wavelet framework to discard only information encoded by detail coefficients that do not exceed $\epsilon_c$. \revtwo{With the refinement approach,} the maximum detail coefficients during the evolution of the equations are guaranteed to never exceed $\epsilon_r$, since we would refine when that happens. \revtwo{Specifically,} if we refine a block at level $L$ for which $\normi{\gamma^{L-1}}^b > \epsilon_r$, we create new blocks at level $L+1$, each one of which is characterized by $\normi{\gamma^{L}}^b = 0$. This means we have to make sure we do not coarsen blocks just after they have been refined, even though technically their detail coefficients are smaller than the coarsening threshold. \revtwo{The implementation of this requirement is discussed in \cref{sect_implementation_global}.}

\paragraph{Ratio between compression and refinement thresholds}
\revtwo{From \eqqref{eq:interp_error_bound}} we know that $\normi{\gamma}^b \propto h^N$ with $N$ the interpolation order of the wavelet. \revtwo{This means} that coarsening a block will \revtwo{generally} increase its detail \revtwo{coefficient} by a factor $2^N$. \revtwo{Consequently,} if \revtwo{we choose} $\epsilon_r / \epsilon_c < 2^N$ and a block with $\normi{\gamma}^b$ only slightly below the coarsening threshold is coarsened, its details will exceed $\epsilon_r$ after coarsening. \revtwo{In this case, the block} will be flagged for refinement again, leading to flip-flops in the grid adaptation. 
Conversely, choosing $\epsilon_r / \epsilon_c > 2^N$ can lead to a non unique grid: if we consider a block with $\normi{\gamma}^b$ slightly above the coarsening threshold and therefore admissible on the grid, the same block coarsened by one level would also be admissible on the grid. In this case the adaptation is therefore not unique and the obtained grid depends on external factors such as the initial level. 

In practice, we observe $\epsilon_r$ to be the threshold that determines the overall accuracy of the simulation, since this is the threshold that sets the maximum value of detail coefficient admissible on the grid. Then $\epsilon_c$ determines the compression rate, or how much information we are willing to discard for a given $\epsilon_r$. This can be controlled by the ratio $\epsilon_r/\epsilon_c$, which we make sure to set to $\epsilon_r / \epsilon_c > 2^N$ to prevent the flip-flopping described above. The effect of both $\epsilon_r$ and the ratio $\epsilon_r/\epsilon_c$ is discussed through numerical experiments below. 

\paragraph{Adaptation frequency}
\revone{For transient problems we have to choose a frequency of adaptation that balances the need to adjust the grid to dynamically evolving scales against the computational cost of changing the grid. In general, lowering the mesh adaptation frequency will lead to wavelet detail coefficients in the grid that fall more and more below the compression threshold in some regions, and increase more and more beyond the refinement threshold in others. The former will not increase the errors made during the simulation, but the latter could potentially lead to the grid not capturing emerging or transported scales that are relevant to the PDE evolution. A counter point to this is that the block-structured grid contains significant `inertia', which grows with the block size, due to the fact that we refine even when one single detail in a block exceeds the threshold.}

\revone{To put the adaptation frequency into context, we can consider a transport problem with characteristic velocity $U$. In this case, the time it takes for the solution to travel $\alpha$ grid points at level $L$ is $2^{-L} \alpha / U$ (assuming a unit cube as root domain). If we adapt every $N_a$ timesteps and follow a CFL-based timestep constraint, the solution has traveled $\cfl  N_a$ grid points between successive adaptations. Though in practice likely problem-dependent our results in \cref{sec:translation_blob} show that for a smooth signal and suitably small refinement threshold this number is allowed to be of comparable value as the block size $N_b$ leading to $N_a \sim N_b/(\cfl)$.}

\begin{table}[ht!]
\centering
\begin{tabular}{|rr|cc|cc|cc|cc|}
\hline
\multicolumn{2}{|c|}{} & \multicolumn{2}{|c|}{\textbf{coarsening}} & \multicolumn{4}{|c|}{\textbf{update after neighbor has coarsened}} & \multicolumn{2}{|c|}{\textbf{refinement}}  \\
\hline\hline
\multicolumn{2}{|c|}{\textbf{wavelet order}} & \multicolumn{2}{|c|}{\textit{\# ghost points}} & \multicolumn{2}{|c|}{\textit{coarse region extension}} & \multicolumn{2}{|c|}{\textit{\# ghost points}} & \multicolumn{2}{|c|}{\textit{\# ghost points}}  \\
 \hline
$N$	& $\tilde{N}$ & \textit{front} & \textit{back}  & \textit{front} ($1$D / $3$D)	& \textit{back} ($1$D / $3$D)  & \textit{front} ($1$D / $3$D)	& \textit{back} ($1$D / $3$D)  & \textit{front} 	& \textit{back}  \\
 \hline
$2$	& $0$	& $0$ 	& $0$ 	& $0$ / $0$ 	& $0$ / $1$ 	& $0$ / $0$ 	& $0$ / $1$	& $0$ 	& $1$ 	\\
$2$	& $2$	& $2$ 	& $1$	& $1$ 	& $2$  	& $2$  	& $3$ 	& $0$ 	& $1$ 	\\
$4$	& $0$	& $0$ 	& $0$	& $0$ / $2$ 	& $0$ / $3$ 	& $0$ / $4$ 	& $0$ / $5$	& $1$ 	& $2$ 	\\
$4$	& $2$	& $4$ 	& $3$	& $3$	& $4$  	& $6$  	& $7$ 	& $1$ 	& $2$ 	\\
$6$	& $0$	& $0$ 	& $0$	& $0$ / $4$	& $0$ / $5$ 	& $0$ / $8$ 	& $0$ / $9$	& $2$ 	& $3$ 	\\
$6$	& $2$	& $6$ 	& $5$	& $5$ 	& $6$  	& $10$  	& $11$ 	& $2$ 	& $3$ 	\\
\hline
\end{tabular}
\caption{Ghost points requirement for a block when coarsening (see \fref{fig:2.2_coarsen_2_a}); for a block at level $L+1$ performing the update step when a neighbor at level $L$ has coarsened (see \fref{fig:2.2_coarsen_2_b}); and when refining a block (see \fref{fig:w2.2_1d_refinement}).
The \textit{coarse region extension} indicates how many additional scaling and detail coefficients at level $L$ have to be computed outside of the blocks boundaries to perform the update step. The associated ghost points are required at the level $L+1$ of the block whose neighbor has coarsened. The size of these regions can vary between $1$D and $3$D, as explained in \sect{sect_multiple_dim}.
}
\label{table_wavelet_npoints}
\end{table}

\subsection{Extension to multiple dimensions}
\label{sect_multiple_dim}

Here we detail the extension of the above methodology to multiple spatial dimensions, starting with the coarsening operation and subsequently emphasizing some of the implementation details to consider.

The filter application in $3$D relies on the successive application of the corresponding $1$D filters in each dimension. To clarify the notation we use a superscript on all the filters to denote the direction in which the filter is applied.
For coarsening we then obtain %
\be
\lambda^L_{k_x,k_y,k_z} = \left[ H^a_{X} \times H^a_{Y} \times H^a_{Z} \right] \lambda^{L+1} %
\eec
which exclusively relies on the $H^a$ filter applied tensorially on the scaling coefficients at level $L+1$. 

To compute the detail values we alternatively apply the filters $H^a$ or $G^a$ depending on the scaling or detail behavior in the considered dimension, as dictated by whether the associated index in that dimension is even or odd. This means we can distinguish different ``degrees'' of detail coefficients, given by the number of directions in which the coefficient behaves as a detail information. 
Specifically, the first degree detail coefficients, which we collectively denote as $\gamma^{L}_{1^{\circ}}$, have an odd index in one direction only and are given by
\begin{eqd}
\gamma^{L}_{m_x} = \left[ G^a_{X} \times H^a_{Y} \times H^a_{Z} \right] \; \lambda^{L+1} \qquad
\gamma^{L}_{m_y} = \left[ H^a_{X} \times G^a_{Y} \times H^a_{Z} \right] \; \lambda^{L+1} \qquad
\gamma^{L}_{m_z} = \left[ H^a_{X} \times H^a_{Y} \times G^a_{Z} \right] \; \lambda^{L+1} \qquad
\end{eqd}
Similarly the second degree detail coefficients, which we collectively denote as $\gamma^{L}_{2^{\circ}}$, have two odd indices and we obtain
\begin{eqd}
\gamma^{L}_{m_x,m_y} = \left[ G^a_{X} \times G^a_{Y} \times H^a_{Z} \right] \; \lambda^{L+1} \qquad
\gamma^{L}_{m_x,m_z} = \left[ G^a_{X} \times H^a_{Y} \times G^a_{Z} \right] \; \lambda^{L+1} \qquad
\gamma^{L}_{m_y,m_z} = \left[ H^a_{X} \times G^a_{Y} \times G^a_{Z} \right] \; \lambda^{L+1}
\end{eqd}
Finally, the third degree scaling coefficient, which we collectively denote as $\gamma^{L}_{3^{\circ}}$, are obtained as
\begin{eqd}
\gamma^{L}_{m_x,m_y,m_z} = \left[ G^a_{X} \times G^a_{Y} \times G^a_{Z} \right] \; \lambda^{L+1}
\end{eqd}
In order to relax further the notation, for all detail coefficients $\gamma^L = \{\gamma^{L}_{1^{\circ}}, \gamma^{L}_{2^{\circ}}, \gamma^{L}_{3^{\circ}}\}$ we will refer to the fine scaling coefficients located at the same positions as \revtwo{$\{\lambda^{L+1}_{1^{\circ}}, \lambda^{L+1}_{2^{\circ}}, \lambda^{L+1}_{3^{\circ}}\}$}, respectively. We also have the `zeroth' degree scaling coefficients associated with even indices in all three directions, which we denote by $\lambda^{L+1}_{0^{\circ}}$ so that $\lambda^{L+1} = \{\lambda^{L+1}_{0^{\circ}}, \lambda^{L+1}_{1^{\circ}},\lambda^{L+1}_{2^{\circ}},\lambda^{L+1}_{3^{\circ}}\}$. %

Revisiting the coarse-extension assumption we made in $1$D, its extension to three dimensions can be formulated as discarding any first, second and third degree detail coefficients in a fine block adjacent to a resolution jump that are involved in the multi-dimensional refinement scheme. This is reflected by the higher $3$D values in the column \text{coarse region extension} of \tref{table_wavelet_npoints}, which represents the index-space distance from the block boundary to the farthest detail coefficient that is discarded. As before, we discard this information through the substitution procedure, which inverts locally the inverse dual lifting step. In $3$D, the substitution on first degree scaling coefficients is performed as
\begin{eqc}
\gamma_{1^{\circ}}^{L} = \left[ G^a \right]_{(i)\neq (0)}  \; \lambda^{L+1} + \lambda_{1^\circ}^{L+1} \quad \Rightarrow  \quad\lambda_{1^{\circ}}^{L+1} \leftarrow \left[ G^a \right]_{(i)\neq (0)}  \; \lambda^{L+1}
\end{eqc}
where we used that $\left[ G^a \right]_{(i) =(0)} = 1$ for all wavelets considered in this work. For second degree scaling coefficients, this becomes
\begin{eqc}
\gamma_{2^{\circ}}^{L} = \left[ G^a \times G^a \right]_{(i,j)\neq (0,0)}  \; \lambda^{L+1} + \lambda_{2^\circ}^{L+1} \quad \Rightarrow \quad  \lambda_{2^{\circ}}^{L+1} \leftarrow \left[ G^a \times G^a \right]_{(i,j)\neq (0,0)}  \; \lambda^{L+1}
\end{eqc}
and for third degree, we find
\begin{eqd}
\gamma_{3^{\circ}}^{L} = \left[ G^a \times G^a \times G^a \right]_{(i,j,k)\neq (0,0,0)}  \; \lambda^{L+1} + \lambda_{3^\circ}^{L+1} \quad \Rightarrow \quad  \lambda_{3^{\circ}}^{L+1} \leftarrow \left[ G^a \times G^a \times G^a\right]_{(i,j,k)\neq (0,0,0)}  \; \lambda^{L+1}
\end{eqd}
At first it seems that each equation is interdependent since the first, second, and third degree scaling coefficients are all included in $\lambda^{L+1}$. However, analyzing the filter $G^a$ reveals that the first degree scaling coefficients $\lambda_{1^{\circ}}^{L+1}$ only need $\lambda_{0^{\circ}}^{L+1}$ to be updated, the second degree coefficients $\lambda_{2^{\circ}}^{L+1}$ need $\lambda_{0^{\circ}}^{L+1}$ and $\lambda_{1^{\circ}}^{L+1}$, \revtwo{and so forth}. This means we can consistently and explicitly perform the substitution step in $3$D by first updating the first degree, then the second degree, and finally the third degree scaling coefficients according to the values of the respective neglected detail coefficients.

\subsection{Boundary conditions}
All theory described above is for infinite signals and grids, and the validation and results cases below rely on fields with compact support, so that boundary conditions are not relevant. In practice, we do need an implementation of boundary conditions for finite signals, which we created by relying on interpolation and extrapolation at the domain boundary. Our solver currently supports zero value (Dirichlet), zero flux (Neumann), plain zeros filling and extrapolation boundary conditions. They all rely on polynomial interpolation done at the order of the wavelet used, which is compatible with the non-lifted wavelet theory but doesn't conserve the moments when used with lifted wavelets. To improve the numerical properties of the interpolation we \revtwo{use Neville's} algorithm \citep{Press:2007}. Wavelets on the interval \citep{Donoho:1992a} would provide more consistent implementations of such boundary conditions, and the lifting scheme provides avenues for moment conservation \citep{Fernandez:1996}, but we reserve this for future work.

%% file: figures/fwt_iwt_lifting_tikz.tex
\tikzstyle{filter} = [draw,minimum size=2em]
\tikzstyle{scaling} = [draw,minimum size=2em]
\tikzstyle{detail} = [draw,minimum size=2em]
\tikzstyle{branch}=[draw,fill=white,shape=circle,inner sep=2pt]

\begin{tikzpicture}[>=latex']

    \node[scaling] at (-0.1,0) (scaling_jp1) {$\lambda^{L+1}$};
    \node at (1,0) (separation) {split};
    \def\x{1}
    \node[coordinate] at (\x,+1) (block_Ht) {};
    \node[coordinate] at (\x,-1) (block_Gt) {};
    \draw[-] (scaling_jp1) -- (separation);
    \draw[-] (separation) |- (block_Ht);
    \draw[-] (separation) |- (block_Gt);
    
    \def\x{2}
    \node[branch] at (\x,-1) (branch_dbf) {$+$};
    \coordinate (branch_duf) at (\x,+1);
    \node[filter] at (\x,0) (block_St) {$\tilde{S}$};
    
    \def\x{3}
    \node[branch] at (\x,1) (branch_uf) {$+$};
    \coordinate (branch_bf) at (\x,-1);
    \node[filter] at (\x,0) (block_S) {$S$};

    \def\x{4}
    \node[scaling,] at (\x,+1) (scaling_j) {$\lambda^L$};
    \node[detail] at (\x,-1) (detail_j) {$\gamma^L$};
    
    \draw[->] (block_Gt) -- (branch_dbf) -- (branch_bf)  -- (detail_j);
    \draw[->] (block_Ht) -- (branch_duf) -- (branch_uf) -- (scaling_j);
    \draw[->] (branch_bf) -- (block_S);
    \draw[->] (block_S) -- (branch_uf);
    \draw[->] (branch_duf) -- (block_St);
    \draw[->] (block_St) -- (branch_dbf);
    
    \def\x{5}
    \node[branch] at (\x,1) (branch_ub) {$-$};
    \coordinate (branch_bb) at (\x,-1);
    \node[filter] at (\x,0) (block_Sb) {$S$};
    
    \def\x{6}
    \node[branch] at (\x,-1) (branch_dbb) {$-$};
    \coordinate (branch_dub) at (\x,+1);
    \node[filter] at (\x,0) (block_Stb) {$\tilde{S}$};

    \def\x{7}
    \node[coordinate] at (\x,+1) (block_H) {};
    \node[coordinate] at (\x,-1) (block_G) {};

    \draw[-] (scaling_j) -- (branch_ub)  -- (branch_dub) -- (block_H);
    \draw[-] (detail_j) -- (branch_bb) -- (branch_dbb) -- (block_G);
    \draw[->] (branch_dub) -- (block_Stb);
    \draw[->] (block_Stb) -- (branch_dbb);
    \draw[->] (branch_bb) -- (block_Sb);
    \draw[->] (block_Sb) -- (branch_ub);

    \def\x{7}
    \node at (\x,0) (merge) {merge};
    \draw[-] (block_H.east) -| (merge);
    \draw[-] (block_G.east) -| (merge);
    
    \def\x{8.2}
    \node[scaling] at (\x,0) (rescaling_jp1) {$\lambda^{L+1}$};
    \draw[-] (merge) -- (rescaling_jp1);
    
    \draw[->] ($ (block_Ht) + (-1,0.75) $) --  node[above] {Forward} ($ (block_Ht) + (2.5,0.75) $);
    \draw[->] ($ (block_H) + (-2.5,0.75) $) --  node[above] {Inverse} ($ (block_H) + (1,0.75) $);

\end{tikzpicture}

%% file: figures/ftw_iwt_tikz.tex
\tikzstyle{filter} = [draw,minimum size=2em]
\tikzstyle{scaling} = [draw,minimum size=2em]
\tikzstyle{detail} = [draw,minimum size=2em]

\begin{tikzpicture}[>=latex']
    \node[scaling] at (0,0) (scaling_jp1) {$\lambda^{L+1}$};
    \coordinate [right of=scaling_jp1] (separation);
    \node[filter] at (1.8,+1) (block_Ht) {${H}^a$};
    \node[filter] at (1.8,-1) (block_Gt) {${G}^a$};
    \draw[-] (scaling_jp1) -- (separation);
    \draw[->] (separation) |- (block_Ht.west);
    \draw[->] (separation) |- (block_Gt.west);
    
    \node[scaling,] at (3,+1) (scaling_j) {$\lambda^L$};
    \node[detail] at (3,-1) (detail_j) {$\gamma^L$};
    \draw[->] (block_Ht) -- (scaling_j);
    \draw[->] (block_Gt) -- (detail_j);
    
    \node[filter] at (4.2,+1) (block_H) {${H}^s$};
    \node[filter] at (4.2,-1) (block_G) {${G}^s$};
    \draw[->] (scaling_j) -- (block_H);
    \draw[->] (detail_j) -- (block_G);
    \coordinate (merge) at (5,0);
    \draw[-] (block_H.east) -| (merge);
    \draw[-] (block_G.east) -| (merge);
    
    \node[scaling] at (6,0) (rescaling_jp1) {$\lambda^{L+1}$};
    \draw[->] (merge) -- (rescaling_jp1);
    
    \draw[->] ($ (block_Ht) + (-2,0.75) $) --  node[above] {Forward} ($ (block_Ht) + (0.5,0.75) $);
    \draw[->] ($ (block_H) + (-0.5,0.75) $) --  node[above] {Inverse} ($ (block_H) + (2,0.75) $);

\end{tikzpicture}

%% file: figures/w2.2_1d_coarsen.tex
\tikzset{cross/.style={cross out, draw=black, minimum size=0.5em, inner sep=0pt, outer sep=0pt},cross/.default={10pt}}

\tikzstyle{scaling} = [fill=black,shape=rectangle,minimum size=0.75em,inner sep=0pt,outer sep=0pt,thick]
\tikzstyle{detail} = [fill=black,diamond,very thick,minimum size=0.75em,inner sep=0pt,outer sep=0pt]
\tikzstyle{bullet}=[draw,fill=black,shape=triangle,minimum size=0.5em]
\tikzstyle{temp} = [fill=\figgrey,shape=circle,minimum size=0.5em]

\def\colblockzero{asparagus}
\def\colblockone{cardinal}
\def\colblocktwo{airforceblue}

\def\opacity{0.5}

\begin{tikzpicture}[>=latex']
	\def\shift{0.1}
	\def\step{1};
	\def\l{0}
	\def\ll{-1}
	\def\lll{-2}
	\def\llll{-3}

   	\draw[-] (-0.5,\l) -- (13*\step+0.5,\l);
	\draw[-] (-0.5,\ll) -- (13*\step+0.5,\ll);
	\draw[-] (-0.5,\lll) -- (13*\step+0.5,\lll);
    	
	\draw[-,very thick,dashed,color=\figgrey] (4*\step,1) -- (4*\step,-3);
	\draw[-,very thick,dashed,color=\figgrey] (10*\step,1) -- (10*\step,-3);
	
	\draw[->] (-1,\l-\shift) -- node[above,rotate=90] {DL} (-1,\ll+\shift);
	\draw[->] (-1,\ll-\shift) -- node[above,rotate=90] {L} (-1,\lll+\shift);
	
   	 \foreach \pid [evaluate=\pid as \x using {int(\pid*\step)}] in {0,...,6} {
		\def\y{\l}
		\ifthenelse{\x<2}{
			\foreach \id [
				evaluate=\id as \xe using {int(\x*2 + \id)},
				evaluate=\id as \xid using {int(\x*2  + \id)}] in {0,1}{
				\node[scaling,fill=\colblockzero,draw=black,label=above:{\textcolor{\colblockzero}{$\lambda^{1}_{\xid}$}}] at ({\x*2 + \id},\y) (scal_\xe) {};
			};
       		}{
        			\ifthenelse{\x<5}{
				\foreach \id [
					evaluate=\id as \xe using {int(\x*2 + \id)},
					evaluate=\id as \xid using {int(\x*2 + \id)}
					] in {0,1}{
					\node[scaling,fill=\colblockone,draw=black,label=above:{\textcolor{\colblockone}{$\lambda^{1}_{\xid}$}}] at ({\x*2 + \id},\y) (scal_\xe) {};
				};

        			}{
            			\foreach \id [
					evaluate=\id as \xe using {int(\x*2 + \id)},
					evaluate=\id as \xid using {int(\x*2 + \id)}
					] in {0,1}{
					\node[scaling,fill=\colblocktwo,draw=black,label=above:{\textcolor{\colblocktwo}{$\lambda^{1}_{\xid}$}}] at ({\x*2 + \id},\y) (scal_\xe) {};
				};
        			};
        		};
		\def\y{\ll}
		\ifthenelse{\x<2}{
			\foreach \id [
				evaluate=\id as \xe using {int(\x*2)},
				evaluate=\id as \xo using {int(\x*2 + 1)},
				evaluate=\id as \xid using {int(\x*2 + \id)}] in {0}{
				\node[scaling,fill=\colblockzero] at ({\x*2 + \id},\y) (dual_\xe) {};
				\node[detail,fill=\colblockzero] at ({\x*2 + \id + 1},\y) (dual_\xo) {};
				
				\draw[->,color=black] ($(scal_\xe.south) + (\shift,-\shift)$) -- ($(dual_\xo.north) - (\shift,-\shift)$) ;
				\draw[->,color=black] ($(scal_\xe.south) + (2*\step-\shift,-\shift)$) -- ($(dual_\xo.north) - (-\shift,-\shift)$) ;
				\draw[->,color=black] ($(scal_\xo.south) + (0,-\shift)$) -- ($(dual_\xo.north) - (0,-\shift)$) ;
				\draw[->,color=black] ($(scal_\xe.south) + (0,-\shift)$) -- ($(dual_\xe.north) - (0,-\shift)$) ;
			};
       		}{
        			\ifthenelse{\x<5}{
				\ifthenelse{\x=3}{
				\foreach \id [
					evaluate=\id as \xe using {int(\x*2)},
					evaluate=\id as \xo using {int(\x*2 + 1)},
					evaluate=\id as \xid using {int(\x*2 + \id)}] in {0}{
					\node[scaling,fill=\colblockone] at ({\x*2 + \id},\y) (dual_\xe){};
					\node[detail,fill=\colblockone] at ({\x*2 + \id + 1},\y) (dual_\xo) {};
					
					\draw[->,color=black] ($(scal_\xe.south) + (\shift,-\shift)$) -- ($(dual_\xo.north) - (\shift,-\shift)$) ;
					\draw[->,color=black] ($(scal_\xe.south) + (2*\step-\shift,-\shift)$) -- ($(dual_\xo.north) - (-\shift,-\shift)$) ;
					\draw[->,color=black] ($(scal_\xo.south) + (0,-\shift)$) -- ($(dual_\xo.north) - (0,-\shift)$) ;
					\draw[->,color=black] ($(scal_\xe.south) + (0,-\shift)$) -- ($(dual_\xe.north) - (0,-\shift)$) ;
				};}{
				\foreach \id [
					evaluate=\id as \xe using {int(\x*2)},
					evaluate=\id as \xo using {int(\x*2 + 1)},
					evaluate=\id as \xid using {int(\x*2 + \id)}] in {0}{
					\node[scaling,fill=\colblockone] at ({\x*2 + \id},\y) (dual_\xe) {};
					\node[detail,fill=\colblockone] at ({\x*2 + \id + 1},\y) (dual_\xo) {};
					
					\draw[->,color=black] ($(scal_\xe.south) + (\shift,-\shift)$) -- ($(dual_\xo.north) - (\shift,-\shift)$) ;
					\draw[->,color=black] ($(scal_\xe.south) + (2*\step-\shift,-\shift)$) -- ($(dual_\xo.north) - (-\shift,-\shift)$) ;
					\draw[->,color=black] ($(scal_\xo.south) + (0,-\shift)$) -- ($(dual_\xo.north) - (0,-\shift)$) ;
					\draw[->,color=black] ($(scal_\xe.south) + (0,-\shift)$) -- ($(dual_\xe.north) - (0,-\shift)$) ;
				};
				};
        			}{
            			\foreach \id [
					evaluate=\id as \xe using {int(\x*2)},
					evaluate=\id as \xo using {int(\x*2 + 1)},
					evaluate=\id as \xid using {int(\x*2 + \id)}
					] in {0}{
						\node[scaling,fill=\colblocktwo] at ({\x*2 + \id},\y) (dual_\xe) {};
						\node[detail,fill=\colblocktwo] at ({\x*2 + \id + 1},\y) (dual_\xo) {};
						
					\draw[->,color=black] ($(scal_\xe.south) + (\shift,-\shift)$) -- ($(dual_\xo.north) - (\shift,-\shift)$) ;
					\ifthenelse{\x<6}{
					\draw[->,color=black] ($(scal_\xe.south) + (2*\step-\shift,-\shift)$) -- ($(dual_\xo.north) - (-\shift,-\shift)$) ;
					}{};
					\draw[->,color=black] ($(scal_\xo.south) + (0,-\shift)$) -- ($(dual_\xo.north) - (0,-\shift)$) ;
					\draw[->,color=black,dashed] ($(scal_\xe.south) + (0,-\shift)$) -- ($(dual_\xe.north) - (0,-\shift)$) ;
				};
        			};
        		};
		\def\y{\lll}
		\ifthenelse{\x<2}{
			\foreach \id [
				evaluate=\id as \xe using {int(\x*2)},
				evaluate=\id as \xo using {int(\x*2 + 1)},
				evaluate=\id as \xid using {int(\x + \id)}] in {0}{
				\node[scaling,draw=black,shape=circle,fill=\colblockzero,label=below:{\textcolor{\colblockzero}{$\lambda^{0}_{\xid}$}}] at ({\x*2 + \id},\y) (coarse_\xe) {};
				\node[detail,fill=\colblockzero,draw opacity=\opacity,label=below:{\textcolor{\colblockzero}{$\gamma^{0}_{\xid}$}}] at ({\x*2 + \id + 1},\y) (coarse_\xo) {};
				\node[cross out,draw=black,thick] at ({\x*2 + \id + 1},\y){};

				\ifthenelse{\x>0}{
				\draw[->,color=black] ($(dual_\xe.south) + (\shift - \step,-\shift)$) -- ($(coarse_\xe.north) - (\shift,-\shift)$) ;
				}{};
				\draw[->,color=black] ($(dual_\xo.south) + (-\shift,-\shift)$) -- ($(coarse_\xe.north) - (-\shift,-\shift)$) ;
				\draw[->,color=black] ($(dual_\xe.south) + (0,-\shift)$) -- ($(coarse_\xe.north) - (0,-\shift)$) ;
				
				\draw[->,color=black] ($(dual_\xo.south) + (0,-\shift)$) -- ($(coarse_\xo.north) - (0,-\shift)$) ;
			};
       		}{
        			\ifthenelse{\x<5}{
				\ifthenelse{\x=2}{
				\foreach \id [
					evaluate=\id as \xe using {int(\x*2)},
					evaluate=\id as \xo using {int(\x*2 + 1)},
					evaluate=\id as \xid using {int(\x + \id)}] in {0}{
					\node[scaling,shape=circle,fill=\colblockone,label=below:{\textcolor{\colblockone}{$\lambda^{0}_{\xid}$}}] at ({\x*2 + \id},\y) (coarse_\xe) {};
					\node[detail,fill=\colblockone,draw opacity=\opacity,label=below:{\textcolor{\colblockone}{$\gamma^{0}_{\xid}$}}] at ({\x*2 + \id + 1},\y) (coarse_\xo) {};
					\node[cross out,draw=black,thick] at ({\x*2 + \id + 1},\y){};
					
					\draw[->,color=black] ($(dual_\xo.south) + (0,-\shift)$) -- ($(coarse_\xo.north) - (0,-\shift)$) ;
				\draw[->,color=black] ($(dual_\xe.south) + (\shift - \step,-\shift)$) -- ($(coarse_\xe.north) - (\shift,-\shift)$) ;
				\draw[->,color=black] ($(dual_\xo.south) + (-\shift,-\shift)$) -- ($(coarse_\xe.north) - (-\shift,-\shift)$) ;
				\draw[->,color=black] ($(dual_\xe.south) + (0,-\shift)$) -- ($(coarse_\xe.north) - (0,-\shift)$) ;
				};
				}{
				\foreach \id [
					evaluate=\id as \xe using {int(\x*2)},
					evaluate=\id as \xo using {int(\x*2 + 1)},
					evaluate=\id as \xid using {int(\x + \id)}] in {0}{
					\node[scaling,shape=circle,fill=\colblockone,label=below:{\textcolor{\colblockone}{$\lambda^{0}_{\xid}$}}] at ({\x*2 + \id},\y) (coarse_\xe) {};
					\node[detail,fill=\colblockone,label=below:{\textcolor{\colblockone}{$\gamma^{0}_{\xid}$}}] at ({\x*2 + \id + 1},\y) (coarse_\xo) {};
					\ifthenelse{\x=4}{
					\node[cross out,draw=black,thick] at ({\x*2 + \id + 1},\y){};
					}{};
					
					\draw[->,color=black] ($(dual_\xe.south) + (\shift - \step,-\shift)$) -- ($(coarse_\xe.north) - (\shift,-\shift)$) ;
					\draw[->,color=black] ($(dual_\xo.south) + (-\shift,-\shift)$) -- ($(coarse_\xe.north) - (-\shift,-\shift)$) ;
					\draw[->,color=black] ($(dual_\xe.south) + (0,-\shift)$) -- ($(coarse_\xe.north) - (0,-\shift)$) ;
					\draw[->,color=black] ($(dual_\xo.south) + (0,-\shift)$) -- ($(coarse_\xo.north) - (0,-\shift)$) ;
				};

				};
        			}{
            			\foreach \id [
					evaluate=\id as \xe using {int(\x*2)},
					evaluate=\id as \xo using {int(\x*2 + 1)},
					evaluate=\id as \xid using {int(\x  + \id)}
					] in {0}{
						\node[scaling,draw=black,shape=circle,fill=\colblocktwo,label=below:{\textcolor{\colblocktwo}{$\lambda^{0}_{\xid}$}}] at ({\x*2 + \id},\y) (coarse_\xe) {};
						\node[detail,fill=\colblocktwo,label=below:{\textcolor{\colblocktwo}{$\gamma^{0}_{\xid}$}}] at ({\x*2 + \id + 1},\y) (coarse_\xo) {};
						\node[cross out,draw=black,thick] at ({\x*2 + \id + 1},\y){};
						\draw[->,color=black] ($(dual_\xe.south) + (\shift - \step,-\shift)$) -- ($(coarse_\xe.north) - (\shift,-\shift)$) ;
						\draw[->,color=black] ($(dual_\xo.south) + (-\shift,-\shift)$) -- ($(coarse_\xe.north) - (-\shift,-\shift)$) ;
						\draw[->,color=black] ($(dual_\xe.south) + (0,-\shift)$) -- ($(coarse_\xe.north) - (0,-\shift)$) ;
						
						\draw[->,color=black] ($(dual_\xo.south) + (0,-\shift)$) -- ($(coarse_\xo.north) - (0,-\shift)$) ;
				};
        			};
        		};
    }
\end{tikzpicture}

%% file: figures/w2.2_1d_smoothing.tex
\tikzset{cross/.style={cross out, draw=black, minimum size=0.5em, inner sep=0pt, outer sep=0pt},cross/.default={10pt}}

\tikzstyle{scaling} = [fill=black,shape=rectangle,minimum size=0.75em,inner sep=0pt,outer sep=0pt,thick]
\tikzstyle{detail} = [fill=black,diamond,very thick,minimum size=0.75em,inner sep=0pt,outer sep=0pt]
\tikzstyle{bullet}=[draw,fill=black,shape=triangle,minimum size=0.5em]
\tikzstyle{temp} = [fill=\figgrey,shape=circle,minimum size=0.5em]

\def\colblockzero{asparagus}
\def\colblockone{cardinal}
\def\colblocktwo{airforceblue}

\def\opacity{0.5}

\begin{tikzpicture}[>=latex']
	\def\shift{0.1}
	\def\step{1};
	\def\l{0}
	\def\ll{-1}
	\def\lll{-2}
	\def\llll{-3}a

   	\draw[-] (4-0.5,\l) -- (10*\step+0.5,\l);
	\draw[-] (4-0.5,\ll) -- (10*\step+0.5,\ll);
	\draw[-] (-0.5,\lll) -- (13*\step+0.5,\lll);
    	
	\draw[-,very thick,dashed,color=\figgrey] (4*\step,1) -- (4*\step,-3);
	\draw[-,very thick,dashed,color=\figgrey] (10*\step,1) -- (10*\step,-3);

	\draw[<-] (-1,\l-\shift) -- node[above,rotate=90] {iDL} (-1,\ll+\shift);
	\draw[<-] (-1,\ll-\shift) -- node[above,rotate=90] {iL} (-1,\lll+\shift);
	
   	 \foreach \x in {0,...,6} {
		\def\y{\l}
		\ifthenelse{\x<2}{
			\foreach \id [
				evaluate=\id as \xe using {int(\x*2 + \id)},
				evaluate=\id as \xid using {int(\x*2  + \id)}] in {0,1}{
			};
       		}{
        			\ifthenelse{\x<5}{
				\foreach \id [
					evaluate=\id as \xe using {int(\x*2 + \id)},
					evaluate=\id as \xid using {int(\x*2 + \id)}
					] in {0,1}{
					\node[scaling,fill=\colblockone,draw=black,label=above:{\textcolor{\colblockone}{${\lambda}^{1}_{\xid}$}}] at ({\x*2 + \id},\y) (scal_\xe) {};
				};

        			}{
				\ifthenelse{\x<5}{
            			\foreach \id [
					evaluate=\id as \xe using {int(\x*2 + \id)},
					evaluate=\id as \xid using {int(\x*2 + \id)}
					] in {0,1}{
					\node[scaling,fill=\colblockone,label=above:{\textcolor{black}{$\lambda^{L+1}_{\xid}$}}] at ({\x*2 + \id},\y) (scal_\xe) {};
				};
				}{
				};
        			};
        		};
		\def\y{\ll}
		\ifthenelse{\x<2}{
			\foreach \id [
				evaluate=\id as \xe using {int(\x*2)},
				evaluate=\id as \xo using {int(\x*2 + 1)},
				evaluate=\id as \xid using {int(\x*2 + \id)}] in {0}{
				
			};
       		}{
        			\ifthenelse{\x<5}{
				\ifthenelse{\x=4 \OR \x=2}{
				\foreach \id [
					evaluate=\id as \xe using {int(\x*2)},
					evaluate=\id as \xo using {int(\x*2 + 1)},
					evaluate=\id as \xid using {int(\x*2 + \id)}] in {0}{
					\node[scaling,fill=\colblockone] at ({\x*2 + \id},\y) (dual_\xe) {};
					\node[detail,fill=\colblockone] at ({\x*2 + \id + 1},\y) (dual_\xo) {};
					\node[cross out,draw=black,thick] at ({\x*2 + \id + 1},\y){};
					
					\ifthenelse{\x=4}{
					\draw[<-,color=black] ($(scal_\xe.south) + (\shift-\step,-\shift)$) -- ($(dual_\xe.north) - (\shift,-\shift)$) ;
					}{};
					\draw[<-,color=black] ($(scal_\xe.south) + (\step-\shift,-\shift)$) -- ($(dual_\xe.north) - (-\shift,-\shift)$) ;
					
					\draw[<-,color=black] ($(scal_\xe.south) + (0,-\shift)$) -- ($(dual_\xe.north) - (0,-\shift)$) ;
				};}{
				\foreach \id [
					evaluate=\id as \xe using {int(\x*2)},
					evaluate=\id as \xo using {int(\x*2 + 1)},
					evaluate=\id as \xid using {int(\x*2 + \id)}] in {0}{
					\node[scaling,fill=\colblockone] at ({\x*2 + \id},\y) (dual_\xe) {};
					\node[detail,fill=\colblockone,draw opacity=\opacity] at ({\x*2 + \id + 1},\y) (dual_\xo) {};
					
					\draw[<-,color=black] ($(scal_\xe.south) + (\shift-\step,-\shift)$) -- ($(dual_\xe.north) - (\shift,-\shift)$) ;
					\draw[<-,color=black] ($(scal_\xe.south) + (\step-\shift,-\shift)$) -- ($(dual_\xe.north) - (-\shift,-\shift)$) ;
					\draw[<-,color=black] ($(scal_\xe.south) + (0,-\shift)$) -- ($(dual_\xe.north) - (0,-\shift)$) ;
					
					\draw[<-,color=black] ($(scal_\xo.south) + (0,-\shift)$) -- ($(dual_\xo.north) - (0,-\shift)$) ;
					
				};
				};
        			}{
			\ifthenelse{\x<6}{
            			\foreach \id [
					evaluate=\id as \xe using {int(\x*2)},
					evaluate=\id as \xo using {int(\x*2 + 1)},
					evaluate=\id as \xid using {int(\x*2 + \id)}
					] in {0}{
						\node[scaling,fill=\colblockone] at ({\x*2 + \id},\y) (dual_\xe) {};
						\ifthenelse{\x=5}{
						\draw[<-,color=black] ($(scal_\xe.south) + (\shift-\step,-\shift)$) -- ($(dual_\xe.north) - (\shift,-\shift)$) ;
						}{};
					
				};
				}{};
        			};
        		};
		\def\y{\lll}
		\ifthenelse{\x<2}{
			\foreach \id [
				evaluate=\id as \xe using {int(\x*2)},
				evaluate=\id as \xo using {int(\x*2 + 1)},
				evaluate=\id as \xid using {int(\x + \id)}] in {0}{
				\node[scaling,draw=black,shape=circle,fill=\colblockzero] at ({\x*2 + \id},\y) (coarse_\xe) {};
				\ifthenelse{\x=1}{
				\node[detail,fill=\colblockzero,label=below:{\textcolor{\colblockzero}{${\gamma}^{0}_{\xid}$}}] at ({\x*2 + \id + 1},\y) (coarse_\xo) {};
				\draw[<-,color=black!70] ($(dual_\xo.south) + (-\shift + \step,-\shift)$) -- node[cross out,draw=black,thick,rotate=0] {} ($(coarse_\xo.north) - (-\shift,-\shift)$) ;
				}{
				\node[detail,fill=\colblockzero] at ({\x*2 + \id + 1},\y) (coarse_\xo) {};
				};
				\node[cross out,draw=black,thick] at ({\x*2 + \id + 1},\y){};
				
			};
       		}{
        			\ifthenelse{\x<5}{
				\ifthenelse{\x=4 \OR \x=2}{
				\foreach \id [
					evaluate=\id as \xe using {int(\x*2)},
					evaluate=\id as \xo using {int(\x*2 + 1)},
					evaluate=\id as \xid using {int(\x + \id )}] in {0}{
					\node[scaling,shape=circle,fill=\colblockone,label=below:{\textcolor{\colblockone}{${\lambda}^{0}_{\xid}$}}] at ({\x*2 + \id},\y)  (coarse_\xe) {};
					\node[detail,fill=\colblockone,label=below:{\textcolor{\colblockone}{${\gamma}^{0}_{\xid}$}}] at ({\x*2 + \id + 1},\y) (coarse_\xo) {};
					\node[cross out,draw=black,thick] at ({\x*2 + \id + 1},\y){};
					
					\draw[<-,color=black!70] ($(dual_\xo.south) + (-\shift + \step,-\shift)$) -- node[cross out,draw=black,thick,rotate=0] {} ($(coarse_\xo.north) - (-\shift,-\shift)$) ;
					\draw[<-,color=black!70] ($(dual_\xo.south) + (+\shift - \step,-\shift)$) -- node[cross out,draw=black,thick,rotate=0] {} ($(coarse_\xo.north) - (+\shift,-\shift)$) ;
					
					\draw[<-,color=black] ($(dual_\xe.south) + (0,-\shift)$) -- ($(coarse_\xe.north) - (0,-\shift)$) ;
				};
				}{
				\foreach \id [
					evaluate=\id as \xe using {int(\x*2)},
					evaluate=\id as \xo using {int(\x*2 + 1)},
					evaluate=\id as \xid using {int(\x + \id)}] in {0}{
					\node[scaling,shape=circle,fill=\colblockone,label=below:{\textcolor{\colblockone}{${\lambda}^{0}_{\xid}$}}] at ({\x*2 + \id},\y) {};
					\node[detail,fill=\colblockone,label=below:{\textcolor{\colblockone}{${\gamma}^{0}_{\xid}$}}] at ({\x*2 + \id + 1},\y) {};
					
					\draw[<-,color=black] ($(dual_\xo.south) + (\shift - \step,-\shift)$) -- ($(coarse_\xo.north) - (\shift,-\shift)$) ;
					\draw[<-,color=black] ($(dual_\xo.south) + (0,-\shift)$) -- ($(coarse_\xo.north) - (0,-\shift)$) ;
					\draw[<-,color=black] ($(dual_\xo.south) + (-\shift + \step,-\shift)$) -- ($(coarse_\xo.north) - (-\shift,-\shift)$) ;
					
					\draw[<-,color=black] ($(dual_\xe.south) + (0,-\shift)$) -- ($(coarse_\xe.north) - (0,-\shift)$) ;
				};
				};
        			}{
            			\foreach \id [
					evaluate=\id as \xe using {int(\x*2)},
					evaluate=\id as \xo using {int(\x*2 + 1)},
					evaluate=\id as \xid using {int(\x  + \id)}
					] in {0}{
						\ifthenelse{\x=5}{
						\node[scaling,draw=black,shape=circle,fill=\colblocktwo,label=below:{\textcolor{\colblocktwo}{${\lambda}^{0}_{\xid}$}}] at ({\x*2 + \id},\y) (coarse_\xe) {};
						}{
						\node[scaling,draw=black,shape=circle,fill=\colblocktwo] at ({\x*2 + \id},\y) (coarse_\xe) {};
						};
						\node[detail,fill=\colblocktwo,label=below:{\textcolor{\colblocktwo}{${\gamma}^{0}_{\xid}$}}] at ({\x*2 + \id + 1},\y) (coarse_\xo) {};
						\node[cross out,draw=black,thick] at ({\x*2 + \id + 1},\y){};
					\ifthenelse{\x=5}{
					\draw[<-,color=black] ($(dual_\xe.south) + (0,-\shift)$) -- ($(coarse_\xe.north) - (0,-\shift)$) ;
					\draw[<-,color=black!70] ($(dual_\xo.south) + (+\shift - \step,-\shift)$) -- node[cross out,draw=black,thick,rotate=0] {} ($(coarse_\xo.north) - (+\shift,-\shift)$) ;
					}{};
				};
        			};
        		};
		
    }
\end{tikzpicture}

%% file: figures/w2.2_1d_ghost_inverse.tex
\tikzset{cross/.style={cross out, draw=black, minimum size=0.5em, inner sep=0pt, outer sep=0pt},cross/.default={10pt}}

\tikzstyle{scaling} = [fill=black,shape=rectangle,minimum size=0.75em,inner sep=0pt,outer sep=0pt,thick]
\tikzstyle{detail} = [fill=black,diamond,very thick,minimum size=0.75em,inner sep=0pt,outer sep=0pt]
\tikzstyle{bullet}=[draw,fill=black,shape=triangle,minimum size=0.5em]
\tikzstyle{temp} = [fill=\figgrey,shape=circle,minimum size=0.5em]

\def\colblockzero{asparagus}
\def\colblockone{cardinal}
\def\colblocktwo{airforceblue}

\def\opacity{0.25}

\begin{tikzpicture}[>=latex']
	\def\shift{0.1}
	\def\step{1};
	\def\l{0}
	\def\ll{-1}
	\def\lll{-2}
	\def\llll{-3}a

   	\draw[-] (-0.5,\l) -- (13*\step+0.5,\l);
	\draw[-] (-0.5,\ll) -- (13*\step+0.5,\ll);
	\draw[-] (-0.5,\lll) -- (13*\step+0.5,\lll);
    	
	\draw[-,very thick,dashed,color=\figgrey] (4*\step,1) -- (4*\step,-3);
	\draw[-,very thick,dashed,color=\figgrey] (10*\step,1) -- (10*\step,-3);

	\draw[<-] (-1,\l) -- node[above,rotate=90] {iDL} (-1,\ll);
	\draw[<-] (-1,\ll) -- node[above,rotate=90] {iL} (-1,\lll);
	
   	 \foreach \x in {0,...,6} {
		\def\y{\l}
		\ifthenelse{\x<2}{
			\ifthenelse{\x=1}{
			\foreach \id [
				evaluate=\id as \xe using {int(\x*2 + \id)},
				evaluate=\id as \xo using {int(\x*2 +1 + \id)},
				evaluate=\id as \xid using {int(\x*2 + \id)},
				evaluate=\id as \xidd using {int(\x*2+1 + \id)}] in {0}{
				\node[scaling,fill=white,thick,draw=\colblockone,label=above:{\textcolor{\colblockone}{$\lambda^{1}_{\xid}$}}] at ({\x*2 + \id},\y) (scal_\xe) {};
				\node[scaling,fill=white,thick,draw=\colblockone,label=above:{\textcolor{\colblockone}{$\lambda^{1}_{\xidd}$}}] at ({\x*2 +1+ \id},\y) (scal_\xo) {};
			};}{};
       		}{
        			\ifthenelse{\x<5}{
				\foreach \id [
					evaluate=\id as \xe using {int(\x*2 + \id)},
					evaluate=\id as \xid using {int(\x*2 + \id)}
					] in {0,1}{
					\node[scaling,fill=\colblockone,draw=black,label=above:{\textcolor{\colblockone}{$\lambda^{1}_{\xid}$}}] at ({\x*2 + \id},\y) (scal_\xe) {};
				};

        			}{
				\ifthenelse{\x<6}{
            			\foreach \id [
					evaluate=\id as \xe using {int(\x*2 + \id)},
					evaluate=\id as \xid using {int(\x*2 + \id)}
					] in {0,1}{
					\node[scaling,fill=white,thick,draw=\colblockone,label=above:{\textcolor{\colblockone}{$\lambda^{1}_{\xid}$}}] at ({\x*2 + \id},\y) (scal_\xe) {};
				};}{
				\foreach \id [
					evaluate=\id as \xe using {int(\x*2 + \id)},
					evaluate=\id as \xid using {int(\x*2 + \id)}
					] in {0}{
					\node[scaling,fill=white,thick,draw=\colblockone,label=above:{\textcolor{\colblockone}{$\lambda^{1}_{\xid}$}}] at ({\x*2 + \id},\y) (scal_\xe) {};
				};
				};
        			};
        		};
		\def\y{\ll}
		\ifthenelse{\x<2}{
			\foreach \id [
				evaluate=\id as \xe using {int(\x*2)},
				evaluate=\id as \xo using {int(\x*2 + 1)},
				evaluate=\id as \xid using {int(\x*2 + \id)}] in {0}{
				
				\ifthenelse{\x>0}{
				\node[scaling,fill=white,draw=\colblockone] at ({\x*2 + \id},\y) (dual_\xe) {};
				\node[detail,fill=\colblockzero] at ({\x*2 + \id + 1},\y) (dual_\xo) {};
				\node[cross out,draw=black,thick] at ({\x*2 + \id + 1},\y){};
								
				\draw[<-,color=black] ($(scal_\xo.south) + (-\shift,-\shift)$) -- ($(dual_\xo.north) - (\step-\shift,-\shift)$) ;
				\draw[<-,color=black] ($(scal_\xo.south) + (+\shift,-\shift)$) -- ($(dual_\xo.north) - (-\step+\shift,-\shift)$) ;
				\draw[<-,color=black] ($(scal_\xe.south) + (0,-\shift)$) -- ($(dual_\xe.north) - (0,-\shift)$) ;
				}{};
			};
       		}{
        			\ifthenelse{\x<5}{
				\foreach \id [
					evaluate=\id as \xe using {int(\x*2)},
					evaluate=\id as \xo using {int(\x*2 + 1)},
					evaluate=\id as \xid using {int(\x*2 + \id)}] in {0}{
					\ifthenelse{\x=2}{
					\node[scaling,draw=\colblockone,fill=white] at ({\x*2 + \id},\y)  (dual_\xe)  {};
					\node[detail,fill=\colblockzero] at ({\x*2 + \id + 1},\y) (dual_\xo) {};
					\node[cross out,draw=black,thick] at ({\x*2 + \id + 1},\y){};
					\draw[->,color=black] ($(scal_\xe.south) + (0,-\shift)$) -- ($(dual_\xe.north) - (0,-\shift)$) ;
					}{};
					\ifthenelse{\x=4}{
					\node[detail,fill=\colblocktwo] at ({\x*2 + \id + 1},\y) (dual_\xo) {};
					\node[cross out,draw=black,thick] at ({\x*2 + \id + 1},\y){};
					}{};

				};
        			}{
            			\foreach \id [
					evaluate=\id as \xe using {int(\x*2)},
					evaluate=\id as \xo using {int(\x*2 + 1)},
					evaluate=\id as \xid using {int(\x*2 + \id)}
					] in {0}{

						\node[scaling,fill=white,draw=\colblockone] at ({\x*2 + \id},\y) (dual_\xe) {};
						\node[detail,fill=\colblocktwo] at ({\x*2 + \id + 1},\y) (dual_\xo) {};
						\node[cross out,draw=black,thick] at ({\x*2 + \id + 1},\y){};
						
					\ifthenelse{\x<6}{
					\draw[<-,color=black] ($(scal_\xo.south) + (-\shift,-\shift)$) -- ($(dual_\xo.north) - (\step-\shift,-\shift)$) ;
					\draw[<-,color=black] ($(scal_\xo.south) + (+\shift,-\shift)$) -- ($(dual_\xo.north) - (-\step+\shift,-\shift)$) ;
					}{};
					\draw[<-,color=black] ($(scal_\xe.south) + (0,-\shift)$) -- ($(dual_\xe.north) - (0,-\shift)$) ;
				};
        			};
        		};
		\def\y{\lll}
		\ifthenelse{\x<2}{
			\foreach \id [
				evaluate=\id as \xe using {int(\x*2)},
				evaluate=\id as \xo using {int(\x*2 + 1)},
				evaluate=\id as \xid using {int(\x + \id)}] in {0}{
				\ifthenelse{\x=1}{
				\node[scaling,draw=black,shape=circle,fill=\colblockzero,label=below:{\textcolor{\colblockzero}{$\lambda^{0}_{\xid}$}}] at ({\x*2 + \id},\y) (coarse_\xe) {};
				}{
				\node[scaling,draw=black,shape=circle,fill=\colblockzero] at ({\x*2 + \id},\y) (coarse_\xe) {};
				};
				\node[detail,fill=\colblockzero] at ({\x*2 + \id + 1},\y) (coarse_\xo) {};
				\node[cross out,draw=black,thick] at ({\x*2 + \id + 1},\y){};
				
				\ifthenelse{\x >0}{
				\draw[<-,color=black] ($(dual_\xe.south) + (0,-\shift)$) -- ($(coarse_\xe.north) - (0,-\shift)$) ;
				}{};
				
			};
       		}{
        			\ifthenelse{\x<5}{
				\ifthenelse{\x=2}{
				\foreach \id [
					evaluate=\id as \xe using {int(\x*2)},
					evaluate=\id as \xo using {int(\x*2 + 1)},
					evaluate=\id as \xid using {int(\x + \id)}] in {0}{
					\node[scaling,shape=circle,fill=white,draw=\colblockzero,thick,label=below:{\textcolor{\colblockzero}{$\lambda^{0}_{\xid}$}}] at ({\x*2 + \id},\y)  (coarse_\xe) {};
					\node[detail,fill=\colblockzero,label=below:{\textcolor{\colblockzero}{$\gamma^{0}_{\xid}$}}] at ({\x*2 + \id + 1},\y) (coarse_\xo) {};
					\node[cross out,draw=black,thick] at ({\x*2 + \id + 1},\y){};

				};
				}{};
				\ifthenelse{\x=3}{
				\foreach \id [
					evaluate=\id as \xe using {int(\x*2)},
					evaluate=\id as \xo using {int(\x*2 + 1)},
					evaluate=\id as \xid using {int(\x + \id)}] in {0}{
					\node[scaling,shape=circle,fill=white,draw=\colblockzero,thick] at ({\x*2 + \id},\y)  (coarse_\xe) {};

				};
				}{};
				\ifthenelse{\x=4}{
				\foreach \id [
					evaluate=\id as \xe using {int(\x*2)},
					evaluate=\id as \xo using {int(\x*2 + 1)},
					evaluate=\id as \xid using {int(\x + \id)}] in {0}{
					\node[scaling,shape=circle,draw=\colblocktwo,fill=white,thick] at ({\x*2 + \id},\y)  (coarse_\xe) {};
					\node[detail,fill=\colblocktwo,label=below:{\textcolor{\colblocktwo}{$\gamma^{0}_{\xid}$}}] at ({\x*2 + \id + 1},\y) (coarse_\xo) {};
					\node[cross out,draw=black,thick] at ({\x*2 + \id + 1},\y){};

				};
				}{};
        			}{
				\ifthenelse{\x=5}{
            			\foreach \id [
					evaluate=\id as \xe using {int(\x*2)},
					evaluate=\id as \xo using {int(\x*2 + 1)},
					evaluate=\id as \xid using {int(\x  + \id)}
					] in {0}{
						\node[scaling,draw=black,shape=circle,fill=\colblocktwo,label=below:{\textcolor{\colblocktwo}{$\lambda^{0}_{\xid}$}}] at ({\x*2 + \id},\y) (coarse_\xe) {};
						\node[detail,fill=\colblocktwo] at ({\x*2 + \id + 1},\y) (coarse_\xo) {};
						\node[cross out,draw=black,thick] at ({\x*2 + \id + 1},\y){};
						
						\draw[<-,color=black] ($(dual_\xe.south) + (0,-\shift)$) -- ($(coarse_\xe.north) - (0,-\shift)$) ;
				
				};
				}{
				\foreach \id [
					evaluate=\id as \xe using {int(\x*2)},
					evaluate=\id as \xo using {int(\x*2 + 1)},
					evaluate=\id as \xid using {int(\x  + \id)}
					] in {0}{
						\node[scaling,draw=black,shape=circle,fill=\colblocktwo,label=below:{\textcolor{\colblocktwo}{$\lambda^{0}_{\xid}$}}] at ({\x*2 + \id},\y) (coarse_\xe) {};
						\node[detail,fill=\colblocktwo] at ({\x*2 + \id + 1},\y) (coarse_\xo) {};
						\node[cross out,draw=black,thick] at ({\x*2 + \id + 1},\y){};

						\draw[<-,color=black] ($(dual_\xe.south) + (0,-\shift)$) -- ($(coarse_\xe.north) - (0,-\shift)$) ;
				
				};
				};
        			};
        		};
    }
\end{tikzpicture}

%% file: figures/w2.2_1d_refinement.tex
\tikzset{cross/.style={cross out, draw=black, minimum size=0.5em, inner sep=0pt, outer sep=0pt},cross/.default={10pt}}

\tikzstyle{scaling} = [fill=black,shape=rectangle,minimum size=0.75em,inner sep=0pt,outer sep=0pt,thick]
\tikzstyle{detail} = [fill=black,diamond,very thick,minimum size=0.75em,inner sep=0pt,outer sep=0pt]
\tikzstyle{bullet}=[draw,fill=black,shape=triangle,minimum size=0.5em]
\tikzstyle{temp} = [fill=\figgrey,shape=circle,minimum size=0.5em]

\def\colblockzero{asparagus}
\def\colblockone{cardinal}
\def\colblocktwo{airforceblue}

\def\opacity{0.5}

\begin{tikzpicture}[>=latex']
	\def\shift{0.1}
	\def\step{1};
	\def\l{0}
	\def\ll{-1}
	\def\lll{-2}
	\def\llll{-3}a

   	\draw[-] (4-0.5,\l) -- (10*\step+0.5,\l);
	\draw[-] (4-0.5,\ll) -- (10*\step+0.5,\ll);
	\draw[-] (2-0.5,\lll) -- (11*\step+0.5,\lll);
    	
	\draw[-,very thick,dashed,color=\figgrey] (4*\step,1) -- (4*\step,-3);
	\draw[-,very thick,dashed,color=\figgrey] (10*\step,1) -- (10*\step,-3);

	\draw[<-] (12*\step,\l-\shift) -- node[above,rotate=90] {iDL} (12*\step,\ll+\shift);
	\draw[<-] (12*\step,\ll-\shift) -- node[above,rotate=90] {iL} (12*\step,\lll+\shift);
	
   	 \foreach \x in {1,...,5} {
		\def\y{\l}
		\ifthenelse{\x<2}{
			\foreach \id [
				evaluate=\id as \xe using {int(\x*2 + \id)},
				evaluate=\id as \xid using {int(\x*2 + \id)}] in {0,1}{
			};
       		}{
        			\ifthenelse{\x<5}{
				\foreach \id [
					evaluate=\id as \xe using {int(\x*2 + \id)},
					evaluate=\id as \xid using {int(\x*2 + \id)}
					] in {0,1}{
					\node[scaling,fill=\colblockone,draw=black,label=above:{\textcolor{\colblockone}{${\lambda}^{1}_{\xid}$}}] at ({\x*2 + \id},\y) (scal_\xe) {};
				};

        			}{
				\foreach \id [
					evaluate=\id as \xe using {int(\x*2 + \id)},
					evaluate=\id as \xid using {int(\x*2  + \id)}
					] in {0}{
					\ifthenelse{\x=5}{
					}{};
				};
        			};
        		};
		\def\y{\ll}
		\ifthenelse{\x<2}{
			\foreach \id [
				evaluate=\id as \xe using {int(\x*2)},
				evaluate=\id as \xo using {int(\x*2 + 1)},
				evaluate=\id as \xid using {int(\x*2 + \id)}] in {0}{
				
			};
       		}{
        			\ifthenelse{\x<5}{
				\foreach \id [
					evaluate=\id as \xe using {int(\x*2)},
					evaluate=\id as \xo using {int(\x*2 + 1)},
					evaluate=\id as \xid using {int(\x*2 + \id)}] in {0}{
					\node[scaling,fill=\colblockone] at ({\x*2 + \id},\y) (dual_\xe) {};
					\node[detail,fill=\colblockone] at ({\x*2 + \id + 1},\y) (dual_\xo) {};
					\node[cross out,draw=black,thick] at ({\x*2 + \id + 1},\y){};

					\ifthenelse{\x=4}{
					\draw[<-,color=black] ($(scal_\xe.south) + (\shift-\step,-\shift)$) -- ($(dual_\xe.north) - (\shift,-\shift)$) ;
					}{};
					
					\ifthenelse{\x>2}{
					\draw[<-,color=black] ($(scal_\xe.south) + (\shift-\step,-\shift)$) -- ($(dual_\xe.north) - (\shift,-\shift)$) ;
					}{};
					\draw[<-,color=black] ($(scal_\xe.south) + (\step-\shift,-\shift)$) -- ($(dual_\xe.north) - (-\shift,-\shift)$) ;
					
					\draw[<-,color=black] ($(scal_\xe.south) + (0,-\shift)$) -- ($(dual_\xe.north) - (0,-\shift)$) ;
				};
        			}{
            			\foreach \id [
					evaluate=\id as \xe using {int(\x*2)},
					evaluate=\id as \xo using {int(\x*2 + 1)},
					evaluate=\id as \xid using {int(\x*2 + \id)}
					] in {0}{

						\ifthenelse{\x=5}{
						\node[scaling,fill=\colblocktwo] at ({\x*2 + \id},\y) (dual_\xe) {};
						\draw[<-,color=black] ($(scal_\xe.south) + (\shift-\step,-\shift)$) -- ($(dual_\xe.north) - (\shift,-\shift)$) ;
					
						}{};
				};
        			};
        		};
		\def\y{\lll}
		\ifthenelse{\x<2}{
			\foreach \id [
				evaluate=\id as \xe using {int(\x*2)},
				evaluate=\id as \xo using {int(\x*2 + 1)},
				evaluate=\id as \xid using {int(\x -2 + \id)}] in {0}{
				\node[scaling,shape=circle,draw=black,fill=\colblockzero] at ({\x*2 + \id},\y) (coarse_\xe) {};
				\node[detail,fill=\colblockzero] at ({\x*2 + \id + 1},\y) (coarse_\xo) {};
				\ifthenelse{\x=1}{
				\node[cross out,draw=black,thick] at ({\x*2 + \id + 1},\y){};
				}{};
				
			};
       		}{
        			\ifthenelse{\x<5}{
				\ifthenelse{\x=4 \OR \x=2 \OR \x=3}{
				\foreach \id [
					evaluate=\id as \xe using {int(\x*2)},
					evaluate=\id as \xo using {int(\x*2 + 1)},
					evaluate=\id as \xid using {int(\x + \id)}] in {0}{
					\node[scaling,draw,shape=circle,fill=\colblockone,label=below:{\textcolor{\colblockone}{${\lambda}^{0}_{\xid}$}}] at ({\x*2 + \id},\y)  (coarse_\xe) {};
					\node[detail,fill=\colblockone] at ({\x*2 + \id + 1},\y) (coarse_\xo) {};
					\node[cross out,draw=black,thick] at ({\x*2 + \id + 1},\y){};

					\draw[<-,color=black] ($(dual_\xe.south) + (0,-\shift)$) -- ($(coarse_\xe.north) - (0,-\shift)$) ;
				};
				}{
				\foreach \id [
					evaluate=\id as \xe using {int(\x*2)},
					evaluate=\id as \xo using {int(\x*2 + 1)},
					evaluate=\id as \xid using {int(\x  + \id)}] in {0}{
					\node[scaling,draw,shape=circle,fill=\colblockone,label=below:{\textcolor{\colblockone}{${\lambda}^{0}_{\xid}$}}] at ({\x*2 + \id},\y) {};
					\node[detail,fill=\colblockone,label=below:{\textcolor{\colblockone}{${\gamma}^{0}_{\xid}$}}] at ({\x*2 + \id + 1},\y) {};
					
					\draw[<-,color=black] ($(dual_\xo.south) + (\shift - \step,-\shift)$) -- ($(coarse_\xo.north) - (\shift,-\shift)$) ;
					\draw[<-,color=black,dashed] ($(dual_\xo.south) + (0,-\shift)$) -- ($(coarse_\xo.north) - (0,-\shift)$) ;
					\draw[<-,color=black] ($(dual_\xo.south) + (-\shift + \step,-\shift)$) -- ($(coarse_\xo.north) - (-\shift,-\shift)$) ;
					
					\draw[<-,color=black] ($(dual_\xe.south) + (0,-\shift)$) -- ($(coarse_\xe.north) - (0,-\shift)$) ;
				};
				};
        			}{
            			\foreach \id [
					evaluate=\id as \xe using {int(\x*2)},
					evaluate=\id as \xo using {int(\x*2 + 1)},
					evaluate=\id as \xid using {int(\x  + \id)}
					] in {0}{
						\node[scaling,draw=black,shape=circle,fill=\colblocktwo,label=below:{\textcolor{\colblocktwo}{$\lambda^{0}_{\xid}$}}] at ({\x*2 + \id},\y) (coarse_\xe) {};
						\node[detail,fill=\colblocktwo] at ({\x*2 + \id + 1},\y) (coarse_\xo) {};
						\node[cross out,draw=black,thick] at ({\x*2 + \id + 1},\y){};						
						
						\ifthenelse{\x=5}{
						\draw[<-,color=black] ($(dual_\xe.south) + (0,-\shift)$) -- ($(coarse_\xe.north) - (0,-\shift)$) ;
						}{};

				};
        			};
        		};
    }
\end{tikzpicture}

%% file: figures/w2.2_1d_substitution.tex
\tikzset{cross/.style={cross out, draw=black, minimum size=0.5em, inner sep=0pt, outer sep=0pt},cross/.default={10pt}}

\tikzstyle{scaling} = [fill=black,shape=rectangle,minimum size=0.75em,inner sep=0pt,outer sep=0pt,thick]
\tikzstyle{detail} = [fill=black,diamond,very thick,minimum size=0.75em,inner sep=0pt,outer sep=0pt]
\tikzstyle{bullet}=[draw,fill=black,shape=triangle,minimum size=0.5em]
\tikzstyle{temp} = [fill=\figgrey,shape=circle,minimum size=0.5em]

\def\colblockzero{asparagus}
\def\colblockone{cardinal}
\def\colblocktwo{airforceblue}

\def\opacity{0.5}

\begin{tikzpicture}[>=latex']
	\def\shift{0.1}
	\def\step{1};
	\def\l{0}
	\def\ll{-1}
	\def\lll{-2}
	\def\llll{-3}a

   	\draw[-] (4-0.5,\l) -- (10*\step+0.5,\l);
	\draw[-] (4-0.5,\ll) -- (10*\step+0.5,\ll);
	\draw[-] (2-0.5,\lll) -- (11*\step+0.5,\lll);
    	
	\draw[-,very thick,dashed,color=\figgrey] (4*\step,1) -- (4*\step,-3);
	\draw[-,very thick,dashed,color=\figgrey] (10*\step,1) -- (10*\step,-3);

	\draw[<-] (12*\step,\l-\shift) -- node[above,rotate=90] {iDL} (12*\step,\ll+\shift);
	\draw[<-] (12*\step,\ll-\shift) -- node[above,rotate=90] {iL} (12*\step,\lll+\shift);
	
   	 \foreach \x in {1,...,5} {
		\def\y{\l}
		\ifthenelse{\x<2}{
			\foreach \id [
				evaluate=\id as \xe using {int(\x*2 + \id)},
				evaluate=\id as \xid using {int(\x*2  + \id)}] in {0,1}{
			};
       		}{
        			\ifthenelse{\x<5}{
				\foreach \id [
					evaluate=\id as \xe using {int(\x*2 + \id)},
					evaluate=\id as \xid using {int(\x*2 + \id)}
					] in {0,1}{
					\node[scaling,fill=\colblockone,draw=black,label=above:{\textcolor{\colblockone}{${\lambda}^{1}_{\xid}$}}] at ({\x*2 + \id},\y) (scal_\xe) {};
				};

        			}{
				\foreach \id [
					evaluate=\id as \xe using {int(\x*2 + \id)},
					evaluate=\id as \xid using {int(\x*2 + \id)}
					] in {0}{
					\ifthenelse{\x=5}{
					\node[scaling,fill=\colblocktwo,label=above:{\textcolor{\colblocktwo}{${\lambda}^{1}_{\xid}$}}] at ({\x*2 + \id},\y) (scal_\xe) {};
					}{};
				};
        			};
        		};
		\def\y{\ll}
		\ifthenelse{\x<2}{
			\foreach \id [
				evaluate=\id as \xe using {int(\x*2)},
				evaluate=\id as \xo using {int(\x*2 + 1)},
				evaluate=\id as \xid using {int(\x*2 + \id)}] in {0}{
				
			};
       		}{
        			\ifthenelse{\x<5}{
				\foreach \id [
					evaluate=\id as \xe using {int(\x*2)},
					evaluate=\id as \xo using {int(\x*2 + 1)},
					evaluate=\id as \xid using {int(\x*2 + \id)}] in {0}{
					\node[scaling,fill=\colblockone] at ({\x*2 + \id},\y) (dual_\xe) {};
					\node[detail,fill=\colblockone] at ({\x*2 + \id + 1},\y) (dual_\xo) {};
					
					\ifthenelse{\x=4 \OR \x=2}{
					\node[cross out,draw=black,thick] at ({\x*2 + \id + 1},\y){};
					}{};

					\ifthenelse{\x=3}{
					\draw[<-,color=black] ($(scal_\xe.south) + (\shift-\step,-\shift)$) -- ($(dual_\xe.north) - (\shift,-\shift)$) ;
					}{
					};
					\ifthenelse{\x=2 \OR \x=4}{
					\draw[<-,color=black] ($(scal_\xe.south) + (\step-\shift,-\shift)$) -- ($(dual_\xe.north) - (-\shift,-\shift)$) ;
					}{};
					
					\draw[->,color=black] ($(scal_\xe.south) + (0,-\shift)$) -- ($(dual_\xe.north) - (0,-\shift)$) ;
				};
        			}{
            			\foreach \id [
					evaluate=\id as \xe using {int(\x*2)},
					evaluate=\id as \xo using {int(\x*2 + 1)},
					evaluate=\id as \xid using {int(\x*2 + \id)}
					] in {0}{

						\ifthenelse{\x=5}{
						\node[scaling,fill=\colblocktwo] at ({\x*2 + \id},\y) (dual_\xe) {};
						\draw[<-,color=black] ($(scal_\xe.south) + (\shift-\step,-\shift)$) -- ($(dual_\xe.north) - (\shift,-\shift)$) ;
						\draw[->,color=black] ($(scal_\xe.south) + (0,-\shift)$) -- ($(dual_\xe.north) - (0,-\shift)$) ;
					
						}{};
				};
        			};
        		};
		\def\y{\lll}
		\ifthenelse{\x<2}{
			\foreach \id [
				evaluate=\id as \xe using {int(\x*2)},
				evaluate=\id as \xo using {int(\x*2 + 1)},
				evaluate=\id as \xid using {int(\x -2 + \id)}] in {0}{
				\node[scaling,draw,shape=circle,fill=\colblockzero] at ({\x*2 + \id},\y) (coarse_\xe) {};
				\node[detail,fill=\colblockzero] at ({\x*2 + \id + 1},\y) (coarse_\xo) {};
				\ifthenelse{\x=1}{
				\node[cross out,draw=black,thick] at ({\x*2 + \id + 1},\y){};
				}{};
				
			};
       		}{
        			\ifthenelse{\x<5}{
				\ifthenelse{\x=4 \OR \x=2 \OR \x=3}{
				\foreach \id [
					evaluate=\id as \xe using {int(\x*2)},
					evaluate=\id as \xo using {int(\x*2 + 1)},
					evaluate=\id as \xid using {int(\x + \id - 2)}] in {0}{
					\node[scaling,shape=circle,fill=\colblockone] at ({\x*2 + \id},\y)  (coarse_\xe) {};
					\node[detail,fill=\colblockone] at ({\x*2 + \id + 1},\y) (coarse_\xo) {};
					\ifthenelse{\x=4 \OR \x=2}{
					\node[cross out,draw=black,thick] at ({\x*2 + \id + 1},\y){};
					}{};

				};
				}{
				\foreach \id [
					evaluate=\id as \xe using {int(\x*2)},
					evaluate=\id as \xo using {int(\x*2 + 1)},
					evaluate=\id as \xid using {int(\x - 2 + \id)}] in {0}{
					\node[scaling,draw,shape=circle,fill=\colblockone,label=below:{\textcolor{\colblockone}{${\lambda}^{0}_{\xid}$}}] at ({\x*2 + \id},\y) {};
					\node[detail,fill=\colblockone,label=below:{\textcolor{\colblockone}{${\gamma}^{0}_{\xid}$}}] at ({\x*2 + \id + 1},\y) {};

				};
				};
        			}{
            			\foreach \id [
					evaluate=\id as \xe using {int(\x*2)},
					evaluate=\id as \xo using {int(\x*2 + 1)},
					evaluate=\id as \xid using {int(\x  + \id)}
					] in {0}{
						\node[scaling,draw,shape=circle,fill=\colblocktwo] at ({\x*2 + \id},\y) (coarse_\xe) {};
						\node[detail,fill=\colblocktwo] at ({\x*2 + \id + 1},\y) (coarse_\xo) {};
						\node[cross out,draw=black,thick] at ({\x*2 + \id + 1},\y){};						
						
						\ifthenelse{\x=5}{
						}{};

				};
        			};
        		};
    }
\end{tikzpicture}

%% file: wavelet_3_implementation.tex
\section{Implementation and algorithms}
\label{sect_implementation_global}
In this section we discuss high-level implementation choices of the ghost reconstruction and grid adaptation operations, deferring the details to \sect{sec_sm_algo} in the Supplementary Materials. 
Our entire code base relies on the external library \code{p4est} \citep{Burstedde:2011} to handle the meta-data infrastructure of the octree, while all grid adaptation and block operations have been implemented directly in our solver. 
Within \code{p4est}, we define each `tree' as a unit cube domain, which forms the root (level $L=0$) of an octree data structure that can be refined, \revtwo{and} the leaves are \revtwo{uniform} resolution blocks of size $N_b^3$. Following \code{p4est}, the trees can be tiled to create a `forest' of trees, which enables us to create rectangular domains of arbitrary aspect ratios. 
We currently have implemented wavelets with \revtwo{$N \in \{2,4,6\}$ and $\tilde{N}\in\{0,2\}$.}. \revtwo{Extension} to higher $N$ is straightforward if needed; higher $\tilde{N}$ on the other hand will potentially deteriorate the efficiency of the solver as \revtwo{we will have to significantly increase the number of adjacent detail coefficients that needs to be discarded in accordance with our coarse-extension assumption}.
Throughout our implementation we apply the 2:1 constraint on levels of adjacent blocks, enforcing it during grid adaptation and exploiting it during all multiresolution wavelet operations. 

\paragraph{Ghost reconstruction procedure}
The implementation of the ghost reconstruction consists of two parts: the first recovers the value from coarser neighbors and same level neighbors (see \algref{algo:ghost_1} in \sect{sec_sm_algo}), and the second computes the values from finer neighbors (see \algref{algo:ghost_2} in \sect{sec_sm_algo}). The inter-rank communication during both parts is handled using MPI-RMA with a Post-Start-Complete-Wait (PSCW) synchronization strategy, chosen to target massively parallel infrastructures \citep{Gropp:2014}. %
In the first step, we copy (or use \code{MPI\_Get}) the required scaling coefficients from the coarse- and same-level neighbors to the current block, where the actual ghost points are computed locally once the required values are gathered. This choice reduces the size of the communications and the required memory for the buffers. 
For similar reasons, to retrieve ghost values from a finer neighbor block, we first compute the required ghost values from the perspective of the neighboring block, which then copies (or uses \code{MPI\_Put}) the coarsened values to our block.
Throughout these steps, we have implemented the ghost point computation for vector- or tensor fields in a ``component-by-component'' way, so that we can overlap communication and computation by performing wavelet-based refinement/coarsening operations on one component while the ghost exchange is performed for the next one.

\paragraph{Grid adaptation}
The implementation of the grid adaptation follows an iterative procedure. During each iteration we start by computing the maximum detail criterion for each block as explained in section~\sect{subsubsec:adaptation}, which in turn dictates the block's desired actions. Then we enforce global policies that include the 2:1 condition, possible user-defined limits of minimum/maximum level (not used in this manuscript), and the prohibition of coarsening blocks that have been refined at earlier iterations; a detailed description is included in \sect{sec_sm_algo}. This finalizes the adaptation decision on each block, after which we perform the refinement and/or coarsening on the affected blocks and use the update step to adjust the scaling coefficients of blocks whose neighbors have just been coarsened. Finally, we use the \code{p4est} grid partitioning algorithm to distribute the blocks among ranks and ensure load balancing of the current grid. This ends the current iteration, after which we recompute the ghost values. The iterative process ends when, under our policy, no blocks have changed their resolution, which ensures that $\normi{\gamma} < \epsilon_r$ on all blocks.

%% file: wavelet_4_validation.tex
\section{Validation}
\label{sect_validation_global}
We present here the numerical validation of our framework on three different aspects: the grid adaptation and the error control, the moment conservation for the lifted wavelets, and the convergence of finite difference \revtwo{operators on multi-level grids}. For all cases, we set the linear dimension of each block to $N_b = 24$, so that each block contains $N_b^3 = 13\,824$ unknowns.

\subsection{Grid adaptation and error control}
\label{sec:epsilon_test}
The grid adaptation test, referred to as the ``epsilon test'', measures the error between a coarsened field and the original, non-compressed information.
According to the wavelet theory and \eqqref{eq:error_mr} this error must be bounded by $C_1 \epsilon_c$, where we observed in practice $C_1 = \mathcal{O}(1)$.
For a fixed value of $\epsilon_c$, the epsilon test proceeds as follows:
\begin{enumerate}
    \item initialize an analytic field on a fine level $L_{\max}$,
    \item given $\epsilon_c$, coarsen the grid according to the block adaptation policy described above,
    \item refine the grid back to the $L_{\max}$ level and compare the error with the initial condition.
\end{enumerate}

\begin{figure}[ht!]
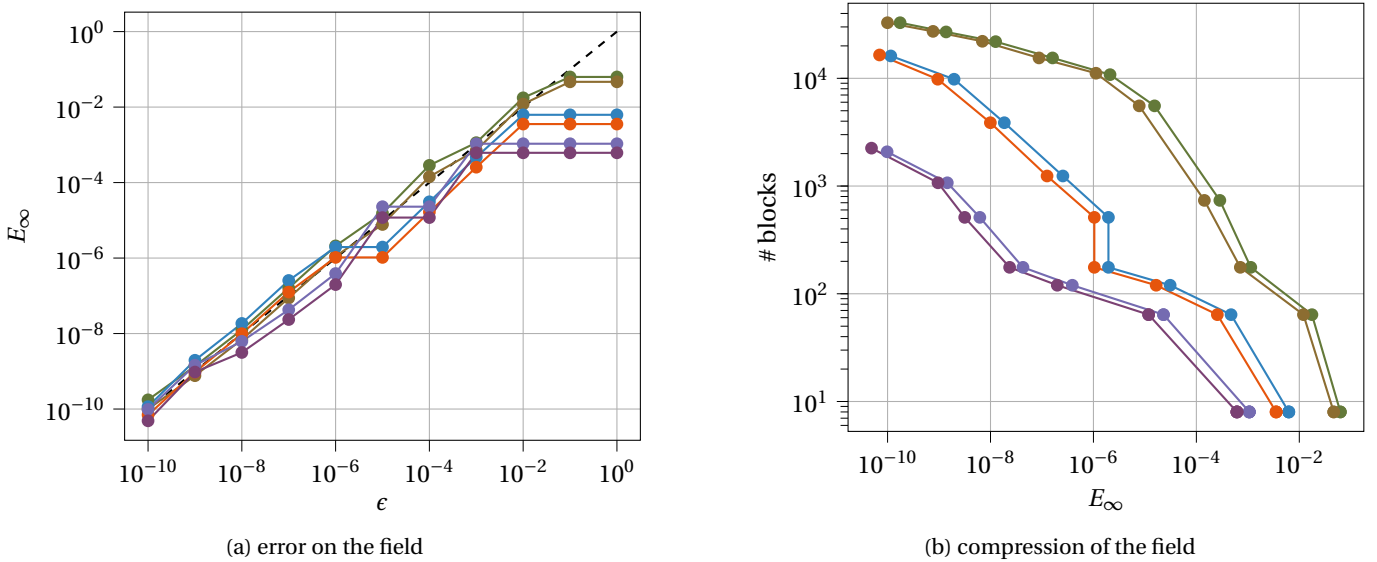

\begin{minipage}[t]{0.48\textwidth}
    \centering
    \InputIfFileExists{figures/epsilon_epsilon}{}{}
    \subcaption{error on the field}
    \label{fig:epsilon_epsilon}
\end{minipage}%
\hfill%
\begin{minipage}[t]{0.48\textwidth}
    \centering
    \InputIfFileExists{figures/epsilon_nblocks}{}{}
    \subcaption{compression of the field}
    \label{fig:epsilon_nblock}
\end{minipage}
 \caption{Epsilon-test: effect of the compression threshold $\epsilon_c$ on error (left) and number of blocks (right) for a static grid adaptation using \wave{2}{0} (\CaptionTwoZero), \waveshort{2}{2} (\CaptionTwoTwo), \waveshort{4}{0} (\CaptionFourZero), \waveshort{4}{2} (\CaptionFourTwo), \waveshort{6}{0} (\CaptionSixZero), \waveshort{6}{2} (\CaptionSixTwo).}
 \label{fig:epsilon_all}
\end{figure}

We have chosen the analytical field to be a scalar Gaussian function centered within a cubic computational domain of size $2$:
\begin{eqc}
f(x,y,z) = \expo{ -\dfrac{r^2}{\sigma^2}} \qquad r^2 = \left( x - 1 \right)^2 + \left( y - 1 \right)^2 + \left( z - 1 \right)^2
\label{eq:compact_scalar_ring}
\end{eqc}
where we set $\sigma = 2/15$. The field is initialized at $L_{\max} = 5$. 
In \fref{fig:epsilon_epsilon} we show the evolution of the infinite norm of the error $E_\infty$ depending on the value of $\epsilon_c$, for a range of different wavelets, where the error is defined as
\begin{eqd}
E_\infty = \norm{ \projval{f}{L_{\max}} (\x) - \projval{f}{L_{\max}}_{\epsilon_c}(\x) }_{\infty} 
\end{eqd}
The results validate that the $\epsilon_c$ is an accurate prediction of the compression error, consistent with the 1D wavelet theory described above. For high values of $\epsilon_c$, the error plateaus as the block granularity in the grid is too low to allow further coarsening. Looking at the number of blocks as a function of the $E_\infty$ in \fref{fig:epsilon_nblock}, we observe that for a given error the number of blocks required to represent the compressed field decreases significantly as the wavelet order $N$ increases. Further, the lifting wavelets characterized by $\tilde{N}=2$ consistently require a slightly smaller number of blocks than their non-lifting counterpart characterized by $\tilde{N}=0$, for the same error and interpolation order $N$.

\subsection{Moment conservation}
\begin{figure}[ht!]
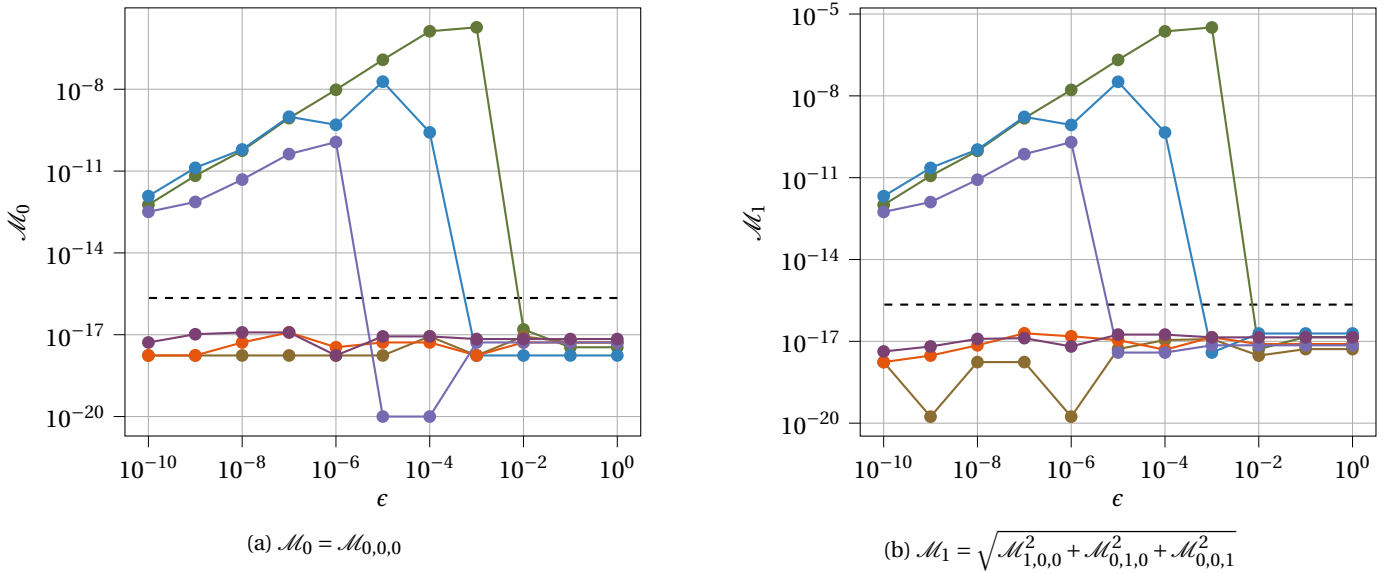

\begin{minipage}{0.48\textwidth}
    \centering
    \InputIfFileExists{figures/epsilon_moments_0}{}{}
    \subcaption{$\mathcal{M}_0 = \mathcal{M}_{0,0,0}$}
    \label{fig:epsilon_moments_0}
\end{minipage}%
\hfill%
\begin{minipage}{0.48\textwidth}
    \centering
    \InputIfFileExists{figures/epsilon_moments_1}{}{}
    \subcaption{$\mathcal{M}_1 =  \sqrt{ \mathcal{M}^2_{1,0,0} + \mathcal{M}^2_{0,1,0} + \mathcal{M}_{0,0,1}^2 }$}
    \label{fig:epsilon_moments_1}
\end{minipage}
\caption{Moment test: effect of the compression threshold $\epsilon_c$ on the conservation of zero (left) and first (right) moments before and after static grid adaptation, for the \wave{2}{0} (\CaptionTwoZero), \waveshort{2}{2} (\CaptionTwoTwo), \waveshort{4}{0} (\CaptionFourZero), \waveshort{4}{2} (\CaptionFourTwo), \waveshort{6}{0} (\CaptionSixZero), \waveshort{6}{2} (\CaptionSixTwo).}
\label{fig:epsilon_moments}
\end{figure}

The moment-preserving properties of lifting wavelets described above can be validated by comparing the moments on a given uniform level $L_{\max}$, both before and after discarding the detail coefficients according to $\epsilon_c$.  
Using the same setup as for the epsilon-test, we compare moments between $\projval{f}{L_{\max}} (x)$ and $\projval{f}{L_{\max}}_{\epsilon_c} (x)$ and define their difference as
\begin{eqc}
\mathcal{M}_{p,q,r} = \int_{\mathbb{R}^3}  \projval{f}{L_{\max}} (\x)  \; x^p y^q z^r \; \dbx - \int_{\mathbb{R}^3}  \projval{f}{L_{\max}}_{\epsilon_c} (\x)  \; x^p y^q z^r \; \dbx
\end{eqc}
where $0 \le {p,q,r} < \tilde{N}$. Each moment can be evaluated from the scaling coefficients using \eqqref{eq:moment_interpolating_1d}.
The results of this test are shown for the different wavelets in \frefs{fig:epsilon_moments_0}{fig:epsilon_moments_1} respectively for the zeroth moment $p=q=r=0$ and the norm of the three first moments. This validates that the lifted interpolating wavelets conserve both the zeroth and the first moment of the scaling coefficients throughout the adaptation process. When considering the non-lifted wavelet family, we notice that the error in the moments is negligible for the higher $\epsilon_c$ values, then suddenly increases at a certain $\epsilon_c$ and gradually decrease when $\epsilon_c$ is reduced. At the largest values for $\epsilon_c$, the adaptation process coarsens the grid uniformly, and by virtue of the dual scaling functions with $\tilde{N}=0$ we retain the original function values on the remaining grid points. Refinement does not affect the moments of the field for any interpolating wavelet (see \eqqref{eq:moments_scaling_function}) and so in this case our test will compare moments between function values on two uniform grids at different resolutions. On a uniform grid the moment integration rule is equivalent to a spectrally accurate trapezoid quadrature due to the compactness and smoothness of the Gaussian function, and it turns out that even for the relatively high values of $\epsilon_c$ considered the coarsening does not affect the error of this approximation, leading to zero values in the moment error. For each of the three wavelets with $\tilde{N} = 0$, there exist a ``critical'' value of $\epsilon_c$ for which the coarse grid first contains multiple levels, thus breaking the favorable convergence properties associated with a uniform grid quadrature and showing the real effect of grid adaptation on the lack of moment conservation for these wavelets.

\subsection{Wavelets and spatial discretization}
\label{sect_conv_weno}
As a third measure of a static validation of our 3D wavelet-based multiresolution grid framework, we consider convergence of various finite-difference operators across a resolution jump as a function of the wavelet order and refinement level. We consider an advection term discretized using conservative upwind finite difference schemes of third order (\cons{3}) and fifth order (\cons{5}), as well as central laplacian operators of second order (\diff{2}) and fourth order (\diff{4}).
More details about the finite-difference schemes used can be found in \sect{sec_finite_diff}.

The analytical field is the same Gaussian blob as in the previous two subsections, and is initialized on a uniform fine level characterized by grid spacing $h_f$. To simplify the measure of convergence, we do not consider the automatic mesh adaptation in this subsection and instead focus on a grid with two levels of resolution that are fixed in space. Starting from the initial condition at level $h_f$, we coarsen one \revtwo{eighth} of the grid by one level, making sure we cover all possible resolution jumps between blocks (jumps across faces, edges, and corners). We then compute the ghost points as described above and evaluate the finite-difference stencil on the entire grid. To compute the error, we compare the discrete values to the analytic solution of applying the continuous differential operators to the analytic field. 

\begin{figure}[ht]
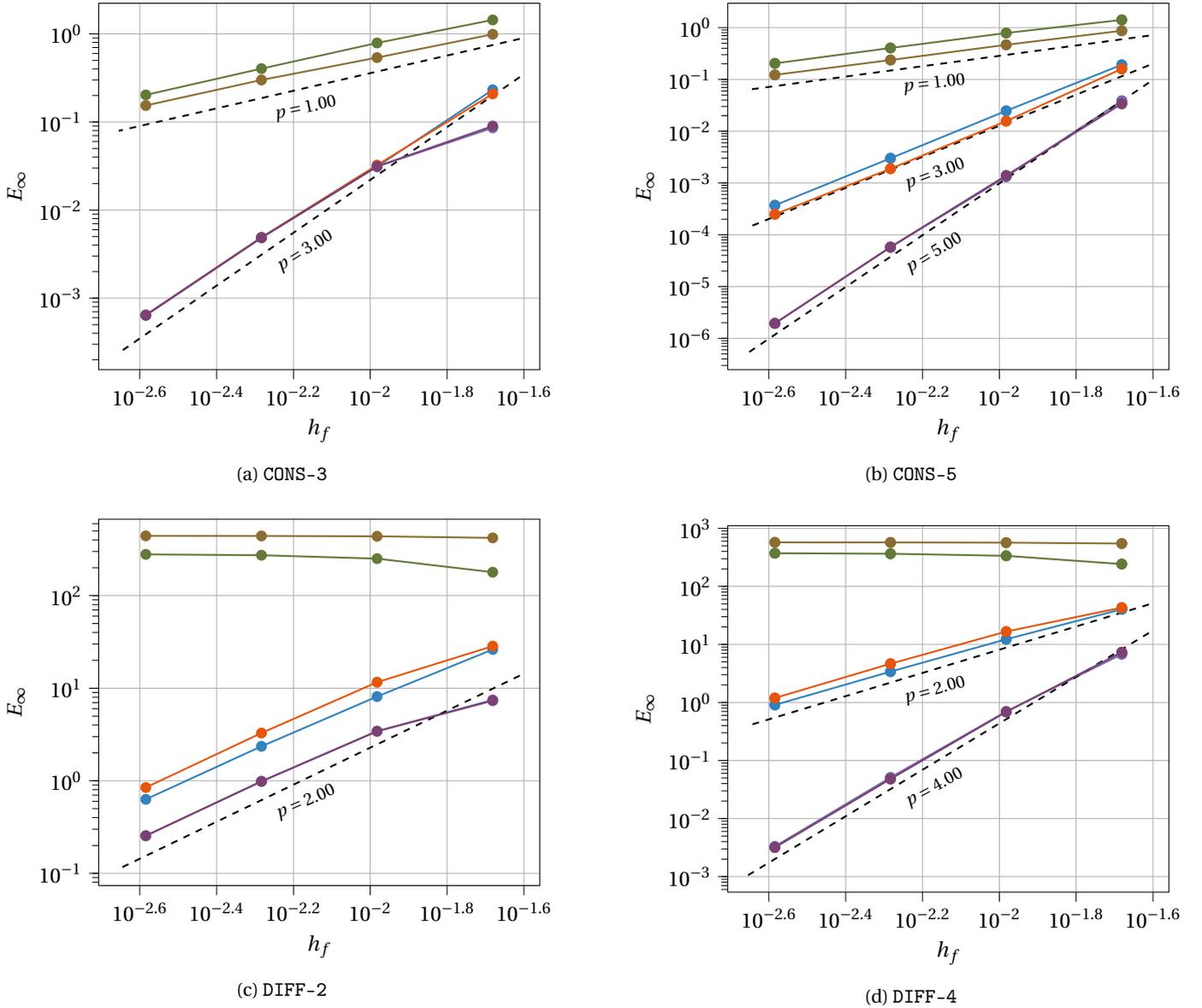

\centering
\begin{minipage}{0.48\textwidth}
    \centering
    \InputIfFileExists{figures/twolevels-cons-3_adapt_normi.tex}{}{}
    \subcaption{\cons{3}}
    \label{fig:twolevels-cons-3_adapt_normi}
\end{minipage}%
\hfill%
\begin{minipage}{0.48\textwidth}
    \centering
    \InputIfFileExists{figures/twolevels-cons-5_adapt_normi.tex}{}{}
    \subcaption{\cons{5}}
    \label{fig:twolevels-cons-5_adapt_normi}
\end{minipage}
\begin{minipage}{0.48\textwidth}
    \centering
    \InputIfFileExists{figures/twolevels-diff-2_adapt_normi.tex}{}{}
    \subcaption{\diff{2}}
    \label{fig:twolevels-diff-2_adapt_normi}
\end{minipage}%
\hfill%
\begin{minipage}{0.48\textwidth}
    \centering
    \InputIfFileExists{figures/twolevels-diff-4_adapt_normi.tex}{}{}
    \subcaption{\diff{4}}
    \label{fig:twolevels-diff-4_adapt_normi}
\end{minipage}%

\caption{Convergence order of spatial derivatives for third and fifth order conservative advection (top row) and a second and fourth order laplacian (bottom row) on a two-level grid with finest level $h_f$, obtained using ghost points reconstructed with \wave{2}{0} (\CaptionTwoZero), \waveshort{2}{2} (\CaptionTwoTwo), \waveshort{4}{0} (\CaptionFourZero), \waveshort{4}{2} (\CaptionFourTwo), \waveshort{6}{0} (\CaptionSixZero), \waveshort{6}{2} (\CaptionSixTwo).}
\label{fig:twolevels}
\end{figure}

The convergence of the infinite norm ($E_{\infty}$) of the error is shown in \fref{fig:twolevels} as a function of $h_f$, for different wavelet orders $N$ and $\tilde{N}$. The results show that if the wavelet order is sufficiently high, the expected convergence order is reached for all finite-difference operators even in the infinity norm, indicating a correct treatment of the resolution jump. 
For lower-order wavelets, the error instead is bound by the polynomial order of the wavelet used to interpolate the fine-level ghost points. Specifically, in this case the error is bound by $N-n$, where $N$ is the wavelet order and $n$ is the order of the derivative operator, consistent with the accuracy order of numerically differentiating an $N$th degree polynomial $n$ times. 
Across all cases, the convergence of the error can thus be given as $h^p$ with $p = \min(k, N - n)$, where $k$ is order of the finite-difference stencil. 

\revtwo{In practice} we \revtwo{should} therefore only consider \revtwo{wavelets} with $N \geq 4$ to obtain a scheme that is at least second-order accurate on first- and second-order PDEs. %

%% file: figures/epsilon_epsilon.tex
\begin{tikzpicture}

\definecolor{color0}{rgb}{0.388235294117647,0.474509803921569,0.223529411764706}
\definecolor{color1}{rgb}{0.549019607843137,0.427450980392157,0.192156862745098}
\definecolor{color2}{rgb}{0.192156862745098,0.509803921568627,0.741176470588235}
\definecolor{color3}{rgb}{0.901960784313726,0.333333333333333,0.0509803921568627}
\definecolor{color4}{rgb}{0.458823529411765,0.419607843137255,0.694117647058824}
\definecolor{color5}{rgb}{0.482352941176471,0.254901960784314,0.450980392156863}

\begin{axis}[
log basis x={10},
log basis y={10},
tick align=outside,
tick pos=left,
x grid style={white!69.0196078431373!black},
xlabel={\(\displaystyle  \epsilon \)},
xmajorgrids,
xmin=3.16227766016838e-11, xmax=3.16227766016838,
xmode=log,
xtick style={color=black},
y grid style={white!69.0196078431373!black},
ylabel={\(\displaystyle E_{\infty}\)},
ymajorgrids,
ymin=1.48690789857195e-11, ymax=3.27797412081759,
ymode=log,
ytick style={color=black}
]
\addplot [thick, black, dashed]
table {%
1 1
0.1 0.1
0.01 0.01
0.001 0.001
0.0001 0.0001
1e-05 1e-05
1e-06 1e-06
1e-07 1e-07
1e-08 1e-08
1e-09 1e-09
1e-10 1e-10
};
\addplot [thick, color0, mark=*, mark size=2, mark options={solid}]
table {%
1 0.0628686814987995
0.1 0.0628686814987995
0.01 0.0175992264638142
0.001 0.00114157529956826
0.0001 0.000285925001987875
1e-05 1.54022242191271e-05
1e-06 2.10289742266518e-06
1e-07 1.60576800675214e-07
1e-08 1.2590034182122e-08
1e-09 1.37055649159158e-09
1e-10 1.7524220184871e-10
};
\addplot [thick, color1, mark=*, mark size=2, mark options={solid}]
table {%
1 0.0466662182045596
0.1 0.0466662182045596
0.01 0.0120284956880223
0.001 0.000714696273196314
0.0001 0.000143051076281875
1e-05 7.75605325036011e-06
1e-06 1.12050297085275e-06
1e-07 8.75896641442959e-08
1e-08 6.95212358795031e-09
1e-09 7.64855083387585e-10
1e-10 9.95444831418701e-11
};
\addplot [thick, color2, mark=*, mark size=2, mark options={solid}]
table {%
1 0.0062453108586551
0.1 0.0062453108586551
0.01 0.0062453108586551
0.001 0.000472134311091499
0.0001 3.0940244805544e-05
1e-05 1.95677929215599e-06
1e-06 1.95677929215599e-06
1e-07 2.54967572138609e-07
1e-08 1.85471604463583e-08
1e-09 1.95990148069147e-09
1e-10 1.15435624763443e-10
};
\addplot [thick, color3, mark=*, mark size=2, mark options={solid}]
table {%
1 0.00354872221107616
0.1 0.00354872221107616
0.01 0.00354872221107616
0.001 0.000255599324239908
0.0001 1.65581991280472e-05
1e-05 1.04047973947807e-06
1e-06 1.04047973947807e-06
1e-07 1.25588684632416e-07
1e-08 9.9364203254293e-09
1e-09 9.45532962690065e-10
1e-10 6.97971424038456e-11
};
\addplot [thick, color4, mark=*, mark size=2, mark options={solid}]
table {%
1 0.00106912203691478
0.1 0.00106912203691478
0.01 0.00106912203691478
0.001 0.00106912203691478
0.0001 2.29643935029111e-05
1e-05 2.29643935029111e-05
1e-06 3.89012644541609e-07
1e-07 4.26999859082233e-08
1e-08 6.20287488040105e-09
1e-09 1.44722595019076e-09
1e-10 9.78992997779926e-11
};
\addplot [thick, color5, mark=*, mark size=2, mark options={solid}]
table {%
1 0.000614449125227901
0.1 0.000614449125227901
0.01 0.000614449125227901
0.001 0.000614449125227901
0.0001 1.18216418776296e-05
1e-05 1.18216418776296e-05
1e-06 1.98253382377445e-07
1e-07 2.35277043341363e-08
1e-08 3.15251869054123e-09
1e-09 9.55037215977228e-10
1e-10 4.87404561155813e-11
};
\end{axis}

\end{tikzpicture}

%% file: figures/epsilon_nblocks.tex
\begin{tikzpicture}

\definecolor{color0}{rgb}{0.388235294117647,0.474509803921569,0.223529411764706}
\definecolor{color1}{rgb}{0.549019607843137,0.427450980392157,0.192156862745098}
\definecolor{color2}{rgb}{0.192156862745098,0.509803921568627,0.741176470588235}
\definecolor{color3}{rgb}{0.901960784313726,0.333333333333333,0.0509803921568627}
\definecolor{color4}{rgb}{0.458823529411765,0.419607843137255,0.694117647058824}
\definecolor{color5}{rgb}{0.482352941176471,0.254901960784314,0.450980392156863}

\begin{axis}[
log basis x={10},
log basis y={10},
tick align=outside,
tick pos=left,
x grid style={white!69.0196078431373!black},
xlabel={\(\displaystyle E_{\infty}\)},
xmajorgrids,
xmin=1.7075064420422e-11, xmax=0.179457490536423,
xmode=log,
xtick style={color=black},
y grid style={white!69.0196078431373!black},
ylabel={\# blocks},
ymajorgrids,
ymin=5.27758104358728, ymax=49756.1283912584,
ymode=log,
ytick style={color=black}
]
\addplot [thick, color0, mark=*, mark size=2, mark options={solid}]
table {%
0.0628686814987995 8
0.0628686814987995 8
0.0175992264638142 64
0.00114157529956826 176
0.000285925001987875 736
1.54022242191271e-05 5552
2.10289742266518e-06 10816
1.60576800675214e-07 15464
1.2590034182122e-08 21904
1.37055649159158e-09 26944
1.7524220184871e-10 32824
};
\addplot [thick, color1, mark=*, mark size=2, mark options={solid}]
table {%
0.0466662182045596 8
0.0466662182045596 8
0.0120284956880223 64
0.000714696273196314 176
0.000143051076281875 736
7.75605325036011e-06 5552
1.12050297085275e-06 11152
8.75896641442959e-08 15464
6.95212358795031e-09 22072
7.64855083387585e-10 27280
9.95444831418701e-11 32824
};
\addplot [thick, color2, mark=*, mark size=2, mark options={solid}]
table {%
0.0062453108586551 8
0.0062453108586551 8
0.0062453108586551 8
0.000472134311091499 64
3.0940244805544e-05 120
1.95677929215599e-06 176
1.95677929215599e-06 512
2.54967572138609e-07 1240
1.85471604463583e-08 3872
1.95990148069147e-09 9808
1.15435624763443e-10 16136
};
\addplot [thick, color3, mark=*, mark size=2, mark options={solid}]
table {%
0.00354872221107616 8
0.00354872221107616 8
0.00354872221107616 8
0.000255599324239908 64
1.65581991280472e-05 120
1.04047973947807e-06 176
1.04047973947807e-06 512
1.25588684632416e-07 1240
9.9364203254293e-09 3872
9.45532962690065e-10 9808
6.97971424038456e-11 16472
};
\addplot [thick, color4, mark=*, mark size=2, mark options={solid}]
table {%
0.00106912203691478 8
0.00106912203691478 8
0.00106912203691478 8
0.00106912203691478 8
2.29643935029111e-05 64
2.29643935029111e-05 64
3.89012644541609e-07 120
4.26999859082233e-08 176
6.20287488040105e-09 512
1.44722595019076e-09 1072
9.78992997779926e-11 2080
};
\addplot [thick, color5, mark=*, mark size=2, mark options={solid}]
table {%
0.000614449125227901 8
0.000614449125227901 8
0.000614449125227901 8
0.000614449125227901 8
1.18216418776296e-05 64
1.18216418776296e-05 64
1.98253382377445e-07 120
2.35277043341363e-08 176
3.15251869054123e-09 512
9.55037215977228e-10 1072
4.87404561155813e-11 2248
};
\end{axis}

\end{tikzpicture}

%% file: figures/epsilon_moments_0.tex
\begin{tikzpicture}

\definecolor{color0}{rgb}{0.388235294117647,0.474509803921569,0.223529411764706}
\definecolor{color1}{rgb}{0.549019607843137,0.427450980392157,0.192156862745098}
\definecolor{color2}{rgb}{0.192156862745098,0.509803921568627,0.741176470588235}
\definecolor{color3}{rgb}{0.901960784313726,0.333333333333333,0.0509803921568627}
\definecolor{color4}{rgb}{0.458823529411765,0.419607843137255,0.694117647058824}
\definecolor{color5}{rgb}{0.482352941176471,0.254901960784314,0.450980392156863}

\begin{axis}[
log basis x={10},
log basis y={10},
tick align=outside,
tick pos=left,
x grid style={white!69.0196078431373!black},
xlabel={\(\displaystyle  \epsilon \)},
xmajorgrids,
xmin=3.16227766016838e-11, xmax=3.16227766016838,
xmode=log,
xtick style={color=black},
y grid style={white!69.0196078431373!black},
ylabel={\(\displaystyle \mathcal{M}_0\)},
ymajorgrids,
ymin=1.93413401757164e-21, ymax=9.63345713172627e-06,
ymode=log,
ytick style={color=black}
]
\addplot [thick, black, dashed]
table {%
1 2.22044604925031e-16
0.1 2.22044604925031e-16
0.01 2.22044604925031e-16
0.001 2.22044604925031e-16
0.0001 2.22044604925031e-16
1e-05 2.22044604925031e-16
1e-06 2.22044604925031e-16
1e-07 2.22044604925031e-16
1e-08 2.22044604925031e-16
1e-09 2.22044604925031e-16
1e-10 2.22044604925031e-16
};
\addplot [thick, color0, mark=*, mark size=2, mark options={solid}]
table {%
1 3.46944695195361e-18
0.1 3.46944695195361e-18
0.01 1.56125112837913e-17
0.001 1.86323971452899e-06
0.0001 1.33539192098513e-06
1e-05 1.20163403618048e-07
1e-06 9.54888163158529e-09
1e-07 8.7193083656012e-10
1e-08 5.56321308464724e-11
1e-09 6.73600238088046e-12
1e-10 5.81640638430692e-13
};
\addplot [thick, color1, mark=*, mark size=2, mark options={solid}]
table {%
1 5.20417042793042e-18
0.1 5.20417042793042e-18
0.01 8.67361737988404e-18
0.001 1.73472347597681e-18
0.0001 8.67361737988404e-18
1e-05 1.73472347597681e-18
1e-06 1.73472347597681e-18
1e-07 1.73472347597681e-18
1e-08 1.73472347597681e-18
1e-09 1.73472347597681e-18
1e-10 1.73472347597681e-18
};
\addplot [thick, color2, mark=*, mark size=2, mark options={solid}]
table {%
1 1.73472347597681e-18
0.1 1.73472347597681e-18
0.01 1.73472347597681e-18
0.001 1.73472347597681e-18
0.0001 2.61731004924659e-10
1e-05 1.8947370237013e-08
1e-06 4.96885583026185e-10
1e-07 9.80355982158065e-10
1e-08 6.19491541398176e-11
1e-09 1.31893888172252e-11
1e-10 1.21139036302065e-12
};
\addplot [thick, color3, mark=*, mark size=2, mark options={solid}]
table {%
1 5.20417042793042e-18
0.1 5.20417042793042e-18
0.01 5.20417042793042e-18
0.001 1.73472347597681e-18
0.0001 5.20417042793042e-18
1e-05 5.20417042793042e-18
1e-06 3.46944695195361e-18
1e-07 1.21430643318376e-17
1e-08 5.20417042793042e-18
1e-09 1.73472347597681e-18
1e-10 1.73472347597681e-18
};
\addplot [thick, color4, mark=*, mark size=2, mark options={solid}]
table {%
1 5.20417042793042e-18
0.1 5.20417042793042e-18
0.01 5.20417042793042e-18
0.001 5.20417042793042e-18
0.0001 1e-20
1e-05 1e-20
1e-06 1.15692786947785e-10
1e-07 4.20281431184888e-11
1e-08 4.87649330638273e-12
1e-09 7.37661667860046e-13
1e-10 3.20557816402278e-13
};
\addplot [thick, color5, mark=*, mark size=2, mark options={solid}]
table {%
1 6.93889390390723e-18
0.1 6.93889390390723e-18
0.01 6.93889390390723e-18
0.001 6.93889390390723e-18
0.0001 8.67361737988404e-18
1e-05 8.67361737988404e-18
1e-06 1.73472347597681e-18
1e-07 1.21430643318376e-17
1e-08 1.21430643318376e-17
1e-09 1.04083408558608e-17
1e-10 5.20417042793042e-18
};
\end{axis}

\end{tikzpicture}

%% file: figures/epsilon_moments_1.tex
\begin{tikzpicture}

\definecolor{color0}{rgb}{0.388235294117647,0.474509803921569,0.223529411764706}
\definecolor{color1}{rgb}{0.549019607843137,0.427450980392157,0.192156862745098}
\definecolor{color2}{rgb}{0.192156862745098,0.509803921568627,0.741176470588235}
\definecolor{color3}{rgb}{0.901960784313726,0.333333333333333,0.0509803921568627}
\definecolor{color4}{rgb}{0.458823529411765,0.419607843137255,0.694117647058824}
\definecolor{color5}{rgb}{0.482352941176471,0.254901960784314,0.450980392156863}

\begin{axis}[
log basis x={10},
log basis y={10},
tick align=outside,
tick pos=left,
x grid style={white!69.0196078431373!black},
xlabel={\(\displaystyle  \epsilon \)},
xmajorgrids,
xmin=3.16227766016838e-11, xmax=3.16227766016838,
xmode=log,
xtick style={color=black},
y grid style={white!69.0196078431373!black},
ylabel={\(\displaystyle \mathcal{M}_1\)},
ymajorgrids,
ymin=3.35001838708236e-21, ymax=1.66856372045887e-05,
ymode=log,
ytick style={color=black}
]
\addplot [thick, black, dashed]
table {%
1 2.22044604925031e-16
0.1 2.22044604925031e-16
0.01 2.22044604925031e-16
0.001 2.22044604925031e-16
0.0001 2.22044604925031e-16
1e-05 2.22044604925031e-16
1e-06 2.22044604925031e-16
1e-07 2.22044604925031e-16
1e-08 2.22044604925031e-16
1e-09 2.22044604925031e-16
1e-10 2.22044604925031e-16
};
\addplot [thick, color0, mark=*, mark size=2, mark options={solid}]
table {%
1 1.44096954211131e-17
0.1 1.44096954211131e-17
0.01 5.20417042793042e-18
0.001 3.22722585222632e-06
0.0001 2.31296665516625e-06
1e-05 2.08129120290887e-07
1e-06 1.65391481453731e-08
1e-07 1.51022850960816e-09
1e-08 9.63576751563239e-11
1e-09 1.16670973620469e-11
1e-10 1.00742713133749e-12
};
\addplot [thick, color1, mark=*, mark size=2, mark options={solid}]
table {%
1 5.20417042793042e-18
0.1 5.20417042793042e-18
0.01 3.00462919747432e-18
0.001 1.21430643318376e-17
0.0001 1.11076544356722e-17
1e-05 5.20417042793042e-18
1e-06 1.73205080756888e-20
1e-07 1.73478112109426e-18
1e-08 1.73478112109426e-18
1e-09 1.73205080756888e-20
1e-10 1.73478112109426e-18
};
\addplot [thick, color2, mark=*, mark size=2, mark options={solid}]
table {%
1 1.93947980722443e-17
0.1 1.93947980722443e-17
0.01 1.93947980722443e-17
0.001 3.87897250448173e-18
0.0001 4.53331417474888e-10
1e-05 3.2817807935348e-08
1e-06 8.60631096382241e-10
1e-07 1.69802639463889e-09
1e-08 1.07299088465338e-10
1e-09 2.28446995645594e-11
1e-10 2.09819366272381e-12
};
\addplot [thick, color3, mark=*, mark size=2, mark options={solid}]
table {%
1 7.94950163848063e-18
0.1 7.94950163848063e-18
0.01 7.94950163848063e-18
0.001 1.40929601402592e-17
0.0001 4.90654912385889e-18
1e-05 1.11076499342709e-17
1e-06 1.53206629090322e-17
1e-07 1.97026712507607e-17
1e-08 7.152448122691e-18
1e-09 3.00462919747432e-18
1e-10 1.73478112109426e-18
};
\addplot [thick, color4, mark=*, mark size=2, mark options={solid}]
table {%
1 7.15245511330102e-18
0.1 7.15245511330102e-18
0.01 7.15245511330102e-18
0.001 7.15245511330102e-18
0.0001 3.87897250448173e-18
1e-05 3.87897250448173e-18
1e-06 2.00385779053547e-10
1e-07 7.27948882428866e-11
1e-08 8.44635019411432e-12
1e-09 1.27767650142661e-12
1e-10 5.5520539855261e-13
};
\addplot [thick, color5, mark=*, mark size=2, mark options={solid}]
table {%
1 1.39857877853494e-17
0.1 1.39857877853494e-17
0.01 1.39857877853494e-17
0.001 1.39857877853494e-17
0.0001 1.76907777093865e-17
1e-05 1.76907777093865e-17
1e-06 6.49074090789879e-18
1e-07 1.30968750346023e-17
1e-08 1.2266347333467e-17
1e-09 6.49074090789879e-18
1e-10 4.24918736097037e-18
};
\end{axis}

\end{tikzpicture}

%% file: figures/twolevels-cons-3_adapt_normi.tex
\begin{tikzpicture}

\definecolor{color0}{rgb}{0.388235294117647,0.474509803921569,0.223529411764706}
\definecolor{color1}{rgb}{0.549019607843137,0.427450980392157,0.192156862745098}
\definecolor{color2}{rgb}{0.192156862745098,0.509803921568627,0.741176470588235}
\definecolor{color3}{rgb}{0.901960784313726,0.333333333333333,0.0509803921568627}
\definecolor{color4}{rgb}{0.458823529411765,0.419607843137255,0.694117647058824}
\definecolor{color5}{rgb}{0.482352941176471,0.254901960784314,0.450980392156863}

\begin{axis}[
log basis x={10},
log basis y={10},
tick align=outside,
tick pos=left,
x grid style={white!69.0196078431373!black},
xlabel={\(\displaystyle h_{f}\)},
xmajorgrids,
xmin=0.00196279596454476, xmax=0.0276409119776713,
xmode=log,
xtick style={color=black},
y grid style={white!69.0196078431373!black},
ylabel={\(\displaystyle E_{\infty}\)},
ymajorgrids,
ymin=0.000154330338435815, ymax=2.2271826015378,
ymode=log,
ytick style={color=black}
]
\addplot [thick, color0, mark=*, mark size=2, mark options={solid}]
table {%
0.0208333320915699 1.44110798835754
0.0104166707023978 0.789189696311951
0.00520833162590861 0.403477519750595
0.00260416627861559 0.202857971191406
};
\addplot [thick, color1, mark=*, mark size=2, mark options={solid}]
table {%
0.0208333320915699 0.989709436893463
0.0104166707023978 0.538349628448486
0.00520833162590861 0.299998760223389
0.00260416627861559 0.153937295079231
};
\addplot [thick, black, dashed]
table {%
0.0245098024606705 0.878754496574402
0.00221354188397527 0.0793625339865685
};
\addplot [thick, color2, mark=*, mark size=2, mark options={solid}]
table {%
0.0208333320915699 0.232036203145981
0.0104166707023978 0.0311903655529022
0.00520833162590861 0.00485720951110125
0.00260416627861559 0.000640325830318034
};
\addplot [thick, black, dashed]
table {%
0.0245098024606705 0.323792397975922
0.00221354188397527 0.000238512177020311
};
\addplot [thick, color3, mark=*, mark size=2, mark options={solid}]
table {%
0.0208333320915699 0.208594813942909
0.0104166707023978 0.0323955789208412
0.00520833162590861 0.00489996466785669
0.00260416627861559 0.000641702790744603
};
\addplot [thick, color4, mark=*, mark size=2, mark options={solid}]
table {%
0.0208333320915699 0.0858029946684837
0.0104166707023978 0.0311903655529022
0.00520833162590861 0.00485720951110125
0.00260416627861559 0.000640325830318034
};
\addplot [thick, color5, mark=*, mark size=2, mark options={solid}]
table {%
0.0208333320915699 0.0899876430630684
0.0104166707023978 0.0312640890479088
0.00520833162590861 0.00485789449885488
0.00260416627861559 0.000640331418253481
};
\draw (axis cs:0.005208333,0.186735348393574) ++(2pt,-2pt) node[
  scale=0.833333333333333,
  anchor=north west,
  text=black,
  rotate=11.6
]{$p = 1.00$};
\draw (axis cs:0.005208333,0.0031070165348227) ++(2pt,-2pt) node[
  scale=0.833333333333333,
  anchor=north west,
  text=black,
  rotate=31.7
]{$p = 3.00$};
\end{axis}

\end{tikzpicture}

%% file: figures/twolevels-cons-5_adapt_normi.tex
\begin{tikzpicture}

\definecolor{color0}{rgb}{0.388235294117647,0.474509803921569,0.223529411764706}
\definecolor{color1}{rgb}{0.549019607843137,0.427450980392157,0.192156862745098}
\definecolor{color2}{rgb}{0.192156862745098,0.509803921568627,0.741176470588235}
\definecolor{color3}{rgb}{0.901960784313726,0.333333333333333,0.0509803921568627}
\definecolor{color4}{rgb}{0.458823529411765,0.419607843137255,0.694117647058824}
\definecolor{color5}{rgb}{0.482352941176471,0.254901960784314,0.450980392156863}

\begin{axis}[
log basis x={10},
log basis y={10},
tick align=outside,
tick pos=left,
x grid style={white!69.0196078431373!black},
xlabel={\(\displaystyle h_{f}\)},
xmajorgrids,
xmin=0.00196279596454476, xmax=0.0276409119776713,
xmode=log,
xtick style={color=black},
y grid style={white!69.0196078431373!black},
ylabel={\(\displaystyle E_{\infty}\)},
ymajorgrids,
ymin=2.48786829861969e-07, ymax=2.95399356385824,
ymode=log,
ytick style={color=black}
]
\addplot [thick, color0, mark=*, mark size=2, mark options={solid}]
table {%
0.0208333320915699 1.40876173973083
0.0208333320915699 1.40876173973083
0.0104166707023978 0.786880612373352
0.0104166707023978 0.786880612373352
0.00520833162590861 0.404545783996582
0.00520833162590861 0.404545783996582
0.00260416627861559 0.203983202576637
0.00260416627861559 0.203983202576637
};
\addplot [thick, color1, mark=*, mark size=2, mark options={solid}]
table {%
0.0208333320915699 0.869587481021881
0.0208333320915699 0.869587481021881
0.0104166707023978 0.465701997280121
0.0104166707023978 0.465701997280121
0.00520833162590861 0.236700743436813
0.00520833162590861 0.236700743436813
0.00260416627861559 0.122459083795547
0.00260416627861559 0.122459083795547
};
\addplot [thick, black, dashed]
table {%
0.0245098024606705 0.699060440063477
0.00221354188397527 0.0631339102983475
};
\addplot [thick, color2, mark=*, mark size=2, mark options={solid}]
table {%
0.0208333320915699 0.191596671938896
0.0208333320915699 0.191596671938896
0.0104166707023978 0.0247161909937859
0.0104166707023978 0.0247161909937859
0.00520833162590861 0.00300077977590263
0.00520833162590861 0.00300077977590263
0.00260416627861559 0.000371026631910354
0.00260416627861559 0.000371026631910354
};
\addplot [thick, black, dashed]
table {%
0.0245098024606705 0.187616363167763
0.00221354188397527 0.000138202085508965
};
\addplot [thick, color3, mark=*, mark size=2, mark options={solid}]
table {%
0.0208333320915699 0.159758612513542
0.0208333320915699 0.159758612513542
0.0104166707023978 0.0156425796449184
0.0104166707023978 0.0156425796449184
0.00520833162590861 0.00188441586215049
0.00520833162590861 0.00188441586215049
0.00260416627861559 0.00024865719024092
0.00260416627861559 0.00024865719024092
};
\addplot [thick, color4, mark=*, mark size=2, mark options={solid}]
table {%
0.0208333320915699 0.0386949554085732
0.0208333320915699 0.0386949554085732
0.0104166707023978 0.0013227480230853
0.0104166707023978 0.0013227480230853
0.00520833162590861 5.73821307625622e-05
0.00520833162590861 5.73821307625622e-05
0.00260416627861559 1.93843720808218e-06
0.00260416627861559 1.93843720808218e-06
};
\addplot [thick, black, dashed]
table {%
0.0245098024606705 0.0868279039859772
0.00221354188397527 5.21673939601897e-07
};
\addplot [thick, color5, mark=*, mark size=2, mark options={solid}]
table {%
0.0208333320915699 0.033903568983078
0.0208333320915699 0.033903568983078
0.0104166707023978 0.00139467976987362
0.0104166707023978 0.00139467976987362
0.00520833162590861 5.80637461098377e-05
0.00520833162590861 5.80637461098377e-05
0.00260416627861559 1.94402787201398e-06
0.00260416627861559 1.94402787201398e-06
};
\draw (axis cs:0.005208333,0.148550368899958) ++(2pt,-2pt) node[
  scale=0.833333333333333,
  anchor=north west,
  text=black,
  rotate=6.9
]{$p = 1.00$};
\draw (axis cs:0.005208333,0.00180031151599631) ++(2pt,-2pt) node[
  scale=0.833333333333333,
  anchor=north west,
  text=black,
  rotate=19.9
]{$p = 3.00$};
\draw (axis cs:0.005208333,3.76230704083013e-05) ++(2pt,-2pt) node[
  scale=0.833333333333333,
  anchor=north west,
  text=black,
  rotate=31.2
]{$p = 5.00$};
\end{axis}

\end{tikzpicture}

%% file: figures/twolevels-diff-2_adapt_normi.tex
\begin{tikzpicture}

\definecolor{color0}{rgb}{0.388235294117647,0.474509803921569,0.223529411764706}
\definecolor{color1}{rgb}{0.549019607843137,0.427450980392157,0.192156862745098}
\definecolor{color2}{rgb}{0.192156862745098,0.509803921568627,0.741176470588235}
\definecolor{color3}{rgb}{0.901960784313726,0.333333333333333,0.0509803921568627}
\definecolor{color4}{rgb}{0.458823529411765,0.419607843137255,0.694117647058824}
\definecolor{color5}{rgb}{0.482352941176471,0.254901960784314,0.450980392156863}

\begin{axis}[
log basis x={10},
log basis y={10},
tick align=outside,
tick pos=left,
x grid style={white!69.0196078431373!black},
xlabel={\(\displaystyle h_{f}\)},
xmajorgrids,
xmin=0.00196279596454476, xmax=0.0276409119776713,
xmode=log,
xtick style={color=black},
y grid style={white!69.0196078431373!black},
ylabel={\(\displaystyle E_{\infty}\)},
ymajorgrids,
ymin=0.0737334839516918, ymax=669.873361488921,
ymode=log,
ytick style={color=black}
]
\addplot [thick, color0, mark=*, mark size=2, mark options={solid}]
table {%
0.0208333320915699 179.284072875977
0.0104166707023978 251.235305786133
0.00520833162590861 273.422546386719
0.00260416627861559 279.272216796875
};
\addplot [thick, color1, mark=*, mark size=2, mark options={solid}]
table {%
0.0208333320915699 421.03515625
0.0104166707023978 437.86767578125
0.00520833162590861 441.720275878906
0.00260416627861559 442.658508300781
};
\addplot [thick, color2, mark=*, mark size=2, mark options={solid}]
table {%
0.0208333320915699 26.1998519897461
0.0104166707023978 8.12287330627441
0.00520833162590861 2.35493803024292
0.00260416627861559 0.629950642585754
};
\addplot [thick, color3, mark=*, mark size=2, mark options={solid}]
table {%
0.0208333320915699 28.4905738830566
0.0104166707023978 11.6087503433228
0.00520833162590861 3.27941703796387
0.00260416627861559 0.845310270786285
};
\addplot [thick, color4, mark=*, mark size=2, mark options={solid}]
table {%
0.0208333320915699 7.47031736373901
0.0104166707023978 3.43151950836182
0.00520833162590861 0.984661817550659
0.00260416627861559 0.254623264074326
};
\addplot [thick, black, dashed]
table {%
0.0245098024606705 13.6802101135254
0.00221354188397527 0.11158063262701
};
\addplot [thick, color5, mark=*, mark size=2, mark options={solid}]
table {%
0.0208333320915699 7.33367872238159
0.0104166707023978 3.42303204536438
0.00520833162590861 0.984479665756226
0.00260416627861559 0.25462019443512
};
\draw (axis cs:0.005208333,0.6177471152939) ++(2pt,-2pt) node[
  scale=0.833333333333333,
  anchor=north west,
  text=black,
  rotate=23.4
]{$p = 2.00$};
\end{axis}

\end{tikzpicture}

%% file: figures/twolevels-diff-4_adapt_normi.tex
\begin{tikzpicture}

\definecolor{color0}{rgb}{0.388235294117647,0.474509803921569,0.223529411764706}
\definecolor{color1}{rgb}{0.549019607843137,0.427450980392157,0.192156862745098}
\definecolor{color2}{rgb}{0.192156862745098,0.509803921568627,0.741176470588235}
\definecolor{color3}{rgb}{0.901960784313726,0.333333333333333,0.0509803921568627}
\definecolor{color4}{rgb}{0.458823529411765,0.419607843137255,0.694117647058824}
\definecolor{color5}{rgb}{0.482352941176471,0.254901960784314,0.450980392156863}

\begin{axis}[
log basis x={10},
log basis y={10},
tick align=outside,
tick pos=left,
x grid style={white!69.0196078431373!black},
xlabel={\(\displaystyle h_{f}\)},
xmajorgrids,
xmin=0.00196279596454476, xmax=0.0276409119776713,
xmode=log,
xtick style={color=black},
y grid style={white!69.0196078431373!black},
ylabel={\(\displaystyle E_{\infty}\)},
ymajorgrids,
ymin=0.000540990163379037, ymax=1108.68298989857,
ymode=log,
ytick style={color=black}
]
\addplot [thick, color0, mark=*, mark size=2, mark options={solid}]
table {%
0.0208333320915699 242.334045410156
0.0208333320915699 242.334045410156
0.0104166707023978 336.194305419922
0.0104166707023978 336.194305419922
0.00520833162590861 364.896575927734
0.00520833162590861 364.896575927734
0.00260416627861559 372.448028564453
0.00260416627861559 372.448028564453
};
\addplot [thick, color1, mark=*, mark size=2, mark options={solid}]
table {%
0.0208333320915699 547.457763671875
0.0208333320915699 547.457763671875
0.0104166707023978 567.095764160156
0.0104166707023978 567.095764160156
0.00520833162590861 571.590270996094
0.00520833162590861 571.590270996094
0.00260416627861559 572.684875488281
0.00260416627861559 572.684875488281
};
\addplot [thick, color2, mark=*, mark size=2, mark options={solid}]
table {%
0.0208333320915699 40.2517700195312
0.0208333320915699 40.2517700195312
0.0104166707023978 12.2541999816895
0.0104166707023978 12.2541999816895
0.00520833162590861 3.39879512786865
0.00520833162590861 3.39879512786865
0.00260416627861559 0.906432509422302
0.00260416627861559 0.906432509422302
};
\addplot [thick, black, dashed]
table {%
0.0245098024606705 48.7001342773438
0.00221354188397527 0.397215515375137
};
\addplot [thick, color3, mark=*, mark size=2, mark options={solid}]
table {%
0.0208333320915699 43.0493812561035
0.0208333320915699 43.0493812561035
0.0104166707023978 16.6352977752686
0.0104166707023978 16.6352977752686
0.00520833162590861 4.64636993408203
0.00520833162590861 4.64636993408203
0.00260416627861559 1.19439101219177
0.00260416627861559 1.19439101219177
};
\addplot [thick, color4, mark=*, mark size=2, mark options={solid}]
table {%
0.0208333320915699 6.76992797851562
0.0208333320915699 6.76992797851562
0.0104166707023978 0.699062943458557
0.0104166707023978 0.699062943458557
0.00520833162590861 0.0509313829243183
0.00520833162590861 0.0509313829243183
0.00260416627861559 0.00330790248699486
0.00260416627861559 0.00330790248699486
};
\addplot [thick, black, dashed]
table {%
0.0245098024606705 15.7430572509766
0.00221354188397527 0.001047323923558
};
\addplot [thick, color5, mark=*, mark size=2, mark options={solid}]
table {%
0.0208333320915699 7.33755874633789
0.0208333320915699 7.33755874633789
0.0104166707023978 0.687022805213928
0.0104166707023978 0.687022805213928
0.00520833162590861 0.0479819849133492
0.00520833162590861 0.0479819849133492
0.00260416627861559 0.00316848838701844
0.00260416627861559 0.00316848838701844
};
\draw (axis cs:0.005208333,2.19911580699155) ++(2pt,-2pt) node[
  scale=0.833333333333333,
  anchor=north west,
  text=black,
  rotate=15.2
]{$p = 2.00$};
\draw (axis cs:0.005208333,0.0321014790632511) ++(2pt,-2pt) node[
  scale=0.833333333333333,
  anchor=north west,
  text=black,
  rotate=28.5
]{$p = 4.00$};
\end{axis}

\end{tikzpicture}

%% file: wavelet_5_results.tex
\section{Convergence analysis for a linear advection equation}
\label{sec:results_fd}
Having validated the correct implementation of the grid adaptation for non-lifted and lifted wavelets, as well as the ghost point reconstruction and finite-difference operators, we focus here on the \revtwo{behavior} of grid adaptation during the evolution of a partial differential equation. We consider the transport of a scalar field in a divergence-free flow field $\nabla \cdot \bu =0$ as a simple case of a hyperbolic conservation law:
\be
\dfrac{\partial \phi}{\partial t} + \nabla \cdot \left( \bu \phi \right) = 0
\label{eq:transport_eq}
\eed

Through this section, we compute the right-hand side using the third-order finite-difference scheme \cons{3}, perform time integration using a third-order RK3-TVD scheme \citep{Gottlieb:1998,Gottlieb:2001} (detailed in \sect{sec_rk3}), and fix the block size to $24^3$.

\revtwo{We note that this problem poses a sufficiently challenging test case to allow us to analyze our methodology and differentiate between the wavelets. Nevertheless, as shown above, the presented software framework in its current form is also able to handle time-dependent problems involving diffusion and reaction terms, as well as vector-based quantities.}

\subsection{Translation of a Gaussian blob}
\label{sec:translation_blob}
\begin{figure}[ht!]
\begin{minipage}{0.19\textwidth}
\centering
\includegraphics[width=\textwidth]{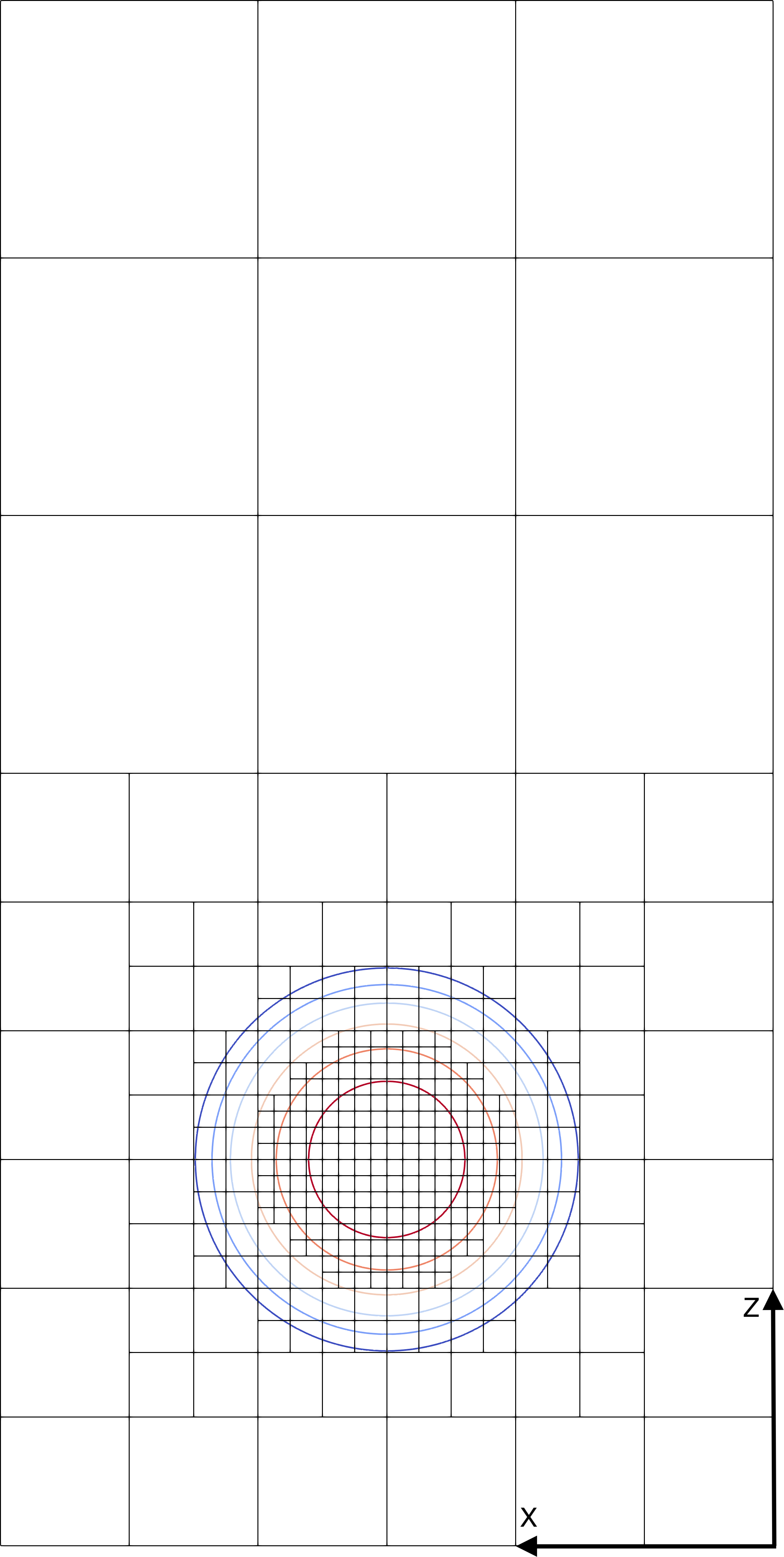}
\subcaption{$t = 0$} 
\end{minipage}%
\hfill%
\begin{minipage}{0.19\textwidth}
\centering
\includegraphics[width=\textwidth]{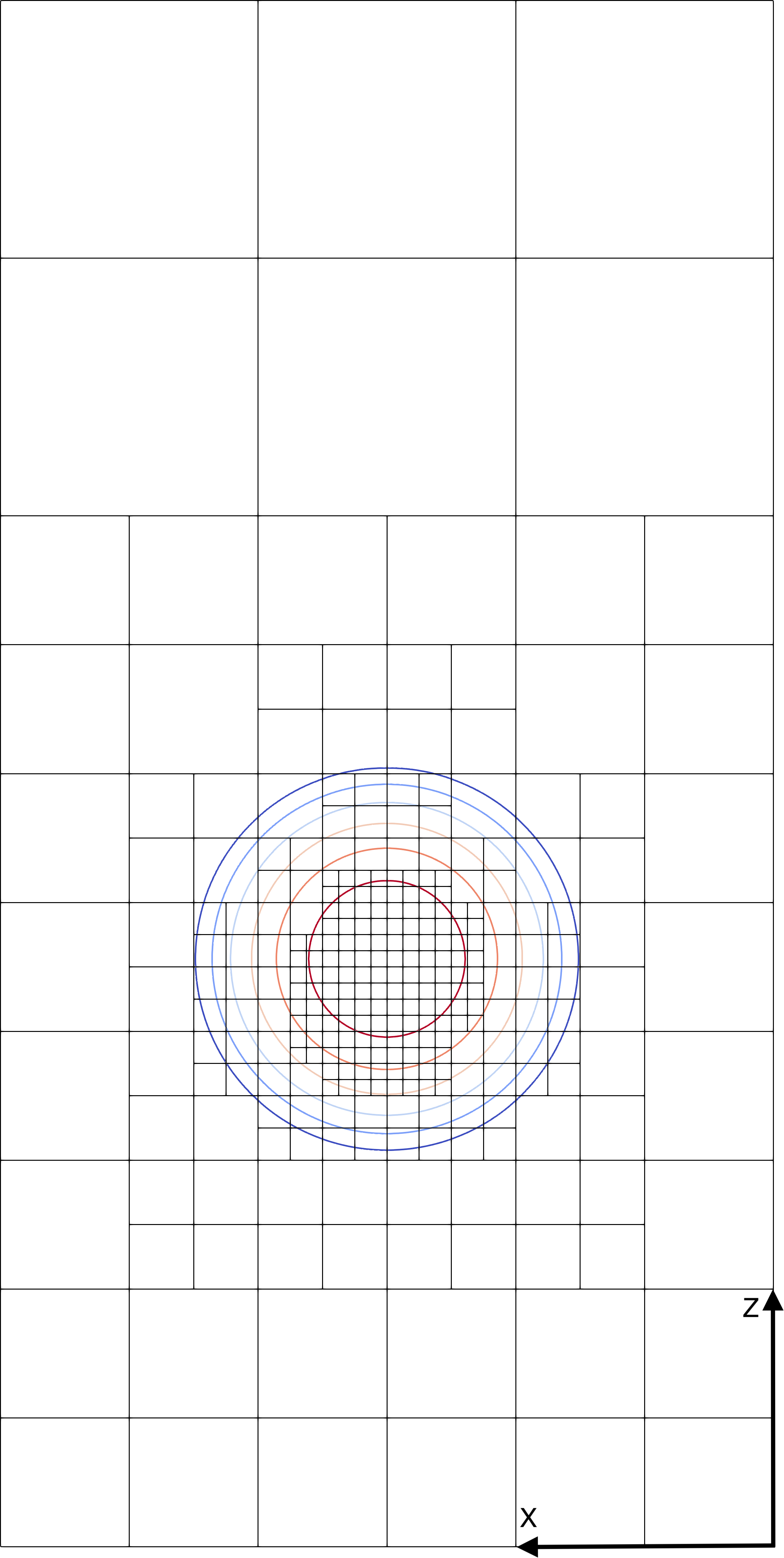}
\subcaption{$t = 0.78$} 
\end{minipage}%
\hfill%
\begin{minipage}{0.19\textwidth}
\centering
\includegraphics[width=\textwidth]{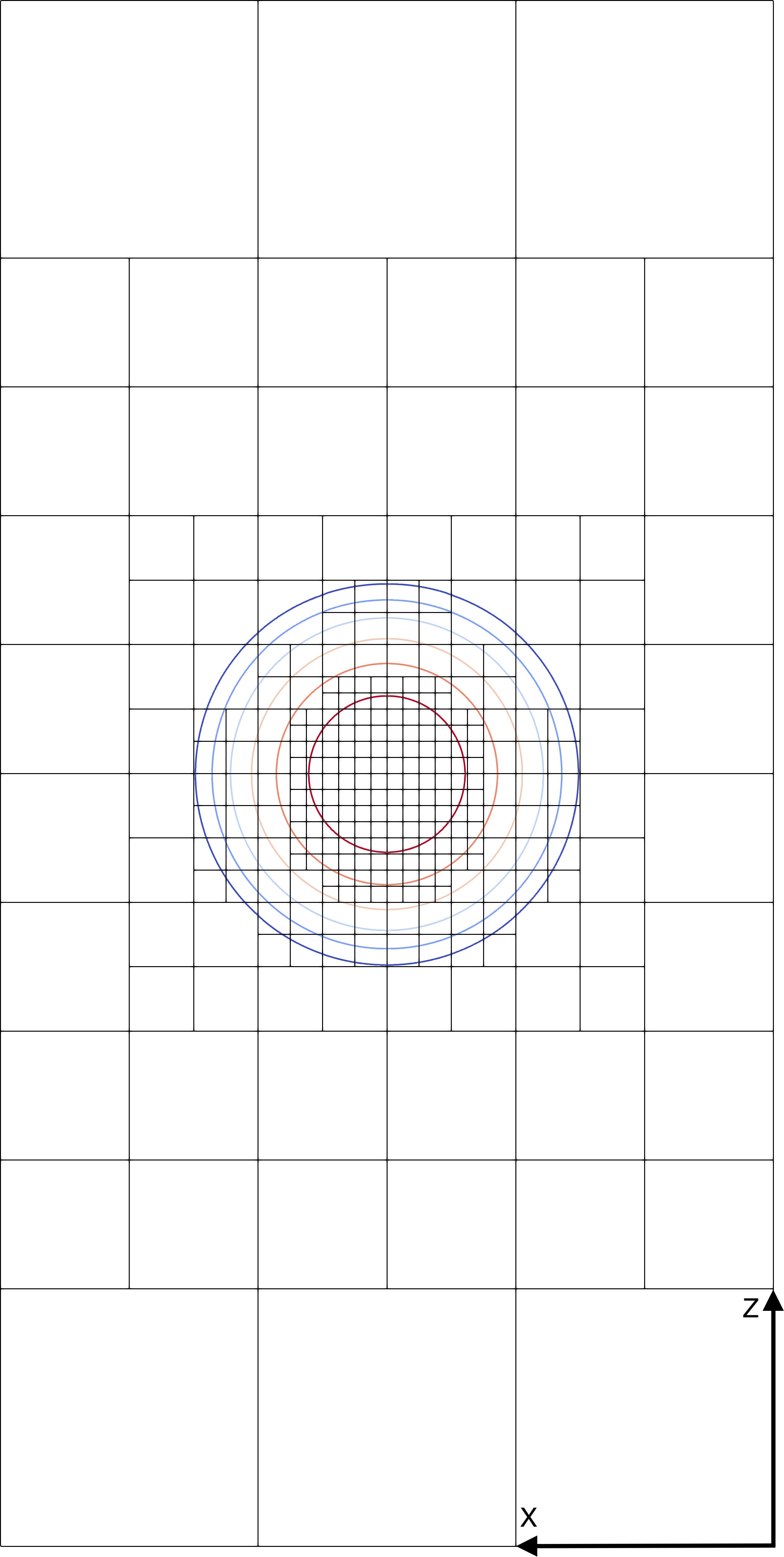}
\subcaption{$t = 1.5$} 
\end{minipage}%
\hfill%
\begin{minipage}{0.19\textwidth}
\centering
\includegraphics[width=\textwidth]{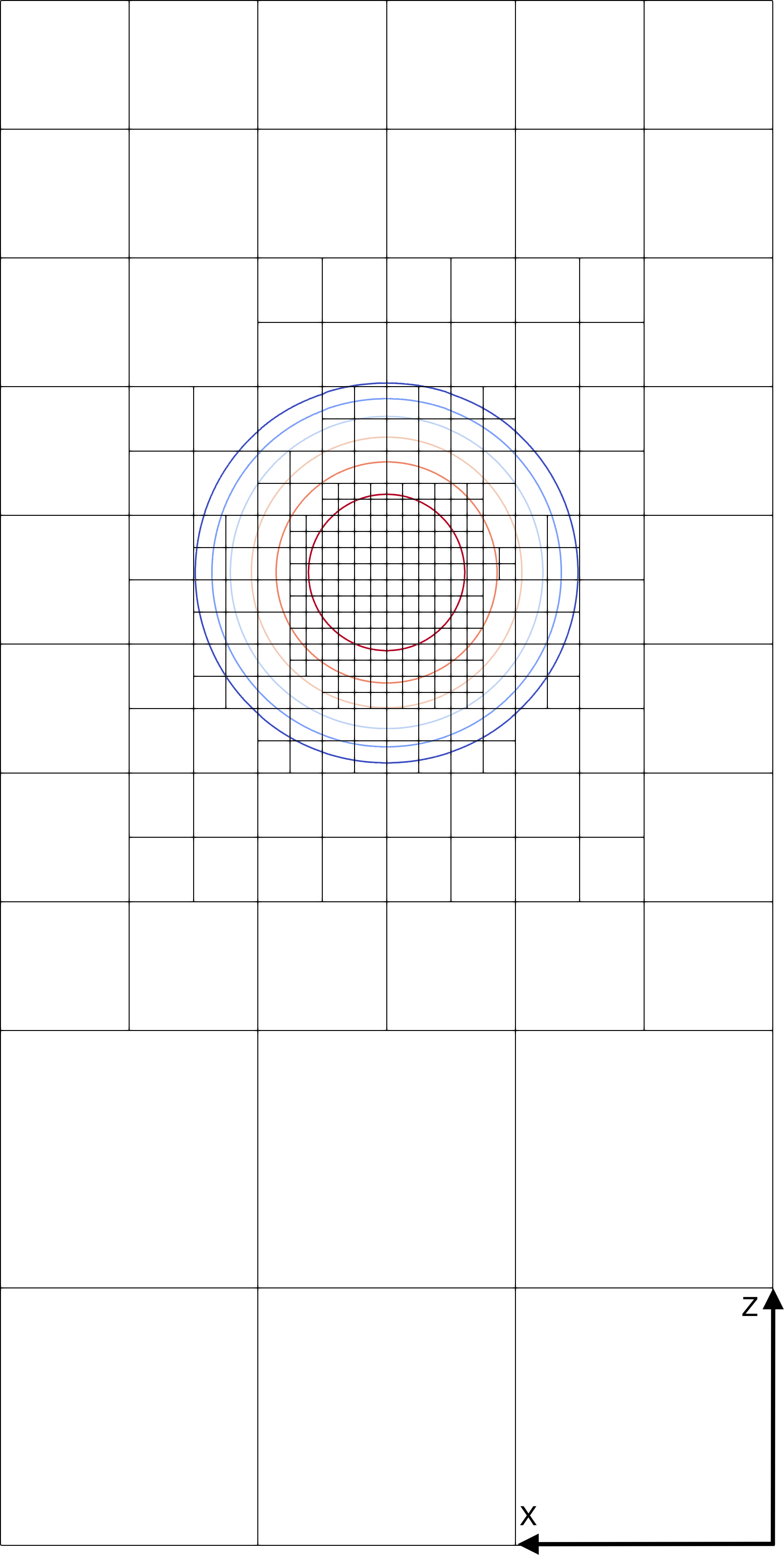}
\subcaption{$t = 2.27$} 
\end{minipage}%
\hfill%
\begin{minipage}{0.19\textwidth}
\centering
\includegraphics[width=\textwidth]{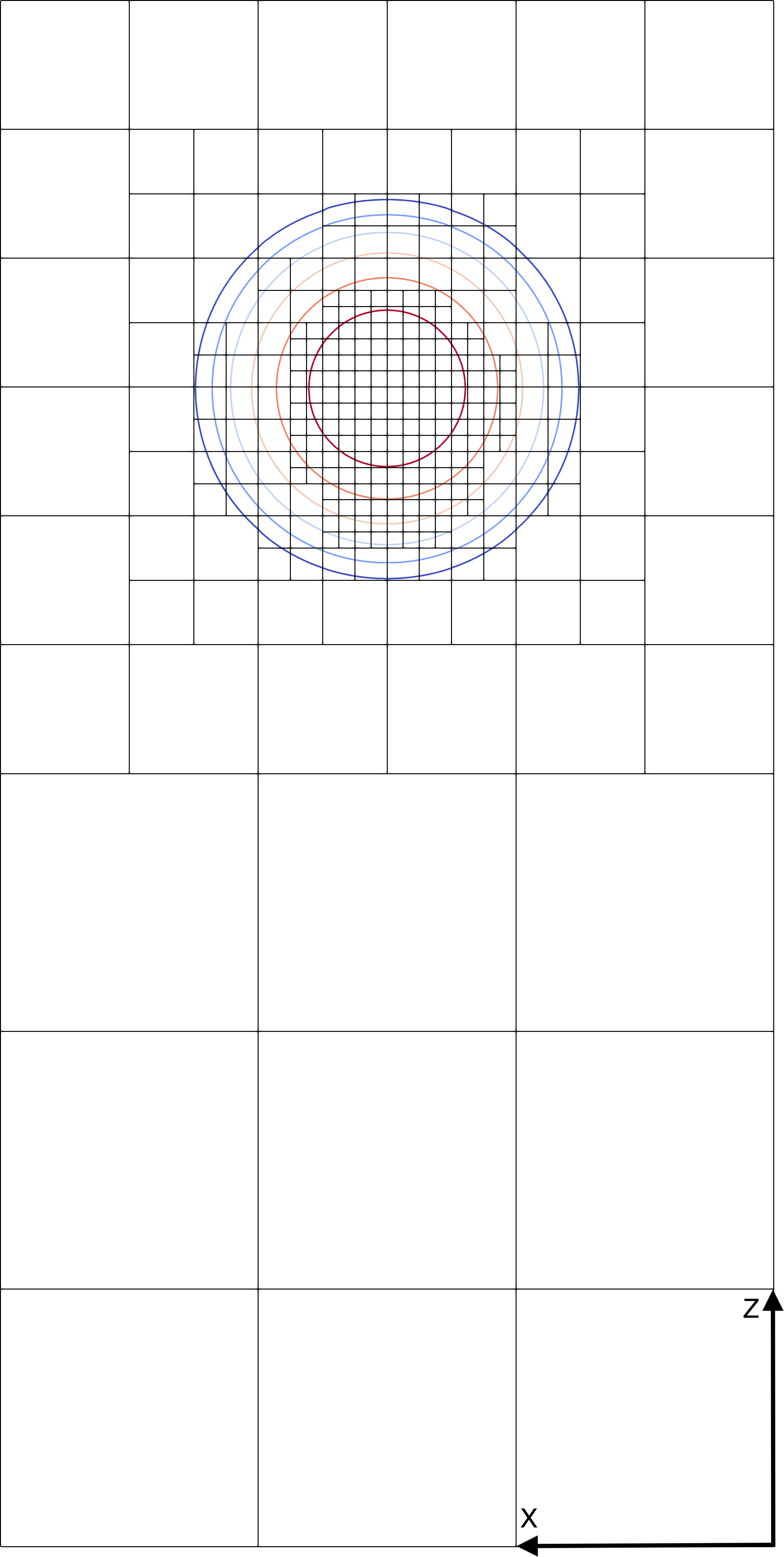}
\subcaption{$t = 3$} 
\end{minipage}%
\caption{\revone{2D Projection in the $XZ$-plane of the gaussian blob for \wave{4}{2} and $\epsilon_r  = 10^{-6}$ with $\epsilon_r / \epsilon_c = 100$. The 2D projections are illustrated with isolevels at $\left[ 10^{-6},10^{-5},10^{-4},10^{-3},10^{-2},10^{-1} \right]$.}}
 \label{fig:new_enright}
\end{figure}

To assess the convergence behavior of our algorithm and implementation, we consider a simple case of the advection of a Gaussian blob in a uniform velocity field. The computational domain is chosen as a rectangular box of size $3 \times 3 \times 6$ \revone{with each unit cube represented by a separate tree, leading to 54 trees in the domain.} We set the velocity as $\left[ 0 \;;\; 0 \;;\; 1 \right]$, and advect a Gaussian blob ($\sigma = 1/5$) initially centered at $(3/2 , 3/2 , 3/2)$ over a distance of $3$, so that we can evaluate the exact solution as a mirror of the initial condition. The time-step is controlled by setting the $\cfl = 1/4$ (based on the finest-level grid spacing), which is small enough so that the spatial discretization errors dominate the time integration errors. \revtwo{We} adapt the grid every $6$ time steps, \revtwo{so that the information travels} at most $1/16$th of the finest-level block between adaptation steps. Within this setting we vary $\epsilon_r$ and $\epsilon_c$ to control the grid adaptation during the evolution of the PDE, focusing on \wave{4}{0} and \wave{4}{2} only.
\revone{For context of the discussion, in \fref{fig:new_enright} we illustrate the obtained grid for the case of \wave{4}{2}, with $\epsilon_r=10^{-6}$ and $\epsilon_c=10^{-8}$. In this case the maximum level during the simulation is 4, leading to an effective grid spacing of $h_f = 2.6\times 10^{-3}$ or, if the grid was uniformly refined to this level, a domain with about 3 trillion grid points. The 2D projections highlight the front/back asymmetry in the grid refinement which is due to the difference between the refinement and coarsening threshold in combination with the moving field: the mesh coarsening will be triggered at larger distances behind the blob than the refinement in front of the blob.}

\paragraph{Effect of refinement threshold}
We first consider the effect of varying the refinement threshold $\epsilon_r$, keeping the ratio $\epsilon_r / \epsilon_c = 100$ fixed. 

\begin{figure}[ht!]
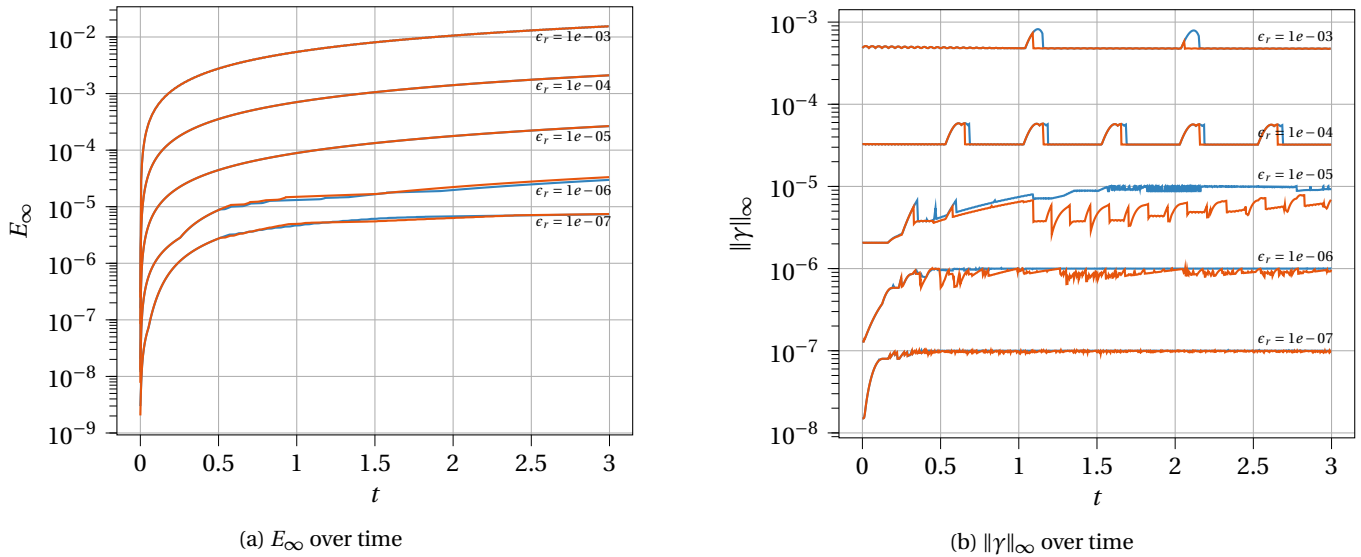

\begin{minipage}{0.49\textwidth}
	\centering
	\InputIfFileExists{figures/advection_time-w40_time_vs_erri0.25}{}{}
	\subcaption{$E_{\infty}$ over time}
	\label{fig:w40_time_vs_erri0.25}
\end{minipage}%
\hfill%
\begin{minipage}{0.49\textwidth}
	\centering
	\InputIfFileExists{figures/advection_time-w40_time_vs_detmax0.25}{}{}
	\subcaption{$\normi{\gamma}$ over time}
	\label{fig:w40_time_vs_detmax0.25}
\end{minipage}
\caption{Effect of $\epsilon_r$ (using $\epsilon_r/\epsilon_c = 100$) on the time evolution of the maximum error (left) and max detail coefficient (right) for a linear advection testcase with \cons{3}, using \wave{4}{0} (\CaptionFourZero) and \wave{4}{2} (\CaptionFourTwo).}
 \label{fig:perf_vs_acc_time_baseline_time}
\end{figure}

The time evolution of the error during the advection is presented at \fref{fig:w40_time_vs_erri0.25} for \wave{4}{0} (in blue) and \wave{4}{2} (in orange) across a range of $\epsilon_r$ values, showing that the error decreases with $\epsilon_r$ \revtwo{without significant differences} between the two wavelets. The evolution of the maximum detail coefficient on the grid ($\normi{\gamma}$) is shown in \fref{fig:w40_time_vs_detmax0.25}, which confirms that the maximum detail is always bound by $\epsilon_r$. Further, the time evolution shows that the maximum detail varies over time in a non-smooth manner, which is explained by noting that the location at which $\normi{\gamma}$ is computed can jump in space as individual blocks refine or coarsen.

\begin{figure}[ht!]
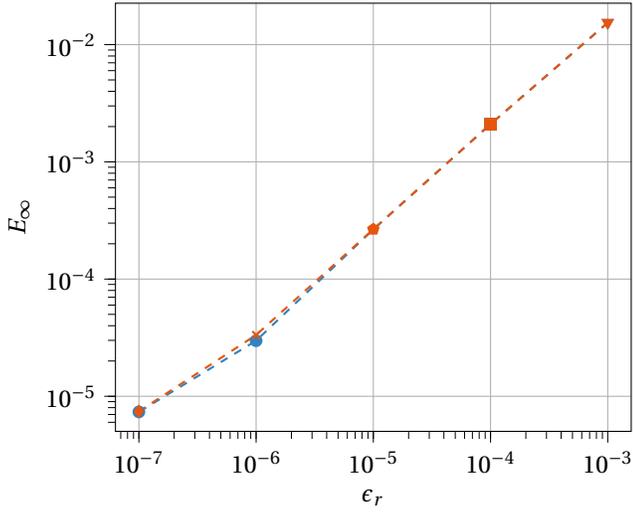
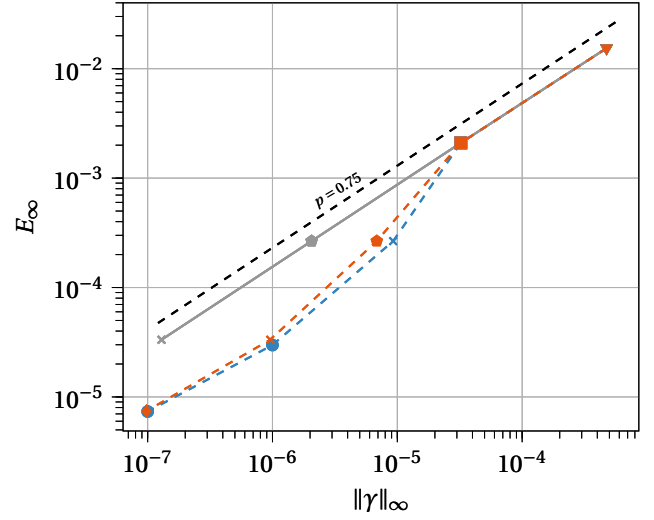
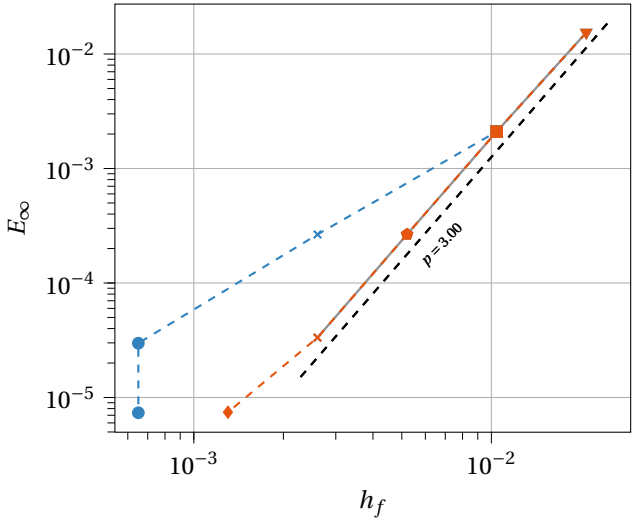
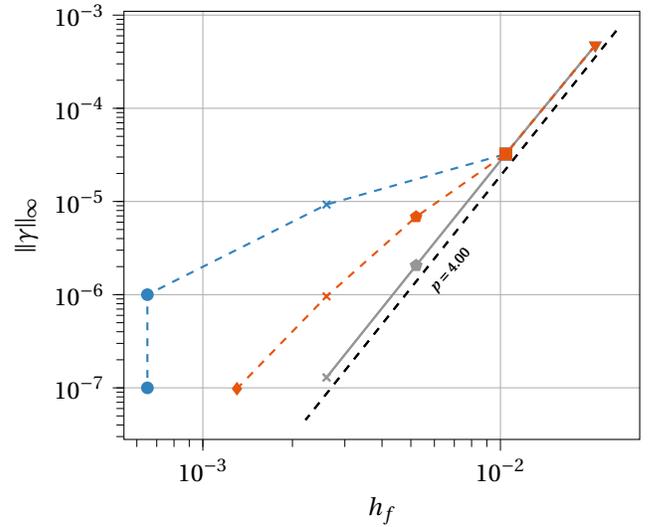

\begin{minipage}{0.49\textwidth}
	\centering
	\InputIfFileExists{figures/advection_anal-w40_eps_vs_erri3.0}{}{}	
	\subcaption{Convergence of $E_{\infty}$ as a function of $\epsilon_r$}
	\label{fig:w40_eps_vs_erri3.0}
\end{minipage}%
\hfill%
\begin{minipage}{0.49\textwidth}
	\centering
	\InputIfFileExists{figures/advection_anal-w40_det_vs_erri3.0}{}{}	
	\subcaption{Convergence of $E_{\infty}$ as a function of $\normi{\gamma}$}
	\label{fig:w40_det_vs_erri3.0}
\end{minipage}%

\begin{minipage}{0.49\textwidth}
	\centering
	\InputIfFileExists{figures/advection_anal-w40_h_vs_erri3.0}{}{}
	\subcaption{Convergence of $E_{\infty}$ as a function of the finest grid spacing $h_f$, showing third-order due to the use of \cons{3}}
	\label{fig:w40_h_vs_erri3.0}
\end{minipage}%
\hfill%
\begin{minipage}{0.49\textwidth}
	\centering
	\InputIfFileExists{figures/advection_anal-w40_h_vs_det3.0}{}{}	
	\subcaption{Variation of $h_f$ as a function of  $\normi{\gamma}$, showing fourth-order due to the use of \wave{4}{X}}
	\label{fig:w40_h_vs_det3.0}
\end{minipage}

\caption{Convergence characteristics for linear advection with \cons{3} for \wave{4}{0} (\CaptionFourZero), \wave{4}{2} (\CaptionFourTwo), and a uniform grid (\CaptionUni), as a function of  $\epsilon_r$ (using $\epsilon_r/\epsilon_c = 100$), with the error evaluated at end time $t=3$. In each plot, the subsequent data points for multiresolution simulations are obtained by systematically varying $\epsilon_r$, whereas the uniform resolution data points are obtained by systematically varying the grid spacing $h_f$. Each marker symbol is associated with a unique value of $h_f$ to facilitate the comparison across subgraphs.}
 \label{fig:w40_h}
\end{figure}

Studying the convergence behavior of a simulation on a multi-level grid is not trivial. On a uniform grid, one would show convergence evaluating the maximum error $E_{\infty}$ as a function of the grid spacing $h$. \revone{On adaptive grids}, however, the grid spacing $h$ varies in space and time, we have no direct control on the minimum grid spacing $h$, and there is no guarantee that the maximum error $E_{\infty}$ is measured at a single physical location when adapting the grid. \revone{Alternatively, one can use a parameter like the effective number of degrees-of-freedom in the simulation to measure convergence. We show associated results briefly at the end of this subsection. However, we consider this convergence metric less relevant to the point of this work, because the effective number of degrees-of-freedom is an outcome of the simulation and will likely vary in time.} Instead, we control the error primarily by varying $\epsilon_r$, and so a \revone{more suitable convergence analysis relates } the maximum error at the end time of the simulation as a function of \revone{the input parameter} $\epsilon_r$ (\fref{fig:w40_eps_vs_erri3.0}). The result demonstrates that $\epsilon_r$ is successful at controlling the error and moreover its value provides an estimate of the error made in a simulation, albeit with a problem-specific prefactor of $\mathcal{O}(10)$ in this case.

We can decompose this convergence behavior into different components. First, through our adaptation policy we guarantee that $\epsilon_r$ bounds the maximum detail coefficient $\normi{\gamma}$. In \fref{fig:w40_det_vs_erri3.0} we show the relation between the error and the maximum detail coefficient, where we included a uniform resolution line (in gray) obtained by varying the constant grid spacing $h$, and computing for each $h$ the maximum error $E_{\infty}$ as well as the maximum detail coefficient evaluated by a single-level \wave{4}{0} transform. Both uniform and multiresolution results show a clear $3/4$ slope, where the wavelet results vary as the locations of the maximum scaling coefficient and the maximum error jump independently across different locations in the grid between individual simulations. The $3/4$ slope can be explained by two observations. First, we know that the order of the \cons{3} spatial discretization is third, so that $E_{\infty} \propto h^3$. This is confirmed in \fref{fig:w40_h_vs_erri3.0}, showing the error as a function of the finest-level grid spacing $h_f$; both the uniform and the \wave{4}{2} lines follow a third-order slope (\wave{4}{0} will be discussed below). Second, according to \eqqref{eq:interp_error_bound}, the detail coefficients associated with a projection of a given smooth function onto a level with spacing $h$ scale as $\normi{\gamma} \propto h^N$, where here $N = 4$ is the polynomial interpolation order of the wavelet. For our data this convergence is confirmed in \fref{fig:w40_h_vs_det3.0}, where we plot the maximum detail coefficient as a function of the finest grid spacing $h_f$. Combining these relations we find that $E_{\infty} \propto \normi{\gamma}^{3/4}$.

Looking more closely at \frefs{fig:w40_h_vs_erri3.0}{fig:w40_h_vs_det3.0} \revtwo{indicates} that though the overall behavior in \fref{fig:w40_eps_vs_erri3.0} is consistent between \wave{4}{0} and \wave{4}{2}, the associated grid adaptation strategies are different. \Fref{fig:w40_h_vs_erri3.0} shows that \wave{4}{2} behaves similarly to the uniform resolution grid, which is impressive as the uniform grid result represents the smallest possible error for any given $h_f$; the \wave{4}{2} does not compromise that error despite the continuous grid adaptation during the simulation. The behavior of \wave{4}{0} instead demonstrates that this wavelet transform generates detail coefficients that do not predict the error committed by the PDE and therefore cause a spurious coarsening and belated refinement. This is emphasized by the last two points on the left (associated with $\epsilon_r = 10^{-6}$ and $\epsilon_r = 10^{-7}$) where the error goes down even though the finest grid spacing stays the same. This means that at $\epsilon_r = 10^{-6}$ the error was associated with a coarser level than the maximum, and refining that level without affecting the finest level successfully reduced the error. Similarly, in \fref{fig:w40_h_vs_det3.0} the gaps between the uniform grid and the adapted grids \revtwo{are} associated with coarsening, which increases the maximum detail coefficients. We see again that the behavior between \wave{4}{0} and \wave{4}{2} is different: the reduced aliasing of the lifted \wave{4}{2} leads to a better correlation of the maximum detail coefficients, the finest grid spacing, and the maximum error, compared with the non-lifted \wave{4}{0}. We emphasize however, that despite the different strategies both wavelets successfully control the error during the evolution of this PDE as a function of $\epsilon_r$, as evidenced by the overlapping lines in \fref{fig:w40_eps_vs_erri3.0}.

\begin{figure}[ht!]
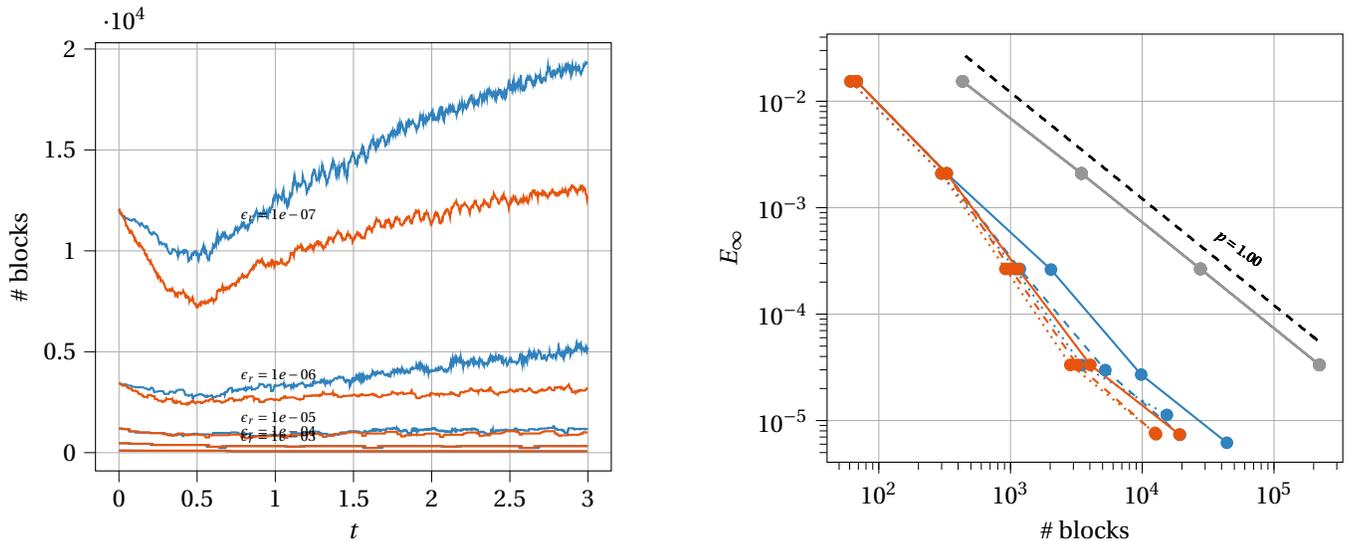

\centering
\centering
\begin{minipage}[t]{0.49\textwidth}
	\centering
	\InputIfFileExists{figures/advection_time-w40_time_vs_nblock0.25}{}{}	
	\subcaption{evolution of the number of blocks over time for various values of $\epsilon_r$, for $\epsilon_r/\epsilon_c = 100$}
	\label{fig:w40_time_vs_nblock}
\end{minipage}%
\hfill%
\begin{minipage}[t]{0.49\textwidth}
	\centering
	\InputIfFileExists{figures/advection_anal-ec_nblock_vs_erri3.0}{}{}	
	\subcaption{wavelet-based convergence: $E_{\infty}$ vs the number of blocks (at $t=3.0$), with $\epsilon_r/\epsilon_c = 16$ (\protect\ThickLine{black}{dotted}), $\epsilon_r/\epsilon_c = 100$ (\protect\ThickLine{black}{dashed}), and $\epsilon_r/\epsilon_c = 10^4$ (\protect\ThickLine{black}{solid}).}
	\label{fig:w40_nblock_vs_erri3.0}
\end{minipage}%
\caption{Evolution of the number of blocks in the grid over time (left), and maximum error against number of blocks (right) for the linear advection testcase with \cons{3} using \wave{4}{0} (\CaptionFourZero), \wave{4}{2} (\CaptionFourTwo), and a uniform grid (\CaptionUni), for various values of $\epsilon_r$'s. In (b), we also show the effect of varying $\epsilon_r/\epsilon_c$ using different line styles.}
 \label{fig:w40_nblock}
\end{figure}

The different strategies of \wave{4}{0} and \wave{4}{2} are \revtwo{further} reflected by the number of blocks required throughout the simulation, as a function of $\epsilon_r$. The evolution of the number of blocks over time is shown in \fref{fig:w40_time_vs_nblock}. For early times, as the Gaussian blob translates through the grid the \revtwo{trailing side of the blob} gets coarsened since the small details there fall below $\epsilon_c$. The \revtwo{leading side} does not get refined yet as the small details do not yet exceed $\epsilon_r$, so the number of blocks decreases, \revone{leading to the front/back asymmetric grid structure as shown in \fref{fig:new_enright}}. Once the \revtwo{leading side of the blob} gets picked up by $\epsilon_r$ the number of blocks increase\revtwo{s} again \revone{, though the asymmetry persists.}. \revtwo{At later times, the number of blocks} plateaus for \wave{4}{2} across all values of $\epsilon_r$, \revtwo{indicating that the grid structure is largely constant}. For \wave{4}{0}, on the other hand, the number of blocks increases throughout the simulation. This is consistent with the convergence analysis above, where we observed that \wave{4}{0} produces detail coefficients that do not accurately reflect the PDE error, and thus refines the grid in locations without strongly reducing the error. %

Plotting the error made in the simulation as a function of the number of blocks (\fref{fig:w40_nblock_vs_erri3.0}) shows \revone{a slope of $1$ with respect to the effective number of degrees-of-freedom for both uniform and multiresolution simulations. Since the number of blocks is inversely proportional to an `effective' linear grid spacing to the power of $3$, this thus recovers the third-order convergence of the discretization scheme. Comparing the two wavelets, we find that } \wave{4}{0} requires up to \revtwo{twice more blocks} compared to \wave{4}{2} at the lowest error values. Both wavelets provide significant gains over the uniform resolution simulation (fewer blocks by a factor of $\approx30-40$ for \wave{4}{0} and $\approx 70$ for \wave{4}{2}), though this metric is heavily dependent on the scale separation of the simulation. A rough estimate implies that the ratio of the volume occupied by a sphere of radius $3 \sigma$ and the volume of the rectangular domain is $\approx 70$, similar to the compression rate of \wave{4}{2}. 

\paragraph{Effect of coarsening threshold}
In the previous section we varied $\epsilon_r$ while keeping the ratio $\epsilon_r/\epsilon_c$ fixed. Repeating the analysis for a range of values for $\epsilon_r/\epsilon_c$ does not significantly change the results, as shown in \fref{fig:w40_nblock_vs_erri3.0}. Here the dotted, dashed, and solid lines correspond to $\epsilon_r/\epsilon_c = 16$, $\epsilon_r/\epsilon_c = 100$, and $\epsilon_r/\epsilon_c = 10^4$ respectively, and each data point for each simulation is associated with a given value of $\epsilon_r$.  For both wavelets, the number of blocks associated with a given error decreases slightly when $\epsilon_r/\epsilon_c$ decreases, reflecting the more aggressive coarsening of the grid when $\epsilon_c$ increases. At all points, except the finest $\epsilon_r$ for \wave{4}{0}, increasing $\epsilon_c$ for a given $\epsilon_r$ decreases the number of blocks without significantly changing the error. This emphasizes the capability of wavelets to detect where to compress information without degrading the overall accuracy of the solution, and shows that generally a simulation should take $\epsilon_r/\epsilon_c$ close to the lower bound of $2^N$ (explained in \sect{subsubsec:adaptation}), with $N$ the wavelet order.

\paragraph{Effect of adaptation frequency}
\revone{For a fixed  $\epsilon_r = 10^{-5}$, $\epsilon_r/\epsilon_c = 100$, and \wave{4}{2} we vary the adaptation frequencies ranging from every $6$ (as in the cases above) to every $768$ time steps. With constant value of $\cfl = 1/4$ this corresponds to the signal traveling from $1/16$th of a block to $8$ blocks at the finest scale between adaptations. As mentioned in \cref{subsubsec:adaptation}, there are two expected effects that occur when increasing the adaptation frequency: a delayed coarsening of blocks where the detail coefficient falls below the coarsening threshold, and a delayed refinement of blocks where the detail coefficient exceeds the refinement threshold. Compared to a more frequently adapted grid, the former will reduce the compression rate of the adapted grid whereas the latter will increase the error of the PDE solution.}

\begin{figure}[ht!]
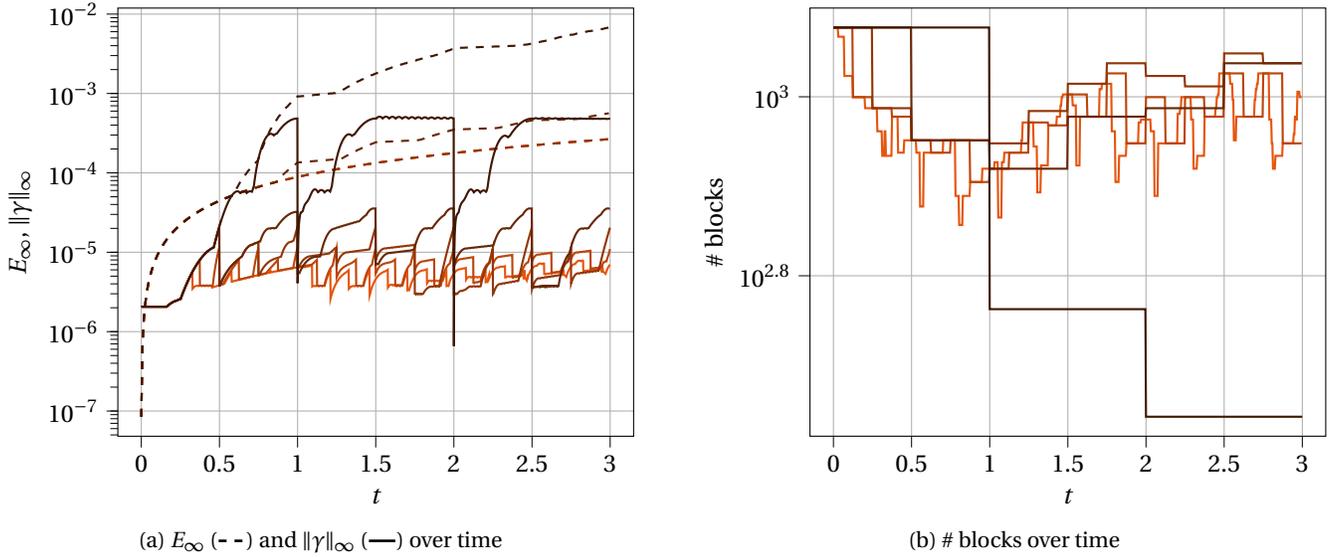

\centering
\begin{minipage}[t]{0.49\textwidth}
	\centering
	\InputIfFileExists{figures/sadv_time_err}{}{}	
	\subcaption{\revone{$E_{\infty}$ (\protect\ThickLine{black}{dashed}) and $\normi{\gamma}$ (\protect\ThickLine{black}{solid}) over time}}
	\label{fig:adapt_freq_err}
\end{minipage}%
\hfill%
\begin{minipage}[t]{0.49\textwidth}
	\centering
	\InputIfFileExists{figures/sadv_time_nblocks}{}{}
	\subcaption{\revone{\# blocks over time}}
	\label{fig:adapt_freq_blocks}
\end{minipage}%
\definecolor{black64233}{RGB}{64,23,3}
\definecolor{darkgray176}{RGB}{176,176,176}
\definecolor{firebrick1886910}{RGB}{188,69,10}
\definecolor{maroon105395}{RGB}{105,39,5}
\definecolor{orangered2308513}{RGB}{230,85,13}
\definecolor{saddlebrown147548}{RGB}{147,54,8}
\caption{\revone{Effect of adaptation frequency for the linear advection testcase with \cons{3} using \wave{4}{2}. The number of timesteps between successive adaptations shown is $6$ (\protect\ThickLine{orangered2308513}{solid}), $96$  (\protect\ThickLine{firebrick1886910}{solid}), $192$  (\protect\ThickLine{saddlebrown147548}{solid}), $384$  (\protect\ThickLine{maroon105395}{solid}), and $768$  (\protect\ThickLine{black64233}{solid}). With $\cfl=1/4$, this corresponds to the information traveling $1/16$th of a block, $1$ block, $2$ blocks, $4$ blocks, and $8$ blocks respectively.}}
 \label{fig:adapt_freq}
\end{figure}

\revone{The results of our numerical experiment are shown in \fref{fig:adapt_freq} with the error and maximum detail evolution for all cases over time (left) and the evolution of the number of blocks (right). Increasing the adaptation frequency up to and including every $192$ time steps increases the number of blocks as grid compression is delayed, but does not strongly affect the error. When the adaptation frequency is every $384$ time steps or lower, we observe that there are prolonged periods of time where the maximum detail coefficient exceeds the refinement threshold $10^{-5}$, and the error starts to increase significantly compared to the other cases. For the lowest frequency (adaptation every $768$ time steps), the maximum level in the grid drops to $2$ compared to $3$ for all other cases, as details of the signal are progressively lost and the grid is compressed accordingly. For this particular testcase, we thus find some robustness in the results to the adaptation frequency between values of $6$ and $192$, partly due to the smooth nature of the function and partly due to our 2:1 constraint in adaptation approach. Both features ensure that the grid adapted at a single time is able to capture and evolve the solution well when the signal travels up to 2 fine-level blocks afterwards.}

\subsection{Deformation of a Gaussian blob}
\label{sect_enright}
\revtwo{Here we} present the results of our framework on a more challenging scale separation problem, and showcase the ability of the wavelet-based adaptation to track the need \revtwo{for} computational resources. \revtwo{Specifically, we} use the advection equation to transport a Gaussian blob in a non-linear periodic incompressible flow field defined as
\begin{eqc}
\begin{array}{rl}
u(\x) = \sin^2(\pi \; x )  \; \sin(2 \pi \; y) \; \sin(2\pi \; z) \\
v(\x) = \sin(2 \pi \; x )  \; \sin^2(\pi \; y) \; \sin(2\pi \; z) \\
w(\x) = \sin(2 \pi \; x )  \; \sin(2 \pi \; y) \; \sin^2(\pi \; z)
\end{array}
\end{eqc}
which was proposed originally in \citep{Leveque:1996} and has been used extensively since in the level-set community since \citep{Enright:2002}. Following the latter we multiply the velocity components by $\cos(\pi t/3)$, and here evaluate only the `forward' evolution up to $t = 1.5$. As initial condition we choose a compact spherically symmetric Gaussian blob defined in terms of the radial coordinate $r$ as
\begin{eqc}
\phi_0(r) = \expo{ -\dfrac{r^2/\sigma^2}{1-r^2/\beta^2}}%
\end{eqc}
for $r < \beta$ and with $\sigma = 0.1$ and $\beta = 2 \sigma$, centered at $\left[ 0.35 \;,\;0.35 \;,\;0.35 \right]$ within a unit cube. We apply grid adaptation tolerances of $\epsilon_r=10^{-2}$ and $\epsilon_r/\epsilon_c=100$, use $\cfl = 0.5$, and adapt the grid every $12$ time steps.

\begin{figure}[ht!]
\centering
\begin{minipage}{0.33\textwidth}
	\centering
	\includegraphics[width=\textwidth]{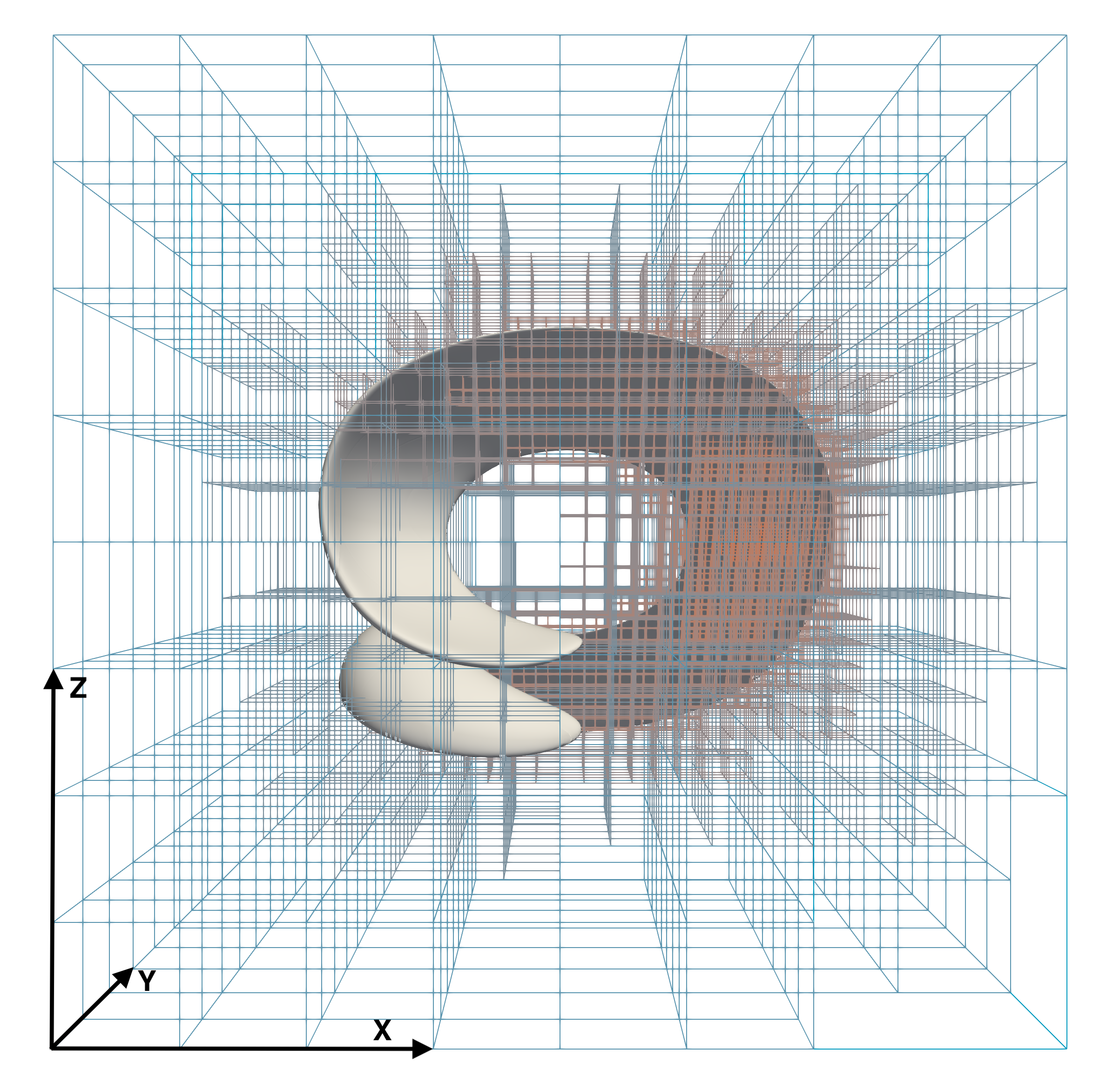}
	\subcaption{$3$D view}
\end{minipage}%
\centering
\begin{minipage}{0.33\textwidth}
	\centering
	\includegraphics[width=\textwidth]{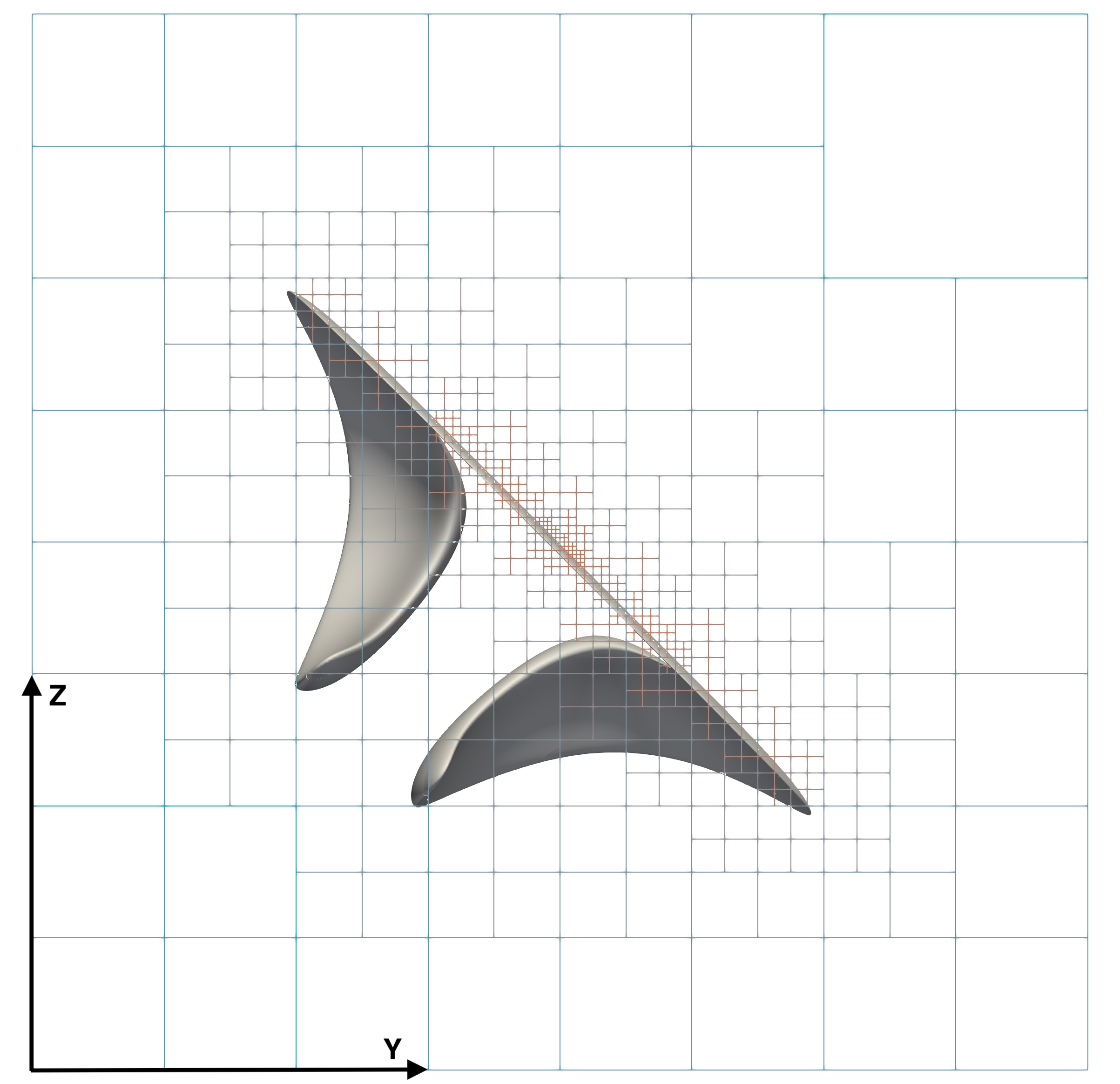}
	\subcaption{plane YZ}
\end{minipage}%
\hfill%
\begin{minipage}{0.33\textwidth}
	\centering
	\includegraphics[width=\textwidth]{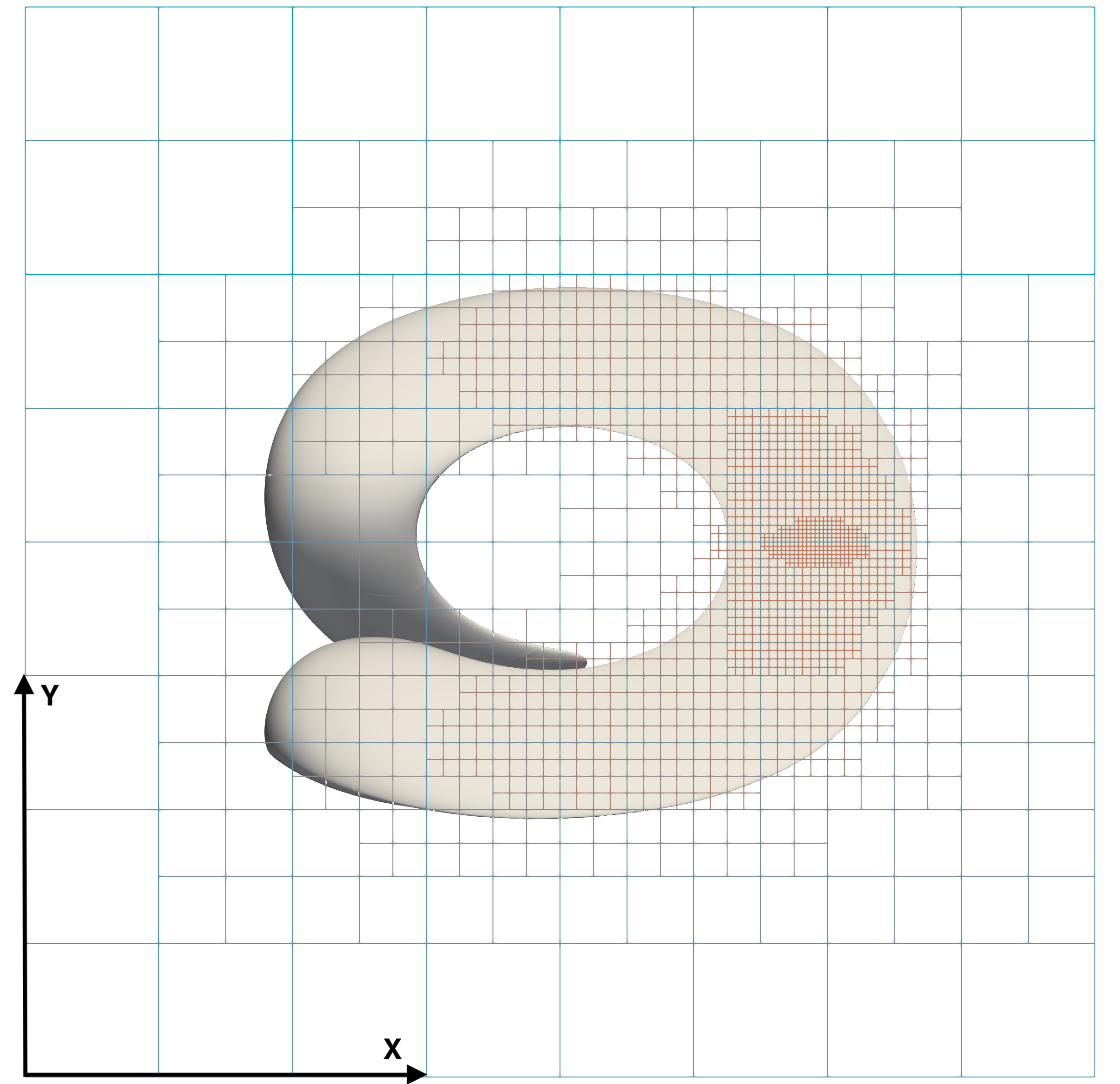}	
	\subcaption{plane XY}
\end{minipage}

\caption{Advection of a compact gaussian scalar field using a deformation velocity field: $3$D and perspective views of the grid together with an isosurface of the field at value $\phi = 0.25$, for $t=1.5$. We have visualized the outlines of individual blocks to highlight the structure of the grid.}
\label{fig:enright_3d}
\end{figure}

At $t=0$, the grid contains only 92 blocks on levels $2$ and $3$, as the Gaussian is smooth and can be captured by the fourth-order wavelets at relatively coarse resolution. As the blob deforms, it thins rapidly near the center of the domain, as shown through visualization of an isosurface in \fref{fig:enright_3d}. This process triggers refinement throughout the evolution leading eventually to \revtwo{approximately} $11,500$ blocks, and levels between $2$ all the way to $8$. The refinement pattern at the final time is interesting, as the grid reflects the subtleties in the wavelet analysis: the maximum levels are localized concentrated exactly where we expect a fourth-order polynomial interpolation to show the largest errors, near the center of the domain. \revtwo{The} `tips' of the deformed shape are still smooth and, despite relatively large gradients, still captured well on much coarser resolutions. From a user perspective, we note that we do not need to bound the maximum level; instead, the consistency of the refinement criterion with the wavelet-based grid adaptation ensures that setting a suitable $\epsilon_r$ is sufficient to keep the levels in the grid within reasonable bounds.

\revone{Finally, we study the ability of the $\epsilon_r$ criterion as an input to control the solution error in this non-linear problem. As the equations have no analytical solutions, we perform a self-convergence study where the result of a uniform-grid simulation at level $5$ (or $768^3$ grid points) at time $t=0.5$ is used as reference solution. For the multiresolution cases, we use \wave{4}{2} and consider four different refinement thresholds $\epsilon_r \in \left\{ 1 \;;\; 10^{-1} \;;\; 10^{-2} \;;\; 10^{-3} \right\}$. We set the time step as $\Delta t = 1/2 \; h_f/U$, with $h_f$ the instantaneous finest spacing in the grid, and $U$ fixed to be $U=1$ at all times.
None of these multiresolution simulations reaches level~5 within this time window. Once the multiresolution simulations has reached $t=0.5$, we refine their multi-level grids uniformly to level 5 using the chosen wavelet, and then compute the maximum error between the refined uniform grid and the reference uniform resolution simulation result.}

\revone{The evolution of the error with the refinement criterion is shown in \fref{fig:convenright_eps} where we observe that the slope matches the theoretical convergence rate of $3/4$ for both the refinement criterion $\epsilon_r$ and the measured maximum detail in the domain $\normi{\gamma}$ reasonably well. In \fref{fig:convenright_nblocks} we show the error as a function of the number of blocks at the final time. This metric is more sensitive as the number of blocks varies in time; in this particular case, the second-last point at $\epsilon_r = 10^{-2}$ just refined the grid prior to $t=0.5$ which has moved the data point to the right compared to the expected slope. Nevertheless, this convergence analysis demonstrates that the main results of the linear transport problem described in \cref{sec:translation_blob} can still be reproduced for this more challenging, non-linear transport problem as well.}

\begin{figure}[ht!]
\begin{minipage}{0.49\textwidth}
	\centering
	\InputIfFileExists{figures/convenright_eps}{}{}	
	\subcaption{$E_{\infty}$ vs $\epsilon_r$ (\protect\ThickLineSquare{col42}{col42}), and vs $\normi{\gamma}$ (\CaptionFourTwo)  with $\epsilon_r/\epsilon_c = 100$.}
	\label{fig:convenright_eps}
\end{minipage}%
\hfill%
\begin{minipage}{0.49\textwidth}
	\centering
	\InputIfFileExists{figures/convenright_blocks}{}{}	
	\subcaption{$E_{\infty}$ vs \# blocks.}
	\label{fig:convenright_nblocks}
\end{minipage}
\caption{Advection of a compact gaussian scalar field using a deformation velocity field:  wavelet-based convergence analysis.}
\label{fig:convenright}
\end{figure}

%% file: figures/advection_anal-w40_eps_vs_erri3.0.tex
\begin{tikzpicture}

\definecolor{darkgray176}{RGB}{176,176,176}
\definecolor{orangered2308513}{RGB}{230,85,13}
\definecolor{steelblue49130189}{RGB}{49,130,189}

\begin{axis}[
log basis x={10},
log basis y={10},
tick align=outside,
tick pos=left,
x grid style={darkgray176},
xlabel={\(\displaystyle \epsilon_{r}\)},
xmajorgrids,
xmin=6.30957344480193e-08, xmax=0.00158489319246111,
xmode=log,
xtick style={color=black},
y grid style={darkgray176},
ylabel={\(\displaystyle E_{\infty}\)},
ymajorgrids,
ymin=5.02155494358481e-06, ymax=0.0225682749682589,
ymode=log,
ytick style={color=black}
]
\addplot [thick, steelblue49130189, dashed]
table {%
0.001 0.0153980366666667
0.0001 0.0020983035
1e-05 0.00026559265
1e-06 2.980041e-05
1e-07 7.3598375e-06
};
\addplot [thick, steelblue49130189, mark=triangle*, mark size=2, mark options={solid,rotate=180}]
table {%
0.001 0.0153980366666667
};
\addplot [thick, steelblue49130189, mark=square*, mark size=2, mark options={solid}]
table {%
0.0001 0.0020983035
};
\addplot [thick, steelblue49130189, mark=x, mark size=2, mark options={solid}]
table {%
1e-05 0.00026559265
};
\addplot [thick, steelblue49130189, mark=*, mark size=2, mark options={solid}]
table {%
1e-06 2.980041e-05
};
\addplot [thick, steelblue49130189, mark=*, mark size=2, mark options={solid}]
table {%
1e-07 7.3598375e-06
};
\addplot [thick, orangered2308513, dashed]
table {%
0.001 0.0153981433333333
0.0001 0.002098303
1e-05 0.000265951666666667
1e-06 3.33405483333333e-05
1e-07 7.43860666666667e-06
};
\addplot [thick, orangered2308513, mark=triangle*, mark size=2, mark options={solid,rotate=180}]
table {%
0.001 0.0153981433333333
};
\addplot [thick, orangered2308513, mark=square*, mark size=2, mark options={solid}]
table {%
0.0001 0.002098303
};
\addplot [thick, orangered2308513, mark=pentagon*, mark size=2, mark options={solid}]
table {%
1e-05 0.000265951666666667
};
\addplot [thick, orangered2308513, mark=x, mark size=2, mark options={solid}]
table {%
1e-06 3.33405483333333e-05
};
\addplot [thick, orangered2308513, mark=diamond*, mark size=2, mark options={solid}]
table {%
1e-07 7.43860666666667e-06
};
\end{axis}

\end{tikzpicture}

%% file: figures/advection_anal-w40_det_vs_erri3.0.tex
\begin{tikzpicture}

\definecolor{darkgray150}{RGB}{150,150,150}
\definecolor{darkgray176}{RGB}{176,176,176}
\definecolor{orangered2308513}{RGB}{230,85,13}
\definecolor{steelblue49130189}{RGB}{49,130,189}

\begin{axis}[
log basis x={10},
log basis y={10},
tick align=outside,
tick pos=left,
x grid style={darkgray176},
xlabel={\(\displaystyle \vert\!\vert \gamma \vert\!\vert_{\infty}\)},
xmajorgrids,
xmin=6.3548222537649e-08, xmax=0.000860316218575168,
xmode=log,
xtick style={color=black},
y grid style={darkgray176},
ylabel={\(\displaystyle E_{\infty}\)},
ymajorgrids,
ymin=4.88659294284375e-06, ymax=0.0399918531607427,
ymode=log,
ytick style={color=black}
]
\addplot [thick, darkgray150]
table {%
0.00047454685 0.0153979966666667
3.22709e-05 0.00209830266666667
2.05770583333333e-06 0.000265951666666667
1.2919935e-07 3.33405483333333e-05
};
\addplot [thick, darkgray150, mark=triangle*, mark size=2, mark options={solid,rotate=180}]
table {%
0.00047454685 0.0153979966666667
};
\addplot [thick, darkgray150, mark=square*, mark size=2, mark options={solid}]
table {%
3.22709e-05 0.00209830266666667
};
\addplot [thick, darkgray150, mark=pentagon*, mark size=2, mark options={solid}]
table {%
2.05770583333333e-06 0.000265951666666667
};
\addplot [thick, darkgray150, mark=x, mark size=2, mark options={solid}]
table {%
1.2919935e-07 3.33405483333333e-05
};
\addplot [thick, black, dashed]
table {%
0.000558290411764706 0.0265527475880451
0.00047454685 0.0235049700250962
3.22709e-05 0.00312843967246703
2.05770583333333e-06 0.000396753264712192
1.2919935e-07 4.97382534365693e-05
1.098194475e-07 4.40291970633413e-05
};
\addplot [thick, steelblue49130189, dashed]
table {%
0.000474611716666667 0.0153980366666667
3.227122e-05 0.0020983035
9.26468033333333e-06 0.00026559265
9.99811216666667e-07 2.980041e-05
9.99999566666667e-08 7.3598375e-06
};
\addplot [thick, steelblue49130189, mark=triangle*, mark size=2, mark options={solid,rotate=180}]
table {%
0.000474611716666667 0.0153980366666667
};
\addplot [thick, steelblue49130189, mark=square*, mark size=2, mark options={solid}]
table {%
3.227122e-05 0.0020983035
};
\addplot [thick, steelblue49130189, mark=x, mark size=2, mark options={solid}]
table {%
9.26468033333333e-06 0.00026559265
};
\addplot [thick, steelblue49130189, mark=*, mark size=2, mark options={solid}]
table {%
9.99811216666667e-07 2.980041e-05
};
\addplot [thick, steelblue49130189, mark=*, mark size=2, mark options={solid}]
table {%
9.99999566666667e-08 7.3598375e-06
};
\addplot [thick, darkgray150]
table {%
0.000474759916666667 0.0153979966666667
3.227193e-05 0.00209830266666667
2.05771e-06 0.000265951666666667
1.2919935e-07 3.33405483333333e-05
};
\addplot [thick, darkgray150, mark=triangle*, mark size=2, mark options={solid,rotate=180}]
table {%
0.000474759916666667 0.0153979966666667
};
\addplot [thick, darkgray150, mark=square*, mark size=2, mark options={solid}]
table {%
3.227193e-05 0.00209830266666667
};
\addplot [thick, darkgray150, mark=pentagon*, mark size=2, mark options={solid}]
table {%
2.05771e-06 0.000265951666666667
};
\addplot [thick, darkgray150, mark=x, mark size=2, mark options={solid}]
table {%
1.2919935e-07 3.33405483333333e-05
};
\addplot [thick, orangered2308513, dashed]
table {%
0.0004747778 0.0153981433333333
3.22720333333333e-05 0.002098303
6.848055e-06 0.000265951666666667
9.59357516666667e-07 3.33405483333333e-05
9.79267516666667e-08 7.43860666666667e-06
};
\addplot [thick, orangered2308513, mark=triangle*, mark size=2, mark options={solid,rotate=180}]
table {%
0.0004747778 0.0153981433333333
};
\addplot [thick, orangered2308513, mark=square*, mark size=2, mark options={solid}]
table {%
3.22720333333333e-05 0.002098303
};
\addplot [thick, orangered2308513, mark=pentagon*, mark size=2, mark options={solid}]
table {%
6.848055e-06 0.000265951666666667
};
\addplot [thick, orangered2308513, mark=x, mark size=2, mark options={solid}]
table {%
9.59357516666667e-07 3.33405483333333e-05
};
\addplot [thick, orangered2308513, mark=diamond*, mark size=2, mark options={solid}]
table {%
9.79267516666667e-08 7.43860666666667e-06
};
\draw (axis cs:2.05770583333333e-06,0.000396753264712192) ++(2pt,2pt) node[
  scale=0.6,
  anchor=south west,
  text=black,
  rotate=30.6
]{$p = 0.75$};
\end{axis}

\begin{axis}[
log basis x={10},
log basis y={10},
tick align=outside,
tick pos=left,
x grid style={darkgray176},
xlabel={\(\displaystyle \vert\!\vert \gamma \vert\!\vert_{\infty}\)},
xmajorgrids,
xmin=6.3548222537649e-08, xmax=0.000860316218575168,
xmode=log,
xtick style={color=black},
y grid style={darkgray176},
ylabel={\(\displaystyle E_{\infty}\)},
ymajorgrids,
ymin=4.88659294284375e-06, ymax=0.0399918531607427,
ymode=log,
ytick style={color=black}
]
\addplot [thick, darkgray150]
table {%
0.00047454685 0.0153979966666667
3.22709e-05 0.00209830266666667
2.05770583333333e-06 0.000265951666666667
1.2919935e-07 3.33405483333333e-05
};
\addplot [thick, darkgray150, mark=triangle*, mark size=2, mark options={solid,rotate=180}]
table {%
0.00047454685 0.0153979966666667
};
\addplot [thick, darkgray150, mark=square*, mark size=2, mark options={solid}]
table {%
3.22709e-05 0.00209830266666667
};
\addplot [thick, darkgray150, mark=pentagon*, mark size=2, mark options={solid}]
table {%
2.05770583333333e-06 0.000265951666666667
};
\addplot [thick, darkgray150, mark=x, mark size=2, mark options={solid}]
table {%
1.2919935e-07 3.33405483333333e-05
};
\addplot [thick, black, dashed]
table {%
0.000558290411764706 0.0265527475880451
0.00047454685 0.0235049700250962
3.22709e-05 0.00312843967246703
2.05770583333333e-06 0.000396753264712192
1.2919935e-07 4.97382534365693e-05
1.098194475e-07 4.40291970633413e-05
};
\addplot [thick, steelblue49130189, dashed]
table {%
0.000474611716666667 0.0153980366666667
3.227122e-05 0.0020983035
9.26468033333333e-06 0.00026559265
9.99811216666667e-07 2.980041e-05
9.99999566666667e-08 7.3598375e-06
};
\addplot [thick, steelblue49130189, mark=triangle*, mark size=2, mark options={solid,rotate=180}]
table {%
0.000474611716666667 0.0153980366666667
};
\addplot [thick, steelblue49130189, mark=square*, mark size=2, mark options={solid}]
table {%
3.227122e-05 0.0020983035
};
\addplot [thick, steelblue49130189, mark=x, mark size=2, mark options={solid}]
table {%
9.26468033333333e-06 0.00026559265
};
\addplot [thick, steelblue49130189, mark=*, mark size=2, mark options={solid}]
table {%
9.99811216666667e-07 2.980041e-05
};
\addplot [thick, steelblue49130189, mark=*, mark size=2, mark options={solid}]
table {%
9.99999566666667e-08 7.3598375e-06
};
\addplot [thick, darkgray150]
table {%
0.000474759916666667 0.0153979966666667
3.227193e-05 0.00209830266666667
2.05771e-06 0.000265951666666667
1.2919935e-07 3.33405483333333e-05
};
\addplot [thick, darkgray150, mark=triangle*, mark size=2, mark options={solid,rotate=180}]
table {%
0.000474759916666667 0.0153979966666667
};
\addplot [thick, darkgray150, mark=square*, mark size=2, mark options={solid}]
table {%
3.227193e-05 0.00209830266666667
};
\addplot [thick, darkgray150, mark=pentagon*, mark size=2, mark options={solid}]
table {%
2.05771e-06 0.000265951666666667
};
\addplot [thick, darkgray150, mark=x, mark size=2, mark options={solid}]
table {%
1.2919935e-07 3.33405483333333e-05
};
\addplot [thick, orangered2308513, dashed]
table {%
0.0004747778 0.0153981433333333
3.22720333333333e-05 0.002098303
6.848055e-06 0.000265951666666667
9.59357516666667e-07 3.33405483333333e-05
9.79267516666667e-08 7.43860666666667e-06
};
\addplot [thick, orangered2308513, mark=triangle*, mark size=2, mark options={solid,rotate=180}]
table {%
0.0004747778 0.0153981433333333
};
\addplot [thick, orangered2308513, mark=square*, mark size=2, mark options={solid}]
table {%
3.22720333333333e-05 0.002098303
};
\addplot [thick, orangered2308513, mark=pentagon*, mark size=2, mark options={solid}]
table {%
6.848055e-06 0.000265951666666667
};
\addplot [thick, orangered2308513, mark=x, mark size=2, mark options={solid}]
table {%
9.59357516666667e-07 3.33405483333333e-05
};
\addplot [thick, orangered2308513, mark=diamond*, mark size=2, mark options={solid}]
table {%
9.79267516666667e-08 7.43860666666667e-06
};
\draw (axis cs:2.05770583333333e-06,0.000396753264712192) ++(2pt,2pt) node[
  scale=0.6,
  anchor=south west,
  text=black,
  rotate=30.6
]{$p = 0.75$};
\end{axis}

\end{tikzpicture}

%% file: figures/advection_anal-w40_h_vs_erri3.0.tex
\begin{tikzpicture}

\definecolor{darkgray150}{RGB}{150,150,150}
\definecolor{darkgray176}{RGB}{176,176,176}
\definecolor{orangered2308513}{RGB}{230,85,13}
\definecolor{steelblue49130189}{RGB}{49,130,189}

\begin{axis}[
log basis x={10},
log basis y={10},
tick align=outside,
tick pos=left,
x grid style={darkgray176},
xlabel={\(\displaystyle h_{f}\)},
xmajorgrids,
xmin=0.000543028010132773, xmax=0.0293850469902458,
xmode=log,
xtick style={color=black},
y grid style={darkgray176},
ylabel={\(\displaystyle E_{\infty}\)},
ymajorgrids,
ymin=4.97624785149204e-06, ymax=0.0272997642624822,
ymode=log,
ytick style={color=black}
]
\addplot [thick, darkgray150]
table {%
0.0208333333333333 0.0153979966666667
0.0104166666666667 0.00209830266666667
0.00520833333333333 0.000265951666666667
0.00260416666666667 3.33405483333333e-05
};
\addplot [thick, darkgray150, mark=triangle*, mark size=2, mark options={solid,rotate=180}]
table {%
0.0208333333333333 0.0153979966666667
};
\addplot [thick, darkgray150, mark=square*, mark size=2, mark options={solid}]
table {%
0.0104166666666667 0.00209830266666667
};
\addplot [thick, darkgray150, mark=pentagon*, mark size=2, mark options={solid}]
table {%
0.00520833333333333 0.000265951666666667
};
\addplot [thick, darkgray150, mark=x, mark size=2, mark options={solid}]
table {%
0.00260416666666667 3.33405483333333e-05
};
\addplot [thick, black, dashed]
table {%
0.0245098039215686 0.0184583414046052
0.0208333333333333 0.0113434426542086
0.0104166666666667 0.00142205011466635
0.00520833333333333 0.000178272733443254
0.00260416666666667 2.23488378936534e-05
0.00221354166666667 1.37343196486556e-05
};
\addplot [thick, steelblue49130189, dashed]
table {%
0.0208333333333333 0.0153980366666667
0.0104166666666667 0.0020983035
0.00260416666666667 0.00026559265
0.000651041666666667 2.980041e-05
0.000651041666666667 7.3598375e-06
};
\addplot [thick, steelblue49130189, mark=triangle*, mark size=2, mark options={solid,rotate=180}]
table {%
0.0208333333333333 0.0153980366666667
};
\addplot [thick, steelblue49130189, mark=square*, mark size=2, mark options={solid}]
table {%
0.0104166666666667 0.0020983035
};
\addplot [thick, steelblue49130189, mark=x, mark size=2, mark options={solid}]
table {%
0.00260416666666667 0.00026559265
};
\addplot [thick, steelblue49130189, mark=*, mark size=2, mark options={solid}]
table {%
0.000651041666666667 2.980041e-05
};
\addplot [thick, steelblue49130189, mark=*, mark size=2, mark options={solid}]
table {%
0.000651041666666667 7.3598375e-06
};
\addplot [thick, darkgray150]
table {%
0.0208333333333333 0.0153979966666667
0.0104166666666667 0.00209830266666667
0.00520833333333333 0.000265951666666667
0.00260416666666667 3.33405483333333e-05
};
\addplot [thick, darkgray150, mark=triangle*, mark size=2, mark options={solid,rotate=180}]
table {%
0.0208333333333333 0.0153979966666667
};
\addplot [thick, darkgray150, mark=square*, mark size=2, mark options={solid}]
table {%
0.0104166666666667 0.00209830266666667
};
\addplot [thick, darkgray150, mark=pentagon*, mark size=2, mark options={solid}]
table {%
0.00520833333333333 0.000265951666666667
};
\addplot [thick, darkgray150, mark=x, mark size=2, mark options={solid}]
table {%
0.00260416666666667 3.33405483333333e-05
};
\addplot [thick, black, dashed]
table {%
0.0245098039215686 0.0184583414046052
0.0208333333333333 0.0113434426542086
0.0104166666666667 0.00142205011466635
0.00520833333333333 0.000178272733443254
0.00260416666666667 2.23488378936534e-05
0.00221354166666667 1.37343196486556e-05
};
\addplot [thick, orangered2308513, dashed]
table {%
0.0208333333333333 0.0153981433333333
0.0104166666666667 0.002098303
0.00520833333333333 0.000265951666666667
0.00260416666666667 3.33405483333333e-05
0.00130208333333333 7.43860666666667e-06
};
\addplot [thick, orangered2308513, mark=triangle*, mark size=2, mark options={solid,rotate=180}]
table {%
0.0208333333333333 0.0153981433333333
};
\addplot [thick, orangered2308513, mark=square*, mark size=2, mark options={solid}]
table {%
0.0104166666666667 0.002098303
};
\addplot [thick, orangered2308513, mark=pentagon*, mark size=2, mark options={solid}]
table {%
0.00520833333333333 0.000265951666666667
};
\addplot [thick, orangered2308513, mark=x, mark size=2, mark options={solid}]
table {%
0.00260416666666667 3.33405483333333e-05
};
\addplot [thick, orangered2308513, mark=diamond*, mark size=2, mark options={solid}]
table {%
0.00130208333333333 7.43860666666667e-06
};
\draw (axis cs:0.00520833333333333,0.000178272733443254) ++(2pt,-2pt) node[
  scale=0.6,
  anchor=north west,
  text=black,
  rotate=46.0
]{$p = 3.00$};
\draw (axis cs:0.00520833333333333,0.000178272733443254) ++(2pt,-2pt) node[
  scale=0.6,
  anchor=north west,
  text=black,
  rotate=46.0
]{$p = 3.00$};
\end{axis}

\end{tikzpicture}

%% file: figures/advection_anal-w40_h_vs_det3.0.tex
\begin{tikzpicture}

\definecolor{darkgray150}{RGB}{150,150,150}
\definecolor{darkgray176}{RGB}{176,176,176}
\definecolor{orangered2308513}{RGB}{230,85,13}
\definecolor{steelblue49130189}{RGB}{49,130,189}

\begin{axis}[
log basis x={10},
log basis y={10},
tick align=outside,
tick pos=left,
x grid style={darkgray176},
xlabel={\(\displaystyle h_{f}\)},
xmajorgrids,
xmin=0.000543028010132773, xmax=0.0293850469902458,
xmode=log,
xtick style={color=black},
y grid style={darkgray176},
ylabel={\(\displaystyle \vert\!\vert \gamma \vert\!\vert_{\infty}\)},
ymajorgrids,
ymin=2.79490989805575e-08, ymax=0.00109920123634941,
ymode=log,
ytick style={color=black}
]
\addplot [thick, darkgray150]
table {%
0.0208333333333333 0.00047454685
0.0104166666666667 3.22709e-05
0.00520833333333333 2.05770583333333e-06
0.00260416666666667 1.2919935e-07
};
\addplot [thick, darkgray150, mark=triangle*, mark size=2, mark options={solid,rotate=180}]
table {%
0.0208333333333333 0.00047454685
};
\addplot [thick, darkgray150, mark=square*, mark size=2, mark options={solid}]
table {%
0.0104166666666667 3.22709e-05
};
\addplot [thick, darkgray150, mark=pentagon*, mark size=2, mark options={solid}]
table {%
0.00520833333333333 2.05770583333333e-06
};
\addplot [thick, darkgray150, mark=x, mark size=2, mark options={solid}]
table {%
0.00260416666666667 1.2919935e-07
};
\addplot [thick, black, dashed]
table {%
0.0245098039215686 0.000679558393651641
0.0208333333333333 0.000354733728726117
0.0104166666666667 2.21708580453823e-05
0.00520833333333333 1.38567862783639e-06
0.00260416666666667 8.66049142397746e-08
0.00221354166666667 4.52083065138764e-08
};
\addplot [thick, steelblue49130189, dashed]
table {%
0.0208333333333333 0.000474611716666667
0.0104166666666667 3.227122e-05
0.00260416666666667 9.26468033333333e-06
0.000651041666666667 9.99811216666667e-07
0.000651041666666667 9.99999566666667e-08
};
\addplot [thick, steelblue49130189, mark=triangle*, mark size=2, mark options={solid,rotate=180}]
table {%
0.0208333333333333 0.000474611716666667
};
\addplot [thick, steelblue49130189, mark=square*, mark size=2, mark options={solid}]
table {%
0.0104166666666667 3.227122e-05
};
\addplot [thick, steelblue49130189, mark=x, mark size=2, mark options={solid}]
table {%
0.00260416666666667 9.26468033333333e-06
};
\addplot [thick, steelblue49130189, mark=*, mark size=2, mark options={solid}]
table {%
0.000651041666666667 9.99811216666667e-07
};
\addplot [thick, steelblue49130189, mark=*, mark size=2, mark options={solid}]
table {%
0.000651041666666667 9.99999566666667e-08
};
\addplot [thick, darkgray150]
table {%
0.0208333333333333 0.000474759916666667
0.0104166666666667 3.227193e-05
0.00520833333333333 2.05771e-06
0.00260416666666667 1.2919935e-07
};
\addplot [thick, darkgray150, mark=triangle*, mark size=2, mark options={solid,rotate=180}]
table {%
0.0208333333333333 0.000474759916666667
};
\addplot [thick, darkgray150, mark=square*, mark size=2, mark options={solid}]
table {%
0.0104166666666667 3.227193e-05
};
\addplot [thick, darkgray150, mark=pentagon*, mark size=2, mark options={solid}]
table {%
0.00520833333333333 2.05771e-06
};
\addplot [thick, darkgray150, mark=x, mark size=2, mark options={solid}]
table {%
0.00260416666666667 1.2919935e-07
};
\addplot [thick, black, dashed]
table {%
0.0245098039215686 0.000679558393651641
0.0208333333333333 0.000354733728726117
0.0104166666666667 2.21708580453823e-05
0.00520833333333333 1.38567862783639e-06
0.00260416666666667 8.66049142397746e-08
0.00221354166666667 4.52083065138764e-08
};
\addplot [thick, orangered2308513, dashed]
table {%
0.0208333333333333 0.0004747778
0.0104166666666667 3.22720333333333e-05
0.00520833333333333 6.848055e-06
0.00260416666666667 9.59357516666667e-07
0.00130208333333333 9.79267516666667e-08
};
\addplot [thick, orangered2308513, mark=triangle*, mark size=2, mark options={solid,rotate=180}]
table {%
0.0208333333333333 0.0004747778
};
\addplot [thick, orangered2308513, mark=square*, mark size=2, mark options={solid}]
table {%
0.0104166666666667 3.22720333333333e-05
};
\addplot [thick, orangered2308513, mark=pentagon*, mark size=2, mark options={solid}]
table {%
0.00520833333333333 6.848055e-06
};
\addplot [thick, orangered2308513, mark=x, mark size=2, mark options={solid}]
table {%
0.00260416666666667 9.59357516666667e-07
};
\addplot [thick, orangered2308513, mark=diamond*, mark size=2, mark options={solid}]
table {%
0.00130208333333333 9.79267516666667e-08
};
\draw (axis cs:0.00520833333333333,1.38567862783639e-06) ++(2pt,-2pt) node[
  scale=0.6,
  anchor=north west,
  text=black,
  rotate=48.4
]{$p = 4.00$};
\draw (axis cs:0.00520833333333333,1.38567862783639e-06) ++(2pt,-2pt) node[
  scale=0.6,
  anchor=north west,
  text=black,
  rotate=48.4
]{$p = 4.00$};
\end{axis}

\end{tikzpicture}

%% file: figures/advection_anal-ec_nblock_vs_erri3.0.tex
\begin{tikzpicture}

\definecolor{darkgray150}{RGB}{150,150,150}
\definecolor{darkgray176}{RGB}{176,176,176}
\definecolor{orangered2308513}{RGB}{230,85,13}
\definecolor{steelblue49130189}{RGB}{49,130,189}

\begin{axis}[
log basis x={10},
log basis y={10},
tick align=outside,
tick pos=left,
x grid style={darkgray176},
xlabel={\# blocks},
xmajorgrids,
xmin=40.4910129377893, xmax=333215.274726013,
xmode=log,
xtick style={color=black},
y grid style={darkgray176},
ylabel={\(\displaystyle E_{\infty}\)},
ymajorgrids,
ymin=4.05777767454551e-06, ymax=0.0427542739016854,
ymode=log,
ytick style={color=black}
]
\addplot [thick, darkgray150, mark=*, mark size=2, mark options={solid}]
table {%
432 0.0153979966666667
3456 0.00209830266666667
27648 0.000265951666666667
221184 3.33405483333333e-05
};
\addplot [thick, black, dashed]
table {%
508.235294117647 0.0238533212377447
432 0.0280627308679349
3456 0.00350784135849186
27648 0.000438480169811483
221184 5.48100212264354e-05
188006.4 6.44823779134534e-05
};
\addplot [thick, darkgray150, mark=*, mark size=2, mark options={solid}]
table {%
432 0.0153979966666667
3456 0.00209830266666667
27648 0.000265951666666667
221184 3.33405483333333e-05
};
\addplot [thick, black, dashed]
table {%
508.235294117647 0.0238533212377447
432 0.0280627308679349
3456 0.00350784135849186
27648 0.000438480169811483
221184 5.48100212264354e-05
188006.4 6.44823779134534e-05
};
\addplot [thick, darkgray150, mark=*, mark size=2, mark options={solid}]
table {%
432 0.0153979966666667
3456 0.00209830266666667
27648 0.000265951666666667
221184 3.33405483333333e-05
};
\addplot [thick, black, dashed]
table {%
508.235294117647 0.0238533212377447
432 0.0280627308679349
3456 0.00350784135849186
27648 0.000438480169811483
221184 5.48100212264354e-05
188006.4 6.44823779134534e-05
};
\addplot [thick, steelblue49130189, dotted, mark=*, mark size=2, mark options={solid}]
table {%
61 0.0153973433333333
299 0.0020983035
1083 0.000265371916666667
3463 3.319949e-05
15384 1.129471e-05
};
\addplot [thick, steelblue49130189, dashed, mark=*, mark size=2, mark options={solid}]
table {%
68 0.0153980366666667
327 0.0020983035
1176.33333333333 0.00026559265
5241 2.980041e-05
19334.3333333333 7.3598375e-06
};
\addplot [thick, steelblue49130189, mark=*, mark size=2, mark options={solid}]
table {%
68 0.0153980366666667
327 0.0020983035
2028 0.000262165466666667
9815.5 2.70218783333333e-05
43923 6.182126e-06
};
\addplot [thick, darkgray150, mark=*, mark size=2, mark options={solid}]
table {%
432 0.0153979966666667
3456 0.00209830266666667
27648 0.000265951666666667
221184 3.33405483333333e-05
};
\addplot [thick, black, dashed]
table {%
508.235294117647 0.0238533212377447
432 0.0280627308679349
3456 0.00350784135849186
27648 0.000438480169811483
221184 5.48100212264354e-05
188006.4 6.44823779134534e-05
};
\addplot [thick, darkgray150, mark=*, mark size=2, mark options={solid}]
table {%
432 0.0153979966666667
3456 0.00209830266666667
27648 0.000265951666666667
221184 3.33405483333333e-05
};
\addplot [thick, black, dashed]
table {%
508.235294117647 0.0238533212377447
432 0.0280627308679349
3456 0.00350784135849186
27648 0.000438480169811483
221184 5.48100212264354e-05
188006.4 6.44823779134534e-05
};
\addplot [thick, darkgray150, mark=*, mark size=2, mark options={solid}]
table {%
432 0.0153979966666667
3456 0.00209830266666667
27648 0.000265951666666667
221184 3.33405483333333e-05
};
\addplot [thick, black, dashed]
table {%
508.235294117647 0.0238533212377447
432 0.0280627308679349
3456 0.00350784135849186
27648 0.000438480169811483
221184 5.48100212264354e-05
188006.4 6.44823779134534e-05
};
\addplot [thick, orangered2308513, dotted, mark=*, mark size=2, mark options={solid}]
table {%
61 0.0153977433333333
299 0.0020983035
915 0.000265951666666667
2844.66666666667 3.33405483333333e-05
12507 7.647336e-06
};
\addplot [thick, orangered2308513, dashed, mark=*, mark size=2, mark options={solid}]
table {%
68 0.0153981433333333
327 0.002098303
999 0.000265951666666667
3218 3.33405483333333e-05
12724 7.43860666666667e-06
};
\addplot [thick, orangered2308513, mark=*, mark size=2, mark options={solid}]
table {%
68 0.0153981433333333
327 0.002098303
1139 0.000265951666666667
4037 3.33405483333333e-05
19171 7.391851e-06
};
\draw (axis cs:27648,0.000438480169811483) ++(2pt,2pt) node[
  scale=0.6,
  anchor=south west,
  text=black,
  rotate=324.0
]{$p = 1.00$};
\draw (axis cs:27648,0.000438480169811483) ++(2pt,2pt) node[
  scale=0.6,
  anchor=south west,
  text=black,
  rotate=324.0
]{$p = 1.00$};
\draw (axis cs:27648,0.000438480169811483) ++(2pt,2pt) node[
  scale=0.6,
  anchor=south west,
  text=black,
  rotate=324.0
]{$p = 1.00$};
\draw (axis cs:27648,0.000438480169811483) ++(2pt,2pt) node[
  scale=0.6,
  anchor=south west,
  text=black,
  rotate=324.0
]{$p = 1.00$};
\draw (axis cs:27648,0.000438480169811483) ++(2pt,2pt) node[
  scale=0.6,
  anchor=south west,
  text=black,
  rotate=324.0
]{$p = 1.00$};
\draw (axis cs:27648,0.000438480169811483) ++(2pt,2pt) node[
  scale=0.6,
  anchor=south west,
  text=black,
  rotate=324.0
]{$p = 1.00$};
\end{axis}

\end{tikzpicture}

%% file: figures/convenright_eps.tex
\begin{tikzpicture}

\definecolor{darkgray176}{RGB}{176,176,176}
\definecolor{orangered2308513}{RGB}{230,85,13}

\begin{axis}[
log basis x={10},
log basis y={10},
tick align=outside,
tick pos=left,
x grid style={darkgray176},
xlabel={\(\displaystyle \epsilon_{r}\), \(\displaystyle \normi{\gamma}\)},
xmajorgrids,
xmin=0.000482299470120457, xmax=1.70560858680107,
xmode=log,
xtick style={color=black},
y grid style={darkgray176},
ylabel={\(\displaystyle E_{\infty}\)},
ymajorgrids,
ymin=0.00308065705376293, ymax=1.20254718730184,
ymode=log,
ytick style={color=black}
]
\addplot [thick, black, dashed]
table {%
1.17647058823529 0.916871653997366
1 0.811657377652099
0.1 0.144335360268491
0.01 0.0256668599306008
0.001 0.00456428485349405
0.00085 0.00404051696755239
};
\addplot [thick, orangered2308513, mark=square*, mark size=2, mark options={solid}]
table {%
1 0.250270718003362
0.1 0.150544605610411
0.01 0.0431423425465493
0.001 0.00504431488189479
};
\addplot [thick, orangered2308513, mark=*, mark size=2, mark options={solid}]
table {%
0.138247 0.250270718003362
0.05715398 0.150544605610411
0.007352768 0.0431423425465493
0.000699222 0.00504431488189479
};
\draw (axis cs:0.01,0.0256668599306008) ++(2pt,-2pt) node[
  scale=0.6,
  anchor=north west,
  text=black,
  rotate=37.4
]{$p = 0.75$};
\end{axis}

\end{tikzpicture}

%% file: figures/convenright_blocks.tex
\begin{tikzpicture}

\definecolor{darkgray176}{RGB}{176,176,176}
\definecolor{orangered2308513}{RGB}{230,85,13}

\begin{axis}[
log basis x={10},
log basis y={10},
tick align=outside,
tick pos=left,
x grid style={darkgray176},
xlabel={\#blocks},
xmajorgrids,
xmin=6.20485790688062, xmax=1661.9236338297,
xmode=log,
xtick style={color=black},
y grid style={darkgray176},
ylabel={\(\displaystyle E_{\infty}\)},
ymajorgrids,
ymin=0.00354009237030725, ymax=0.948186608694087,
ymode=log,
ytick style={color=black}
]
\addplot [thick, orangered2308513, mark=*, mark size=2, mark options={solid}]
table {%
8 0.250270718003362
36 0.150544605610411
274 0.0431423425465493
1289 0.00504431488189479
};
\addplot [thick, black, dashed]
table {%
9.41176470588235 0.625107337466345
8 0.735420397019229
36 0.163426754893162
274 0.0214721283801235
1289 0.00456428485349405
1095.65 0.00536974688646359
};
\draw (axis cs:274,0.0214721283801235) ++(2pt,-2pt) node[
  scale=0.6,
  anchor=north west,
  text=black,
  rotate=323.3
]{$p = 1.00$};
\end{axis}

\end{tikzpicture}

%% file: wavelet_6_scalability.tex
\section{Weak scalability and performance analysis}
\label{sec:scalability}
In this section we present the results of the weak scalability campaign we have done on \typo{Cori, a Cray XC40 supercomputer whose CPU partition contains $2388$ nodes of $32$ cores each}\footnote{The MPI implementation used is \code{OpenMPI/4.1.2}}.
The testcase is the execution of $50$ times steps for the advection of a Gaussian tube aligned with the $z$-direction in a rectangular domain. The size of the domain, and therefore the size of the simulation $S$, is adapted to maintain the number of blocks per core constant: $S\left(N_c \right) = \left[ 1 \times 1 \times N_c/32 \right]$, with $N_c$ the number of cores.
To get a significant sample of each operation, we perform the adaptation every time step, although in practice we would adapt less frequently. 
The domain is initialized with 75.5 blocks per rank and one rank per core, which with a block size of $24^3$ leads to about 1 million unknowns per core. We used the \wave{4}{0} and $\cfl = 0.25$.

\begin{figure}[ht!]
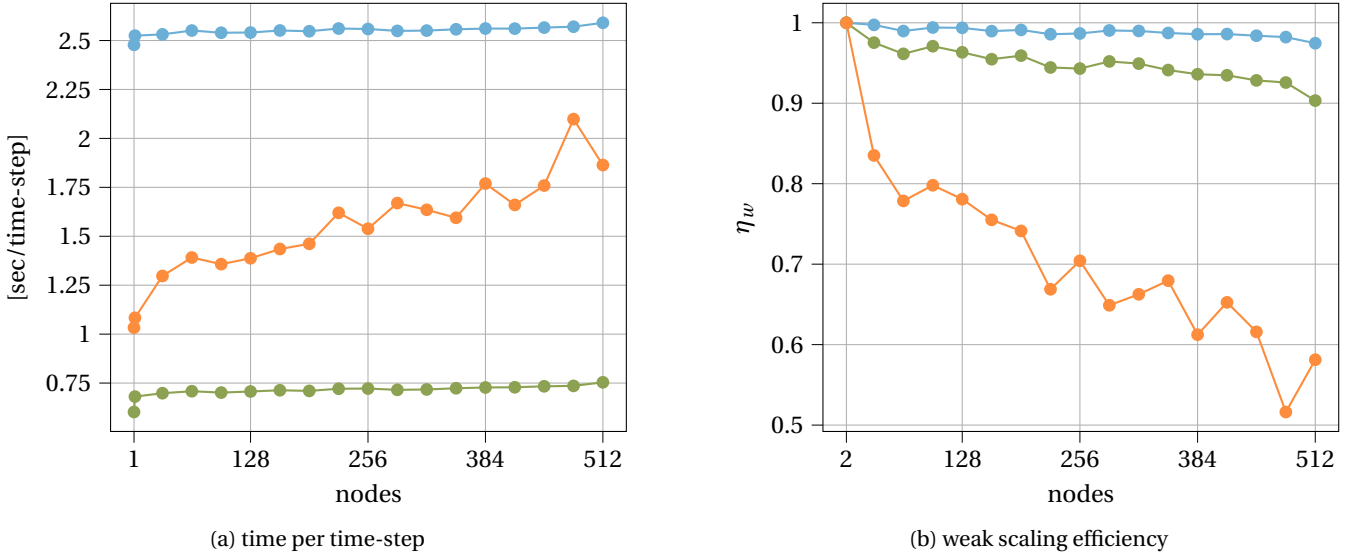

\begin{minipage}{0.49\textwidth}
	\centering
	\InputIfFileExists{figures/time.tex}{}{}
	\subcaption{\typo{time per time-step}}
	\label{fig:weak_time}
\end{minipage}%
\hfill%
\begin{minipage}{0.49\textwidth}
	\centering
	\InputIfFileExists{figures/weak.tex}{}{}
	\subcaption{weak scaling efficiency}
	\label{fig:weak_weak}
\end{minipage}%
\caption{\revone{Weak scaling for the three main operations involved in adaptive grid simulations: the grid adaptation (\CaptionWeakAdapt, performed every time step here), the stencil operations (\CaptionWeakStencil), and the time spent in the ghost computations ( \CaptionWeakGhost ) during the stencil evaluation. The times are given in seconds per time-step, from $1$ node to $512$ nodes ($32$ cores to $16,384$ cores).}}
\end{figure}

We study the performance of our framework from $1$ nodes to \typo{$512$} nodes, with the latter corresponding to \typo{$16,384$} cores and about \typo{$17$} billion unknowns. In \fref{fig:weak_time} we show the evolution of the \revone{time spent per timestep for the grid adaptation (orange) and the stencil computation (blue). We also show separately the time spent in ghost computations (green), which is the communication-heavy part of the stencil computation. Though highly problem-dependent, for this particular problem the grid adaptation takes roughly half the time of the stencil computation within an RK3 timestep, which is very reasonable given that we typically would need to adapt the grid only every 10 to 100 timesteps. Further, the ghosting takes up about a third of the time of a stencil computation, even on very large partitions. }

Based on this timing data, we show in \fref{fig:weak_weak} the weak efficiency $\eta_w$ defined as
\begin{eqc}
\eta_w(N_c) = \dfrac{ T\left(N_{\text{ref}},S\left(N_{\text{ref}} \right) \right) }{ T\left(N_{c},S\left(N_c \right)\right)}
\end{eqc}
where $T\left(N,S\left(N\right)\right)$ is the time taken by the simulation to run a problem of size $S$ on $N$ cores and where we used \typo{$N_{\text{ref}} = 64$, \textit{i.e.} $2$ nodes}. %

Based on \frefs{fig:weak_time}{fig:weak_weak}, the Post-Start-Complete-Wait (PSCW) strategy allows us to reach a perfect scalability in the case of the ghost and the stencil operation. This observation is confirmed by analyzing the break-down of the different operations illustrated in \fref{fig:weak_pie_ghost}.
This result is only possible because both operations are scalable and done on group of ranks with the PSCW calls, instead of using the entire communicator.%
We do pay a price to achieve \revtwo{this scalability} in the ghost computation, because we have to re-initialize the ghost meta-data structure every time we modify the grid. As observed in \fref{fig:weak_weak}, this adaptation process has a lower parallel efficiency. Analyzing the timing of different operations within the adaptation step (\fref{fig:weak_pie_adapt}) shows that the less scalable operations are the synchronization step, which contains a non-blocking \code{MPI\_Allreduce} and various RMA synchronizations, and the \code{reset} operation of the meta-data for the ghosting, which involves the creation and deletion of the window on \code{MPI\_COMM\_WORLD}. For the latter we currently use non-dynamic windows as we do not expect the adaptation to be called often compared to the use of the ghosting and the stencils, and so its computation cost will not dominate in practice. However, in the future we could still improve the implementation by considering a dynamic window allocation, which might be better suited for a policy that requires very frequent adaptations. %
Overall, the observed behavior of this weak scalability test is consistent with the implementation choices, as all the global operations are less scalable by definition, while the \revtwo{local operations} demonstrate \revtwo{excellent} scalability.

%% file: figures/time.tex
\begin{tikzpicture}

\definecolor{coral25314160}{RGB}{253,141,60}
\definecolor{cornflowerblue107174214}{RGB}{107,174,214}
\definecolor{darkgray176}{RGB}{176,176,176}
\definecolor{yellowgreen14016282}{RGB}{140,162,82}

\begin{axis}[
minor xtick={},
minor ytick={},
tick align=outside,
tick pos=left,
x grid style={darkgray176},
xlabel={nodes},
xmajorgrids,
xmin=-24.55, xmax=537.55,
xtick style={color=black},
xtick={1,128,256,384,512},
y grid style={darkgray176},
ylabel={[sec/time-step]},
ymajorgrids,
ymin=0.502756797, ymax=2.690143383,
ytick style={color=black},
ytick={0.5,0.75,1,1.25,1.5,1.75,2,2.25,2.5,2.75}
]
\addplot [thick, cornflowerblue107174214, mark=*, mark size=2, mark options={solid}]
table {%
1 2.47794066
2 2.52483342
32 2.53144448
64 2.551065
96 2.53984212
128 2.54080952
160 2.5513529
192 2.54746938
224 2.56130098
256 2.55891172
288 2.54922242
320 2.55077864
352 2.55728902
384 2.56124454
416 2.56078606
448 2.56611084
480 2.57065244
512 2.59071672
};
\addplot [thick, coral25314160, mark=*, mark size=2, mark options={solid}]
table {%
1 1.03295658
2 1.08336024
32 1.2972414
64 1.3914301
96 1.3574841
128 1.3875524
160 1.43476166
192 1.46142746
224 1.61984044
256 1.53841746
288 1.66948716
320 1.63517504
352 1.5943524
384 1.76918386
416 1.6602332
448 1.75910558
480 2.09841638
512 1.86383186
};
\addplot [thick, yellowgreen14016282, mark=*, mark size=2, mark options={solid}]
table {%
1 0.60218346
2 0.68095382
32 0.69819916
64 0.70831568
96 0.70137944
128 0.7069205
160 0.71326334
192 0.70995386
224 0.72100916
256 0.72204972
288 0.71532344
320 0.71737004
352 0.72346766
384 0.72750512
416 0.7285104
448 0.7334722
480 0.73557714
512 0.7538712
};
\end{axis}

\end{tikzpicture}

%% file: figures/weak.tex
\begin{tikzpicture}

\definecolor{coral25314160}{RGB}{253,141,60}
\definecolor{cornflowerblue107174214}{RGB}{107,174,214}
\definecolor{darkgray176}{RGB}{176,176,176}
\definecolor{yellowgreen14016282}{RGB}{140,162,82}

\begin{axis}[
minor xtick={},
minor ytick={},
tick align=outside,
tick pos=left,
x grid style={darkgray176},
xlabel={nodes},
xmajorgrids,
xmin=-23.5, xmax=537.5,
xtick style={color=black},
xtick={2,128,256,384,512},
y grid style={darkgray176},
ylabel={\(\displaystyle \eta_{w}\)},
ymajorgrids,
ymin=0.492088911829787, ymax=1.02418624229382,
ytick style={color=black},
ytick={0.4,0.5,0.6,0.7,0.8,0.9,1,1.1}
]
\addplot [thick, cornflowerblue107174214, mark=*, mark size=2, mark options={solid}]
table {%
2 1
32 0.997388423861463
64 0.989717400379841
96 0.994090695684659
128 0.993712200826452
160 0.989605718597376
192 0.991114334807039
224 0.985762095011575
256 0.986682502669533
288 0.990432768906842
320 0.989828509776136
352 0.98730859134569
384 0.985783817424946
416 0.985960310952333
448 0.983914404882059
480 0.98217611245805
512 0.974569469718017
};
\addplot [thick, yellowgreen14016282, mark=*, mark size=2, mark options={solid}]
table {%
2 1
32 0.975300256734769
64 0.961370529027396
96 0.970877931637118
128 0.963267892217017
160 0.95470183565021
192 0.959152218708974
224 0.944445449208995
256 0.943084390365805
288 0.95195233641442
320 0.949236491671718
352 0.941236018759982
384 0.936012409094798
416 0.934720794651662
448 0.928397586166183
480 0.925740867912236
512 0.903276076868303
};
\addplot [thick, coral25314160, mark=*, mark size=2, mark options={solid}]
table {%
2 1
32 0.835126168498785
64 0.778594799695651
96 0.798064772913362
128 0.780770686570107
160 0.755080282811572
192 0.741302780775722
224 0.668806762226531
256 0.704204332158321
288 0.64891798269356
320 0.662534721665027
352 0.679498610219422
384 0.612350284497847
416 0.652534981230348
448 0.615858566033313
480 0.516275154123606
512 0.581254276874524
};
\end{axis}

\end{tikzpicture}

%% file: wavelet_7_conclusion.tex
\section{Conclusion}
\label{sec:conclusion}
This work provides a detailed explanation of the mathematical foundation and distributed computational implementation of a 3D block-structured adaptive grid method based on a multiresolution analysis using wavelets, named \code{murphy}. In our approach we apply significant emphasis on the handling of resolution jumps in block-structured grids to provide consistency with non-lifted and lifted interpolating wavelets of second, fourth, and sixth polynomial order. We validate the implementation through rigorous tests of error control and moment conservation on static grids. 

Compared with most existing adaptive mesh refinement approaches, the wavelet-based approach provides explicit grid adaptation metrics that are intrinsically linked to the point-wise error made in the field compared with a polynomial interpolation. \revtwo{The wavelet} framework provides a consistent multiresolution perspective for adaptation metric, adaptation procedure, and ghost reconstruction across resolution jumps, with a formal separation of scales in all operations. When combined with finite-difference schemes we demonstrate the ability of our approach to reach point-wise high order convergence on multi-level grids. Further, the nature of the non-overlapping octree-based grids with constant-size blocks provides excellent opportunities for scalability. We exploit this in our implementation using state-of-the-art one-sided communication strategies that show excellent scalability of the grid adaptation, ghost reconstruction, and stencil computation processes at least up to \typo{$16,384$ cores}.

We tested our multiresolution adaptive grid algorithm on the convergence of simple linear hyperbolic conservation laws in the form of scalar advection using divergence-free velocity fields. The results demonstrate that the refinement threshold, which provides user control over the error permitted as detected by the wavelet analysis, is an excellent indicator for the global field error even during the evolution of this partial differential equation. As such, reducing the refinement threshold leads to convergence of the error, along a slope that can be captured by the ratio of the convergence order of the finite-difference scheme, and the polynomial order of accuracy of the wavelet used. For lifted wavelets we observe that the global maximum error on any multi-level grid is very close to the error on a uniform-resolution grid of the same maximum resolution. Compared to non-lifted wavelets, the reduced aliasing properties of lifted wavelets enable them to detect the need for refinement and opportunity for compression more efficiently, leading ultimately to a smaller number of blocks in the grid by about a factor of two for the same maximum error. 

\revtwo{The implementation and performance benefits of our octree-based block-structured grid implementation does lead to constraints in terms of the adaptation patterns that can be achieved, for instance compared to patch-based adaptive mesh refinement techniques. The use of `forest-of-trees' as offered by \code{p4est} enables us to create arbitrary aspect ratio rectangular domains, rather than remain restricted to cubic domains as in previous codes \citep{Rossinelli:2015}. However, our current implementation is still limited to cubic blocks with isotropic grids. To gain more flexibility, one can build on domain mapping techniques~\citep{Brown-Dymkoski:2017} or consider extending work in computer graphics that generalizes wavelets to work on arbitrary compact surfaces~\citep{Lounsbery:1997}. Closer to our approach, it should also be relatively straightforward to use a larger set of blocks that enable separate refinement or compression in each of the three cartesian directions. We reserve this extension for future work.}

With the fundamentals and implementation of our framework presented, future work will focus on investigating the performance of the wavelet-based refinement criterion to accurately capture refinement requirements for non-linear PDEs. Previous work on wavelet-adapted grids is promising in this regard (e.g.~\citep{Rossinelli:2015}) but systematic investigations are lacking. At the same time, we note that this level of rigor is also often absent from other adaptive mesh refinement methodologies, which typically rely on heuristic \revone{or \textit{post-hoc}} criteria for refinement and compression that make it difficult to compare their ability to capture emerging scales, or discard information that does not affect the overall error. In our case, however, the formal wavelet multiresolution theory provides a useful perspective to frame this discussion and analyze convergence with respect to the refinement threshold, compared with other possible policies. Further, we wish to explore the benefits of moment conservation offered by lifted wavelets as demonstrated in our solver, for the solution of conservative PDEs such as the hyperbolic conservation law considered in this work. Finally, we intend to combine our multiresolution framework with an elliptic solver to handle problems such as the free-space incompressible Navier-Stokes equations, as well as an immersed interface method \citep{Gillis:2019, Gabbard:2022} which will enable us to perform high-order simulations with embedded interfaces for multiphysics problems.

%% file: wavelet_s2_implementation.tex
\section{Implementation details - filters bank}
\label{sect:a0_implementation}

In this section we detail the exact construction of the filter banks, as used in the software. From an implementation perspective, two approaches have been considered in the literature: the step-by-step implementation of the lifting scheme, as in \citep{Fernandez:1996}, or the filter-based approach as in \citep{Rossinelli:2015}. Although the operation count is smaller in the first approach and a factor $2$ has been reported in the complexity \citep{Daubechies:1998}, we follow the second approach which provides significant simplifications especially in multiple dimensions. Further, the memory layout in the implementation is achieved more naturally when interlacing the scaling and detail coefficients:
\begin{eqd}
u^{L} = \left[ \dots \lambda^{L}_{k-1} \;,\; \gamma^{L}_{k-1} \;,\;\lambda^{L}_k \;,\;\gamma^{L}_k \;,\; \lambda^{L}_{k+1} \;,\; \gamma^{L}_{k+1} \;,\; \dots \right]
\end{eqd}
This does not affect the $G^a$ and $H^a$ filters used in the analysis, but requires the definition of two new filters $J^s$ and $K^s$ to replace $H^s$ and $G^s$ during the synthesis operation, such that
\be
\lambda^{L+1}_{2k} = J^s \; u^{L} \qquad \lambda^{L+1}_{2k+1} = K^s \; u^{L}
\eec
where the array $u^{L}$ in the left equation is implicitly assumed to be centered on $\lambda^{L}_k$, and in the right equation on $\gamma^{L}_k$. 
This change of perspective is strictly equivalent to the classical $H^s$/$G^s$ filters and has been done to simplify the implementation. Considering this new approach, the filter coefficients for \wave{2}{0}, \wave{4}{0}, and \wave{6}{0} are given in \tref{table_coef_nonlifted}, and the coefficients for \wave{2}{2}, \wave{4}{2}, and \wave{6}{2} are given in \tref{table_coef_lifted}.

\begin{table}[h!]
\centering
\begin{minipage}{0.34\textwidth}
\centering
\def\arraystretch{1.5}
\begin{tabular}{r|cccc}
$H^a$&&$ 1 $ \\
$G^a$ & $ -\frac{1}{2} $&$ 1 $&$ -\frac{1}{2} $\\
\hline
$J^s$ &&$ 1 $ \\
$K^s$ &$ \frac{1}{2} $&$ 1 $&$ \frac{1}{2} $
\end{tabular}
\subcaption{\wave{2}{0}}
\label{table_coef_interp_2_0}
\end{minipage}%
\hfill%
\begin{minipage}{0.64\textwidth}
\centering
\def\arraystretch{1.5}
\begin{tabular}{r|ccccccccc}
$H^a$&&&&$1$ \\
$G^a$ &$\frac{1}{16} $&$ 0 $&$ -\frac{9}{16} $&$ 1 $&$ -\frac{9}{16} $&$ 0 $&$ \frac{1}{16}$\\
\hline
$J^s$ & & &&$ 1 $\\
$K^s$ &$ -\frac{1}{16} $&$ 0 $&$ \frac{9}{16} $&$ 1 $&$ \frac{9}{16} $&$ 0 $&$ -\frac{1}{16} $
\end{tabular}
\subcaption{\wave{4}{0}}
\label{table_coef_interp_4_0}
\end{minipage}

\begin{minipage}{\textwidth}
\def\arraystretch{1.5}
\centering
\begin{tabular}{r|ccccccccccc}
$H^a$& & & & & &  $1$ \\
$G^a$ & $-\frac{3}{256}$ & $0$ & $\frac{25}{256}$ & $0$ & $-\frac{75}{128}$ & $1$ & $-\frac{75}{128}$ & $0$ & $\frac{25}{256}$ & $0$ & $-\frac{3}{256}$\\
\hline
$J^s$ & & & & & &  $1$\\
$K^s$ & $\frac{3}{256}$ & $0$ & $-\frac{25}{256}$ & $0$ & $\frac{75}{128}$ & $1$ & $\frac{75}{128}$ & $0$ & $-\frac{25}{256}$ & $0$ & $\frac{3}{256}$
\end{tabular}
\subcaption{\wave{6}{0}}
\label{table_coef_interp_6_0}
\end{minipage}

\caption{Non-lifted interpolating wavelet filter coefficients as operated on interleaved scaling and detail coefficients. The filters $H^a$ and $J^s$ are centered on (even) scaling coefficients, and the filters $G^a$ and $K^s$ are centered on odd scaling coefficients and detail coefficients, respectively.}
\label{table_coef_nonlifted}
\end{table}

\begin{table}[h!]
\centering
\begin{minipage}{0.34\textwidth}
\centering
\def\arraystretch{1.5}
\begin{tabular}{r|ccccc}
$H^a$&$ -\frac{1}{8} $&$ \frac{1}{4} $&$ \frac{3}{4} $&$ \frac{1}{4} $&$ -\frac{1}{8} $ \\
$G^a$ & &$ -\frac{1}{2} $&$ 1 $&$ -\frac{1}{2} $\\
\hline
$J^s$ & &$ -\frac{1}{4} $&$ 1 $&$ -\frac{1}{4} $\\
$K^s$ &$ -\frac{1}{8} $&$ \frac{1}{2} $&$ \frac{3}{4} $&$ \frac{1}{2} $&$ -\frac{1}{8} $
\end{tabular}
\subcaption{\wave{2}{2}}
\label{table_coef_interp_2_2}
\end{minipage}%
\hfill%
\begin{minipage}{0.64\textwidth}
\centering
\def\arraystretch{1.5}
\begin{tabular}{r|ccccccccc}
$H^a$&$ \frac{1}{64} $&$ 0 $&$ -\frac{1}{8} $&$ \frac{1}{4} $&$ \frac{23}{32} $&$ \frac{1}{4}  $&$ -\frac{1}{8} $&$ 0 $&$ \frac{1}{64} $ \\
$G^a$ & &$\frac{1}{16} $&$ 0 $&$ -\frac{9}{16} $&$ 1 $&$ -\frac{9}{16} $&$ 0 $&$ \frac{1}{16}$\\
\hline
$J^s$ & & & &$ -\frac{1}{4} $&$ 1 $&$ -\frac{1}{4} $\\
$K^s$ &$ \frac{1}{64} $&$ -\frac{1}{16} $&$ -\frac{1}{8} $&$ \frac{9}{16} $&$ \frac{23}{32} $&$ \frac{9}{16} $&$ -\frac{1}{8} $&$ -\frac{1}{16} $&$ \frac{1}{64} $
\end{tabular}
\subcaption{\wave{4}{2}}
\label{table_coef_interp_4_2}
\end{minipage}

\begin{minipage}{\textwidth}
\def\arraystretch{1.5}
\centering
\begin{tabular}{r|ccccccccccccc}
$H^a$ & $-\frac{3}{1024}$ & $0$ & $\frac{11 }{ 512}$ & $0$ & $-\frac{125 }{ 1024}$ & $\frac{1}{4}$ & $\frac{181 }{ 256}$ & $\frac{1 }{ 4}$ & $-\frac{125 }{ 1024}$ & $0$ & $\frac{11 }{ 512}$ & $ 0$ & $-\frac{3 }{ 1024}$ \\
$G^a$ & & $-\frac{3}{256}$ & $0$ & $\frac{25}{256}$ & $0$ & $-\frac{75}{128}$ & $1$ & $-\frac{75}{128}$ & $0$ & $\frac{25}{256}$ & $0$ & $-\frac{3}{256}$\\
\hline
$J^s$ & & & & & &$ -\frac{1}{4} $&$ 1 $&$ -\frac{1}{4} $\\
$K^s$ & $-\frac{3}{1024}$&$ \frac{3}{256}$&$ \frac{11}{512}$&$ -\frac{25}{256}$&$ -\frac{125}{1024}$&$ \frac{75}{128}$&$ \frac{181}{256}$&$ \frac{75}{128}$&$ -\frac{125}{1024}$&$ -\frac{25}{256}$&$ \frac{11}{512}$&$ \frac{3}{256}$&$ -\frac{3}{1024}$
\end{tabular}
\subcaption{\wave{6}{2}}
\label{table_coef_interp_6_2}
\end{minipage}

\caption{Lifted interpolating wavelet filter coefficients as operated on interleaved scaling and detail coefficients. The filters $H^a$ and $J^s$ are centered on (even) scaling coefficients, and the filters $G^a$ and $K^s$ are centered on odd scaling coefficients and detail coefficients, respectively.}
\label{table_coef_lifted}
\end{table}

%% file: wavelet_s4_algo.tex
\section{Algorithms and implementation}
\label{sec_sm_algo}
In this section we describe in detail the different algorithms as implemented in our presented software framework. We refer to \sect{sect_implementation_global} for additional explanations and context regarding the different steps.
In \cref{algo:ghost_1} and \cref{algo:ghost_2} we detail the implementation of the ghosting procedure, while in \cref{algo:adapt} we expose the implementation of the adaptation of the grid.

\subsection{Ghost computation}
We describe here the technicalities of the ghost computation implementation.
To facilitate the discussion, we introduce some notation inspired by the \code{p4est} naming convention.
Any block $\block$ is part of the total set of blocks $\Block$ within the grid, which are distributed among multiple MPI ranks. In this section we understand $\block$ to contain the cartesian grid data of size $N_b^3$, as well as a ghost region on the same level as $\block$ extending $N_g$ points in each dimension on each side of the block, where $N_g$ is determined from the PDE requirements or wavelet support, depending on the operation. A subset of $\Block$ is the mirror group $\Mirror$, defined as the set of blocks that are the neighbor of at least one other block on another rank, and thus have to be accessed through MPI communications. For each block $\block$, we further define $6$ different neighbor sets. All neighbors that exist on the same level as $\block$ fall in the group $\Ghost^0$, all neighbors that are coarser than $\block$ fall in the group $\Ghost^-$, and all neighbors that are finer than $\block$ fall in the group $\Ghost^+$. Each of these group can either exist on the same rank as $\block$, in which case they are denoted as `local' $\Ghost^{0/-/+}_L$, or on another rank in which case they are denoted as `global' $\Ghost^{0/-/+}_G$. The ghost point computation is divided into two parts. For a block $\block$, the first step combines the values of coarser neighbors $\Ghost^{-}_{G/L}$ and same level neighbors $\Ghost^{0}_{G/L}$ with a coarse perspective of the data in $\block$ to be able to compute a local refinement of the ghost region and retrieve the required ghost points (see \algref{algo:ghost_1}). This step implements the process outlined in the main section in \sect{sec_ghost_reconstruction}. In the second step, block $\block$ uses the refined region of step 1 (including ghosts) to coarsen again, leading to coarse ghost points that block $\block$ then communicates to neighbors $\Ghost^{-}_{G/L}$ (see \algref{algo:ghost_2}). 

\paragraph{Part~1: Coarser and same-level neighbors }
First, due to the continuity in memory needed by the RMA window (our current memory allocation does not guarantee continuity for a given field across blocks), we start by a copy of all the $\block \in \Mirror$ into the buffer reserved for the communication, $\bufferrma$. The ghost values from same-level neighbors are trivially obtained by accessing the values needed by the ghost region through a memory copy for $\Ghost^0_L$ and using \FnGet{} for $\Ghost^0_G$. 

To obtain ghost values from coarser and same-level neighbors of our current block $\block$, we have to use the refinement operation as described in \sect{sec:block_struct}. This operation relies on the grid values of $\Ghost^{-}$ and the even values of our current block $\block$ and $\Ghost^{0}$. %
We proceed by gathering all the required values for the computation of the ghost into a temporary buffer $\tempblock$ (whose size is $1/8^{\text{th}}$ of the block-size, extended with a coarsened representation of the ghost region). The contribution of the coarser neighbors are first obtained through a copy or \FnGet{} of the required values into $\tempblock$, and we gather the even values of the same-level neighbors, together with the even data of the current block. The data inside $\tempblock$ associated with fine-level neighbors is left blank and ignored throughout this process. Once all the needed data are gathered into $\tempblock$, we apply any domain boundary conditions at the coarse level. Finally, we perform the refinement operation locally on the rank associated with $\block$ and read out the computed ghost values within the regions overlapping with $\Ghost^{0}$ and $\Ghost^{-}$ from the refined buffer.

\paragraph{Part 2: Fine-level neighbors }
In the second part of the algorithm, we use the data in $\block$ and its ghost points computed in part~1 to compute coarse-level ghost points for $\Ghost^{-}$, with the individual steps described in \algref{algo:ghost_2}.
We first apply any specified domain boundary conditions on $\block$ if required, and then perform the substitution step as detailed in \sect{sec_mr_substitution}. We then coarsen $\block$ along the lines described in \sect{sec:block_struct} and store the coarsened data into $\tempblock$. Once the coarsening is completed, block $\block$ copies the required coarse-level ghost points to $\Ghost^{-}_L$ or issues a \FnPut{} for $\Ghost^{-}_G$ to write into the coarser neighbor's memory. When completed, we are left with the application of the boundary conditions as all the ghost informations are now complete.

With these two steps, all block interfaces in the domain can be handled. In the way explained above, the ghost reconstruction would correspond to a tree traversal from the coarsest to the finest levels in part~1, and back up to the coarsest level in part~2. However, choosing $N_b$ larger than the support of the wavelet filters together with the 2:1 constraint in resolution jumps guarantees that our algorithm for ghost reconstruction parallelizes efficiently across the different levels, with only independent synchronizations for each block in the  \FnWait{} call of \algref{algo:ghost_1} independent of the level on which it exists.

\begin{algorithm2e}[ht!]
\caption{Ghosting - Part 1}
\label{algo:ghost_1}
\SetAlgoLined
\DontPrintSemicolon
$\bufferrma$ $\leftarrow$ \FnCopy{$\block \in \Mirror$}

\BlankLine\BlankLine

\FnPost{}\; 
\FnStart{}\;
\BlankLine
\ForEach{$\block \in \Block$}{
    $\block$ $\leftarrow$ \FnCopy{$\Ghost_L^0$}, \FnGet{$\Ghost_G^0$}\;
    
    \BlankLine
    
    $\tempblock$ $\leftarrow$ \FnCopy{$\Ghost_L^-$}, \FnSample{$\Ghost_L^0$}\;
    $\tempblock$ $\leftarrow$ \FnGet{$\Ghost_G^-$}, \FnGet{\FnSample{$\Ghost_G^0$}}\; 
}
\BlankLine
\FnComplete{}\;
\FnWait{}\;
\BlankLine\BlankLine
\ForEach{$\block \in \Block$}{
    $\tempblock$ $\leftarrow$ \FnSample{$\block$}\;
    $\tempblock$ $\leftarrow$ \FnBounday{}\;
    $\block$ $\leftarrow$ \FnRefine{$\tempblock$}\;
}
\BlankLine\BlankLine
\end{algorithm2e}

\begin{algorithm2e}[ht!]
\caption{Ghosting - Part 2}
\label{algo:ghost_2}
\SetAlgoLined
\DontPrintSemicolon
\FnPost{}\; 
\FnStart{}\;

\BlankLine

\ForEach{$\block \in \Block$}{
$\block$ $\leftarrow$ \FnBounday{}\;
$\block$ $\leftarrow$ \FnOverwrite{}\;
$\tempblock$ $\leftarrow$ \FnCoarsen{b}\;
$\Ghost^-_L$ $\leftarrow$ \FnCopy{$\tempblock$}\;
$\Ghost^-_G$ $\leftarrow$ \FnPut{$\tempblock$}\;
}
\BlankLine
\FnComplete{}\;
\FnWait{}\;

\BlankLine\BlankLine
\ForEach{$\block \in \Block$}{
$\block$ $\leftarrow$ $\bufferrma$\;
$\block$ $\leftarrow$ \FnBounday{}\;
}
\end{algorithm2e}

\subsection{Grid adaptation}
The adaptation as detailed in \algref{algo:adapt} is an iterative procedure starting on a grid $\Block_k$, where the iteration $k$ is divided in multiple steps:
\begin{enumerate}
\item Ghost values computation and detail computation: the ghost values are updated dimension by dimension first on the field used for the detail computation and then on the rest of the fields present on the grid. %
\item In \FnCriterion{$\block \in \Block_k$} the detail coefficients are communicated and computed, either directly for scalar fields or on a component-by-component basis for vector and tensor fields.%
\item Every block gets a corresponding status that expresses its desired evolution given the detail coefficient values and the refinement $\epsilon_r$ and coarsening $\epsilon_c$ tolerances:
\begin{enumerate}
    \item \code{M\_ADAPT\_FINER}: if any of the detail coefficient is $\geq \epsilon_r$ for any dimension,
    \item \code{M\_ADAPT\_COARSER}: if every detail coefficient is $\leq \epsilon_c$ for every dimension,
    \item \code{M\_ADAPT\_SAME}: if no other status is assigned.
\end{enumerate}
\item Once all blocks have their desired status assigned, we check compliance against our grid adaptation policy to ensure that the grid adaptation is valid and consistent with the multiresolution theory; if needed, we change the status of individual block to enforce compliance. The policy consists of the following rules:
\begin{enumerate}
    \item we give priority to the refinement: if a block's coarser neighbor wants to refine, the latter cannot be coarsened;
    \item the need of a fine block to be refined has to be ``propagated" to its coarser neighbor: if a block has a finer neighbor which wants to refine, the latter must refine first to not break the 2:1 condition; similarly, if a block has a coarser neighbor, that neighbor has to be refined first;
    \item if a block has been refined in a previous iteration, it cannot be coarsened.
\end{enumerate}
In addition, we enable the option to enforce possible user-defined bounds on the minimum and/or maximum level permitted, if needed; these bounds are not used in the results of this work.
\item Once each block has received its final status, we use the \code{p4est} library to create the meta-structure associated with a new set of blocks $\Block_k^*$.
\item Refine or coarsen the blocks as needed to obtain the updated grid $\Block_{k+1}$ through the \FnInterpolate{} function.
\item Partition the grid using \code{p4est} .
\item Synchronize the executed adaptation step of each block to its neighbors.
\item Update the adapted fields in \FnSmooth{$\block \in \Block_{k+1}$} to enforce that fine blocks receive updated scaling coefficients if one or more neighbors have been coarsened (\sect{sec_mr_coarsening}).%
\item Obtain the number of adapted blocks over the whole grid (implemented through a non-blocking reduce operation). If no blocks have been adapted, we have reached the final grid with $\normi{\gamma} \leq \epsilon_r$. Else, move to the next iteration $k+1$. 
\end{enumerate}

Throughout these steps, we rely on \code{p4est} functions to support the coarsening, the refinement, and the partitioning of the grid; however, this is mostly limited to handling the metadata, as we have implemented the block refinement, coarsening and partitioning ourselves to exploit asynchronous and non-blocking MPI calls. 
\begin{algorithm2e}
\caption{grid adaptation}
\label{algo:adapt}
\SetAlgoLined
\DontPrintSemicolon
\While{adapt}{
    \FnGhost{$\block \in \Block_k$}\;
    \FnCriterion{$\block \in \Block_k$}\;
    \FnSyncStatus{}\;
    \FnPolicy{$\block \in \Block_k$}\;
    \BlankLine
    \FnTreeCoarsen{}\;
    \FnTreeRefine{}\;
	$\Block_{k+1} \leftarrow$ \FnInterpolate{$\block^* \in \Block^*_{k+1}$}\;
    \FnTreePartition{}\;
    \BlankLine

	\FnSyncStatus{}\;
	\FnSmooth{$\block \in \Block_{k+1}$}\;
	\BlankLine
	$\Block_k \leftarrow \Block_{k+1}$
}
\end{algorithm2e}

%% file: wavelet_s3_rkandweno.tex
\section{Time integration scheme}
\label{sec_rk3}
As a time integrator we are using the RK3-TVD \citep{Gottlieb:1998,Gottlieb:2001} also know as RK3-SSP.
This is an explicit in time, 3 step Runge-Kutta scheme, which only requires two temporary buffers. 
The integration for $\dot{u} = f(t,u) $ from time $t^n$ to $t^{n+1}$ with $\Delta t = t^{n+1} - t^n$ is given by the three stage equations:
\begin{eqd}
y_1  = \Delta t f\left(t^n, u^n \right) + u^n ~,\qquad 
y_2  = \frac{1}{4} \left[ \Delta t f\left( t+\Delta\!t , y_1 \right) + y_1 \right] + \frac{3}{4} u^n ~,\qquad
u^{n+1}  = \frac{2}{3} \left[ \Delta t f\left(t+\dfrac{\Delta\!t}{2}  y_2 \right) + y_2 \right] + \frac{1}{3} u^n
\end{eqd}

\section{Finite differences scheme}
\label{sec_finite_diff}
For the spatial discretization of the derivatives in the advection equation, we rely on the \cons{3} scheme, which is a fixed-weight version of a \code{WENO} schemes.
First, given a divergence-free velocity field, the linear advection equation is considered in conservative form:
\be
\dfrac{\partial \phi}{\partial t} + \bu \cdot \nabla \phi = 0 \quad \Leftrightarrow \quad \dfrac{\partial \phi}{\partial t} +  \nabla \cdot \left( \bu \; \phi \right) = 0
\eec
with flux function $f = \bu \; \phi$. 
We consider a $1$D version of this equation, discretized on a uniform grid, \ic{} the grid cell $i$ spans $\left[ x_{i-1/2} \;;\; x_{i+1/2} \right]$:
\be
\dfrac{\partial \phi_i}{\partial t}  = -\left( f_{i+1/2} - f_{i-1/2} \right) \dfrac{1}{h}
\eec
where $f_{i \pm 1/2}$ is the flux at the cell boundary. Here we follow \citep[section 2.1]{Shu:1997} to derive a conservative finite-difference form for the flux reconstruction. 
To maintain stability \citep[section 2.1]{Shu:1997} the flux terms must be decomposed into a positive $f^+$ and negative part $f^-$ such that
\be
\dfrac{\partial f^{+}}{\partial \phi} > 0 \qquad \text{and} \qquad \dfrac{\partial f^{-}}{\partial \phi} < 0
\eed
For our equations in 1D we find $f^+ = u \phi$ if $u >0$ and $f^- = u \phi$ if $u <0$. Since the velocity information at the interface is not directly available, we use a simple reconstruction, $u_{i+1/2} = 1/2 \; \left( u_{i+1} + u_{i-1} \right)$. This approximation gives the sign of the velocity at the interface from the velocity in each grid cell, which determines the choice between $f^+$ and $f^-$ at each cell interface.

The actual flux computation can be done in different ways, leading to essentially non-oscillatory (ENO) and weighted ENO (WENO) interpolation formulas. Using three points in the flux calculation, one can choose between two stable stencils to compute $f^+_{i+1/2}$ 
\be
S_0: \quad f^+_{i+1/2} = -1/2 \; f^+_{i-1} + 3/2 \; f^+_{i} \qquad S_1: \quad f^+_{i+1/2} = 1/2 \; f^+_{i} + 1/2 \; f^+_{i+1}
\eed
Further, one can associate a smoothness indicator to each stencil $\beta_0 = \left(f^+_{i} - f^+_{i-1} \right)^2$ and $\beta_1 = \left(f^+_{i+1} - f^+_{i} \right)^2$. The standard ENO approach relies on the $\beta$ indicators to choose the best stencil, \ic{} we choose either $S_0$ or $S_1$ to evaluate the flux. The WENO technique instead combines the stencils together relying on weights associated to each of them ($w_0$ and $w_1$ respectively) in order to obtain the most accurate evaluation possible.

In the case of the stencils $S_0$ and $S_1$ one can combine them together to reach a third-order stencil using the ``ideal'' weights $w_0 = \gamma_0 = 1/3$ and $w_1 = \gamma_1 = 2/3$. By doing so, one obtains the conservative third-order stencil \cons{3} we used in the manuscript. Along similar lines, the fifth order conservative stencil \cons{5} can be derived from the optimal weights of \weno{5}. Tables \tref{tab:weno3} and \tref{tab:weno5} summarize the flux definitions and optimal weights $\gamma$ for each of these schemes. For non-smooth fields the weights can be adapted locally based on the $\beta$ smoothness indicators, and we implemented the WENO-Z version \cite{Borges:2008,Don:2013} in our code -- however for the results in this manuscript we only use the fixed-weight schemes.

\begin{table}[h!]
    \centering
    \renewcommand{\arraystretch}{2.0}
    \begin{tabular}{l|lll}
        \multirow{2}{*}{$f^+_{i+1/2}$}
        & $S^+_0 = -1/2 f_{i-1} \;+\; 3/2 f_{i}$ & $\beta^+_0 = \left(f_{i}-f_{i-1}\right)^2$ & $\gamma^+_0 = 1/3$ \\
        & $S^+_1 = 1/2 f_{i} \;+\; 1/2 f_{i+1}$ & $\beta^+_1 = \left(f_{i+1}-f_{i}\right)^2$ & $\gamma^+_1 = 2/3$ \\
    \hline
        \multirow{2}{*}{$f^-_{i-1/2}$}
        & $S^-_0 = 1/2 f_{i-1} + 1/2 f_{i}$ & $\beta^-_0 = \beta^+_0 = \left(f_{i}-f_{i-1}\right)^2$ & $\gamma^-_0 = 2/3$ \\
        & $S^-_1 = 3/2 f_{i} - 1/2 f_{i+1}$ & $\beta^-_1 = \beta^+_1 = \left(f_{i+1}-f_{i}\right)^2$ & $\gamma^-_1 = 1/3$ \\
    \end{tabular}
    \caption{WENO stencils - order 3. Setting the weights equal to the $\gamma$'s leads to the third-order \cons{3} scheme. }
    \label{tab:weno3}
\end{table}

\begin{table}[h!]
    \centering
    \renewcommand{\arraystretch}{2.0}
    \begin{tabular}{l|lll}
        \multirow{3}{*}{$f^+_{i+1/2}$}
        & $S^+_0 = 1/3 f_{i-2} \;-\; 7/6 f_{i-1} \;+\; 11/6 f_{i}$
            & $\beta^+_0 = \frac{1}{4} \left(f_{i-2}- 4 f_{i-1} + 3 f_{i}\right)^2 + \frac{13}{12} \left( f_{i-2} - 2 f_{i-1} + f_{i} \right)^2$ & $\gamma^+_0 = 1/10$ \\
        & $S^+_1 = -1/6 f_{i-1} \;+\; 5/6 f_{i} \;+\; 1/3 f_{i+1}$
            & $\beta^+_1 = \frac{1}{4} \left(-f_{i-1} + f_{i+1}\right)^2 + \frac{13}{12} \left( f_{i-1} - 2 f_{i} + f_{i+1} \right)^2$ & $\gamma^+_1 = 3/5$ \\
        & $S^+_2 = 1/3 f_{i} \;+\; 5/6 f_{i+1} \;-\; 1/6 f_{i+2}$
            & $\beta^+_2 = \frac{1}{4} \left(-3 f_{i} + 4 f_{i+1} - f_{i+2}\right)^2 + \frac{13}{12} \left( f_{i} - 2 f_{i+1} + f_{i+2} \right)^2$ & $\gamma^+_2 = 3/10$ \\
    \hline
        \multirow{2}{*}{$f^-_{i-1/2}$}
        & $S^-_0 = -1/6 f_{i-2} \;+\; 5/6 f_{i-1} \;+\; 1/3 f_{i}$
            & $\beta^-_0 = \beta^+_0 = \frac{1}{4} \left(f_{i-2}- 4 f_{i-1} + 3 f_{i}\right)^2 + \frac{13}{12} \left( f_{i-2} - 2 f_{i-1} + f_{i} \right)^2$ & $\gamma^-_0 = 3/10$ \\
        & $S^-_1 = 1/3 f_{i-1} \;+\; 5/6 f_{i} \;-\; 1/6 f_{i+1}$
            & $\beta^-_1 = \beta^+_1 = \frac{1}{4} \left(-f_{i-1} + f_{i+1}\right)^2 + \frac{13}{12} \left( f_{i-1} - 2 f_{i} + f_{i+1} \right)^2$ & $\gamma^-_1 = 3/5$ \\
        & $S^-_2 = 11/6 f_{i} \;-\; 7/6 f_{i+1} \;+\; 1/3 f_{i+2}$
            & $\beta^-_2 = \beta^+_2 = \frac{1}{4} \left(-3 f_{i} + 4 f_{i+1} - f_{i+2}\right)^2 + \frac{13}{12} \left( f_{i} - 2 f_{i+1} + f_{i+2} \right)^2$ & $\gamma^-_2 = 1/10$ \\
    \end{tabular}
    \caption{WENO stencils - order 5. Setting the weights equal to the $\gamma$'s leads to the fifth-order \cons{5} scheme. }
    \label{tab:weno5}
\end{table}

%% file: wavelet_s1_scalability.tex
\section{Weak scalability analysis}
\begin{figure}[ht!]
\begin{minipage}{0.32\textwidth}
	\centering
	\ifarxiv
	\includegraphics[width=\textwidth]{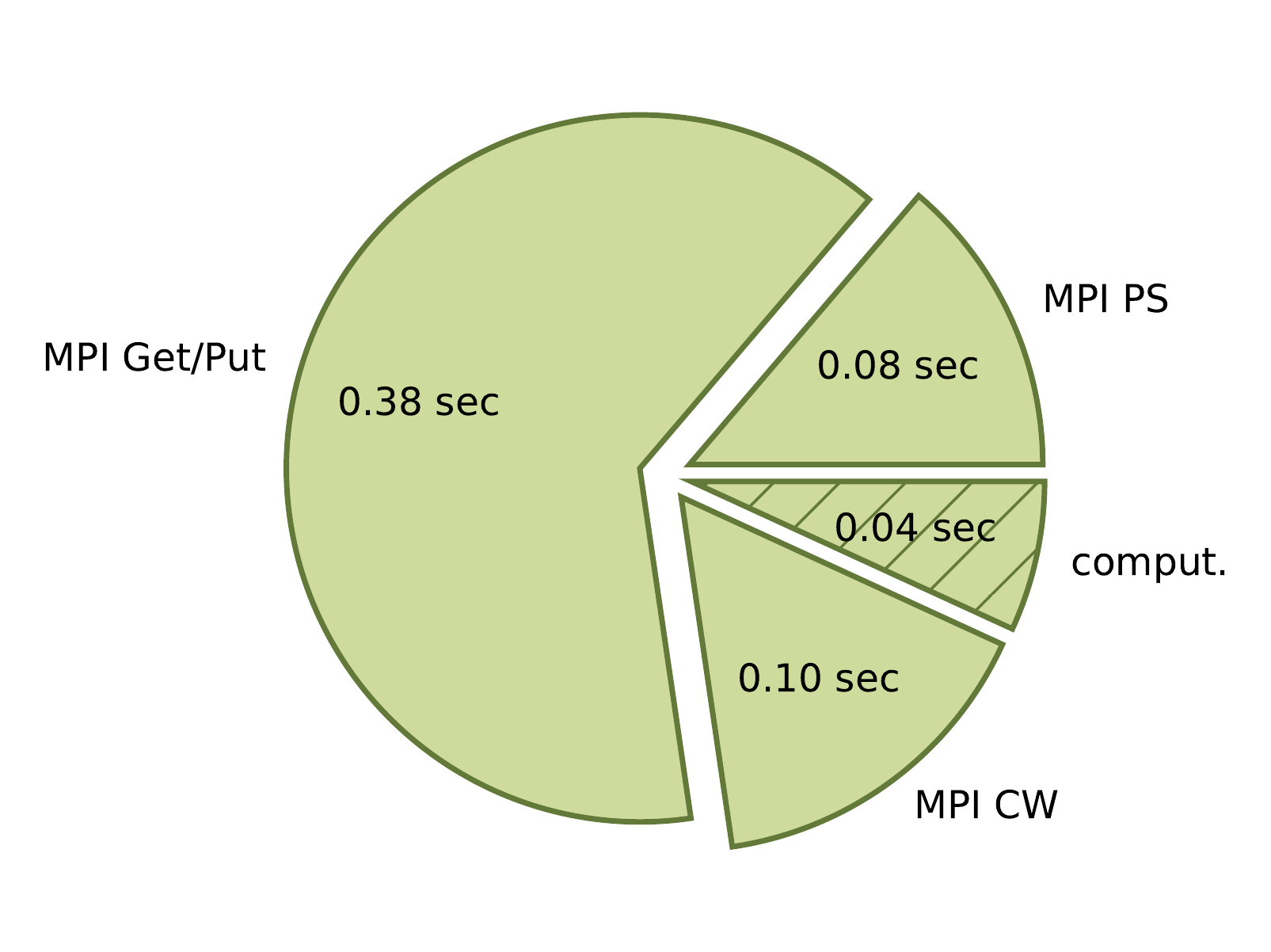}
	\else
	\includegraphics[width=\textwidth]{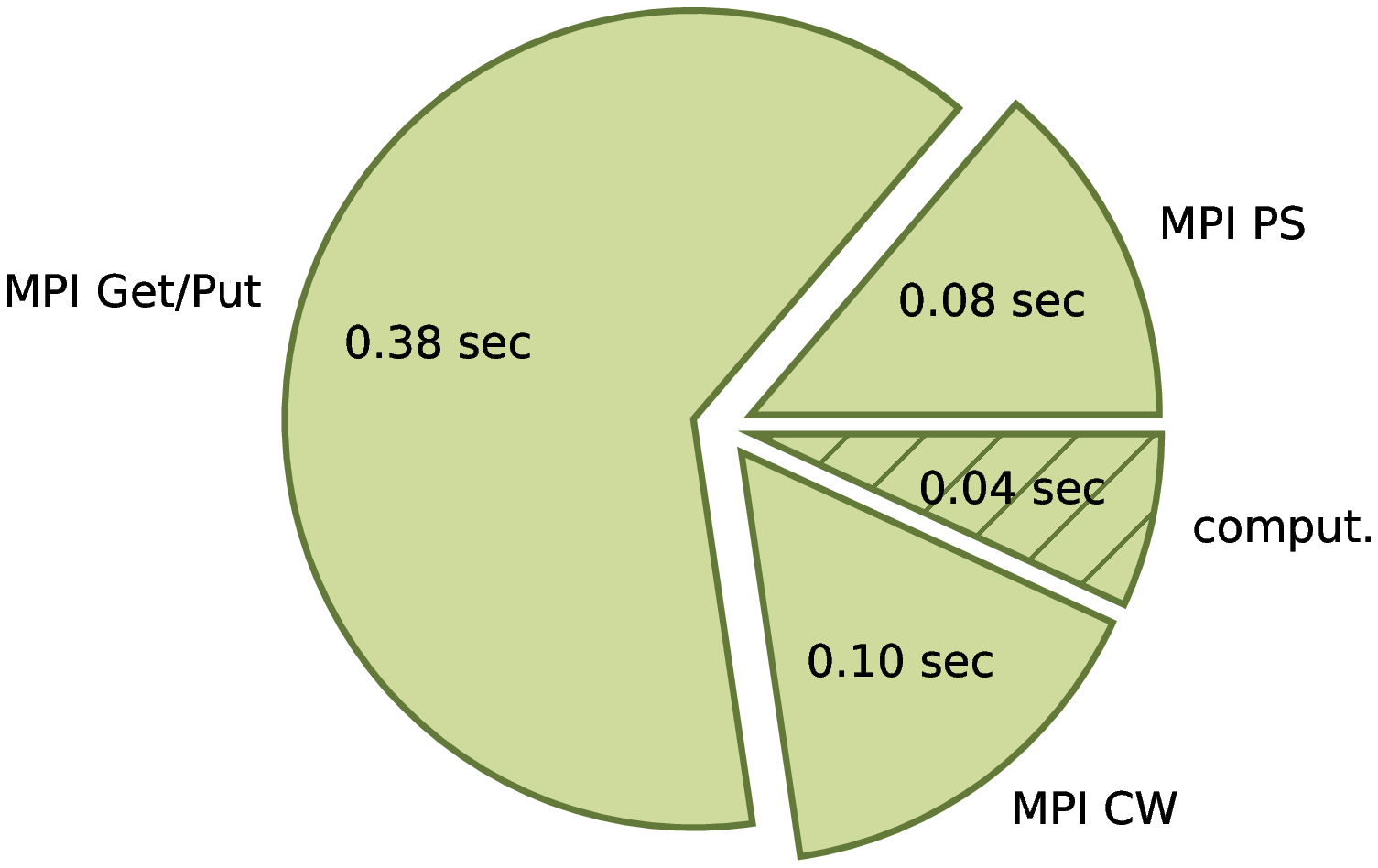}
	\fi
	\subcaption{\typo{$1$ node - $32$ ranks}}
\end{minipage}%
\hfill%
\begin{minipage}{0.32\textwidth}
	\centering
	\ifarxiv
	\includegraphics[width=\textwidth]{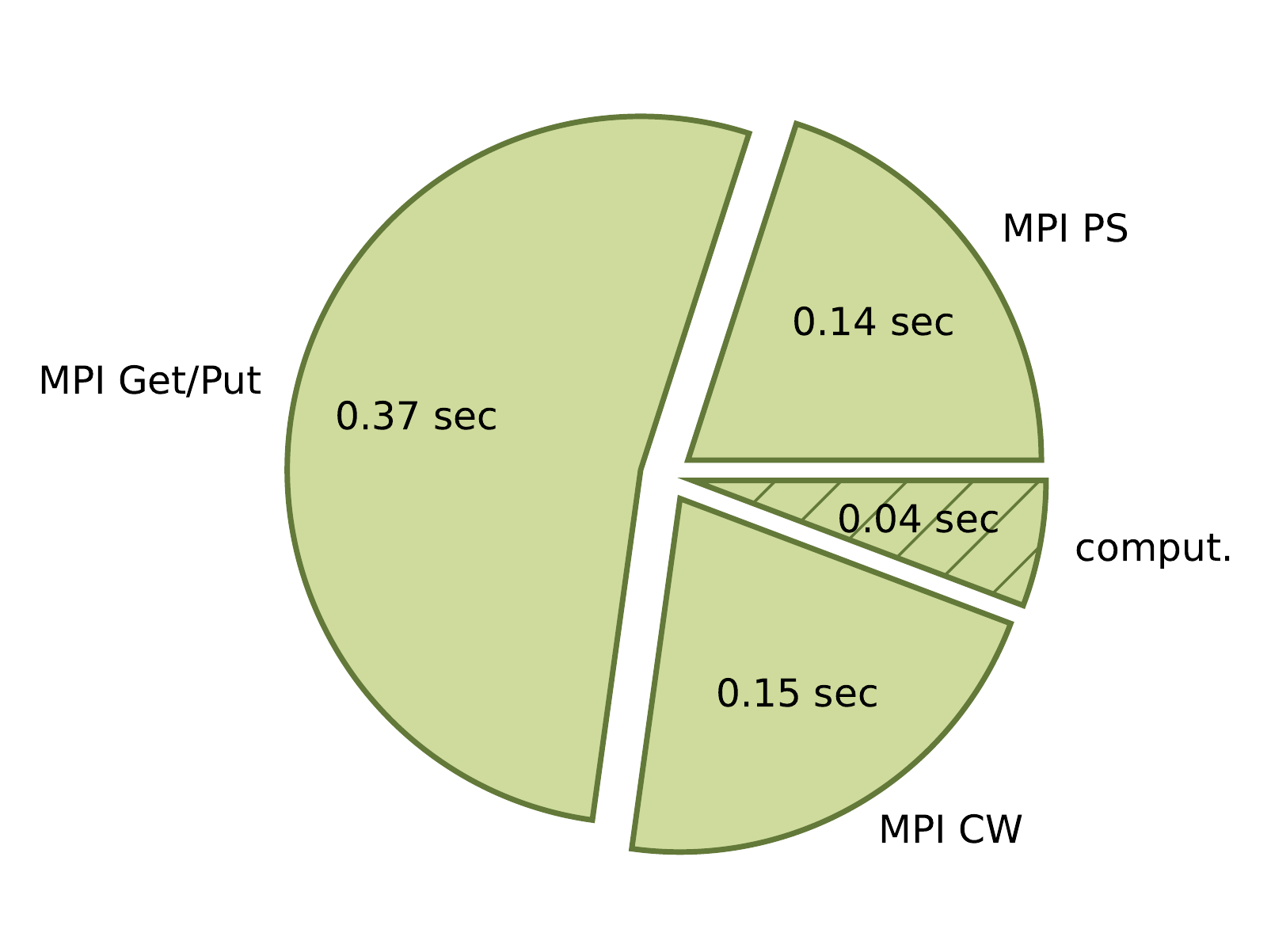}
	\else
	\includegraphics[width=\textwidth]{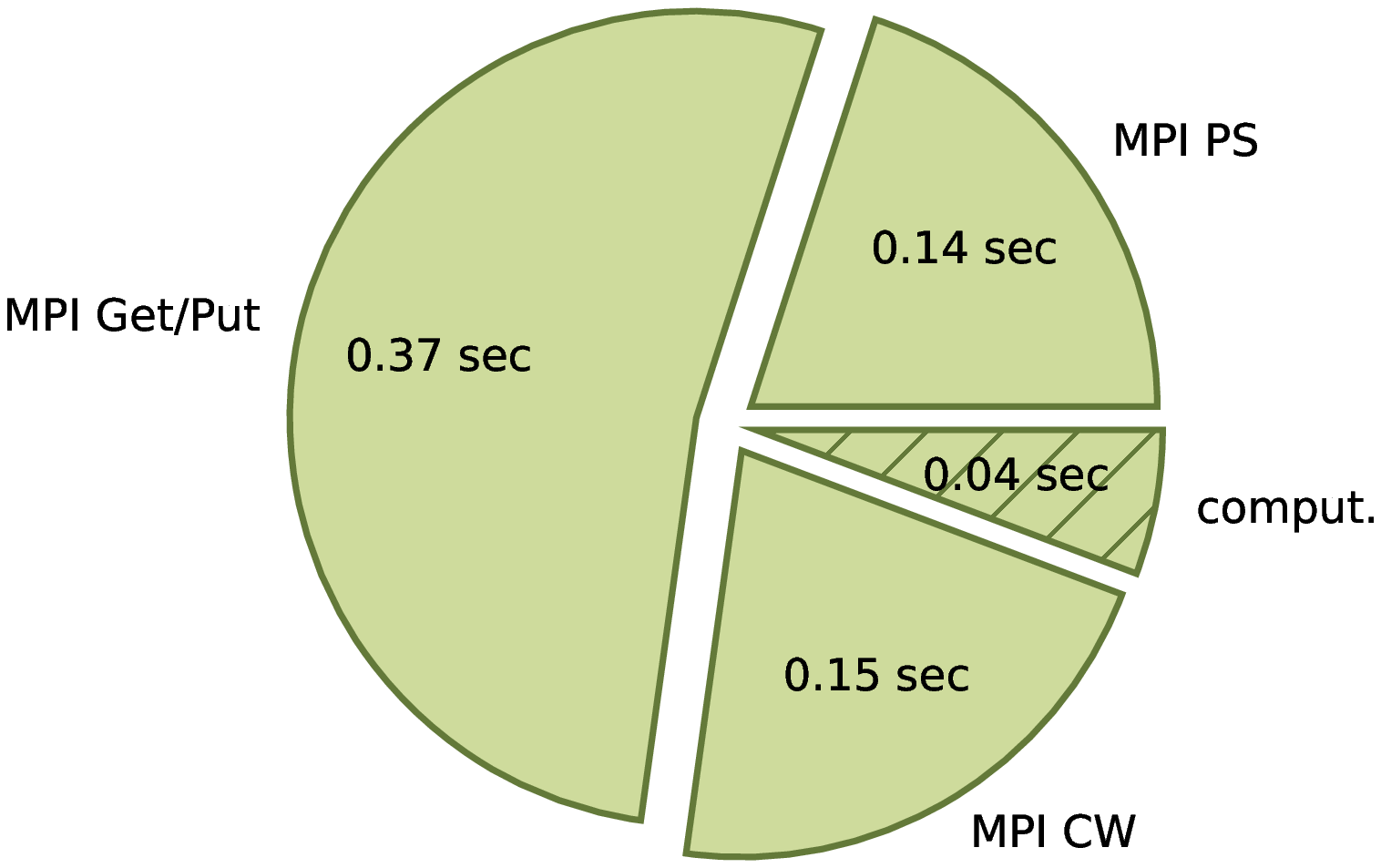}
	\fi
	\subcaption{\typo{$128$ node - $4,096$ ranks}}
\end{minipage}%
\hfill%
\begin{minipage}{0.32\textwidth}
	\centering
	\ifarxiv
	\includegraphics[width=\textwidth]{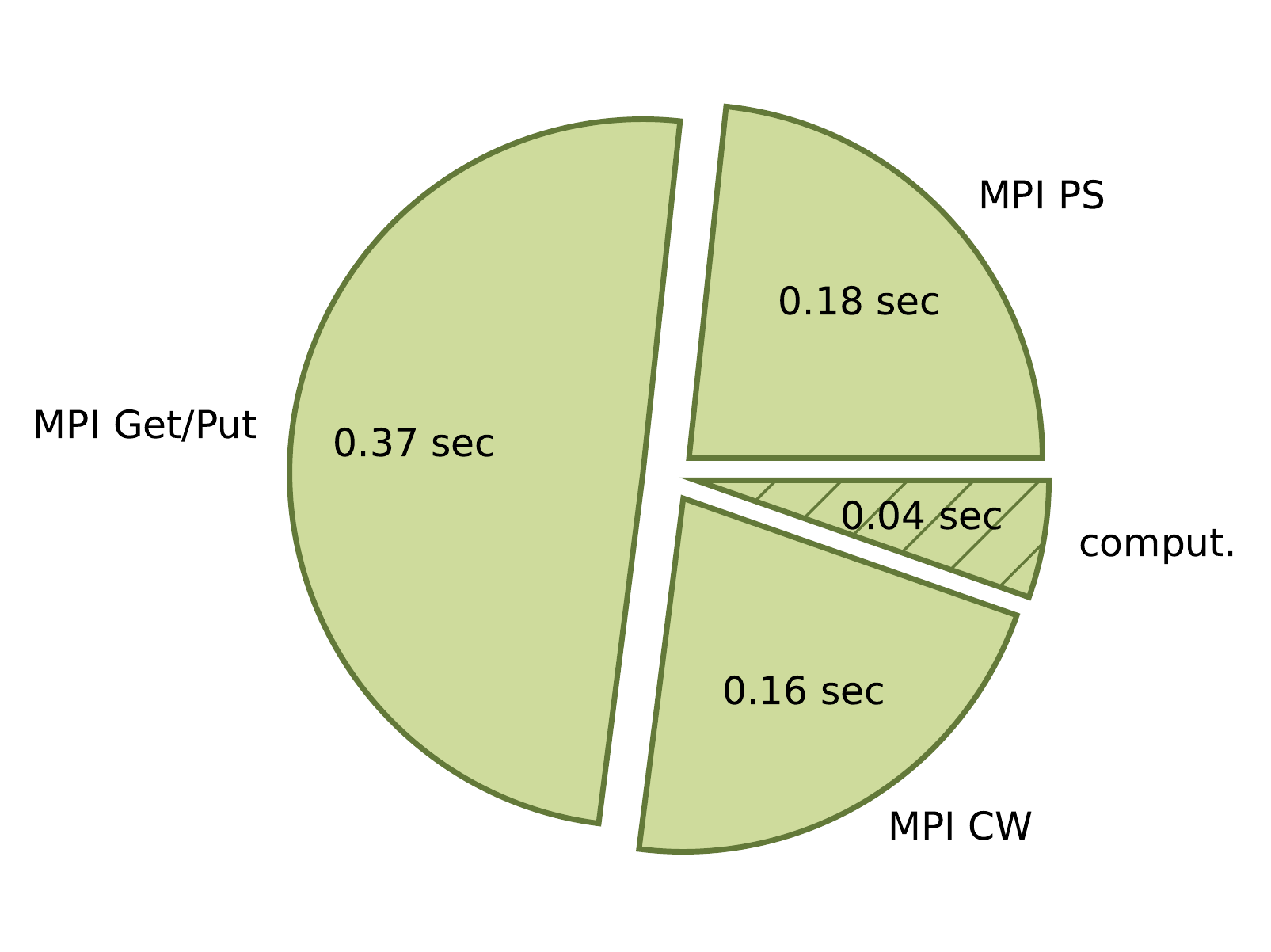}
	\else
	\includegraphics[width=\textwidth]{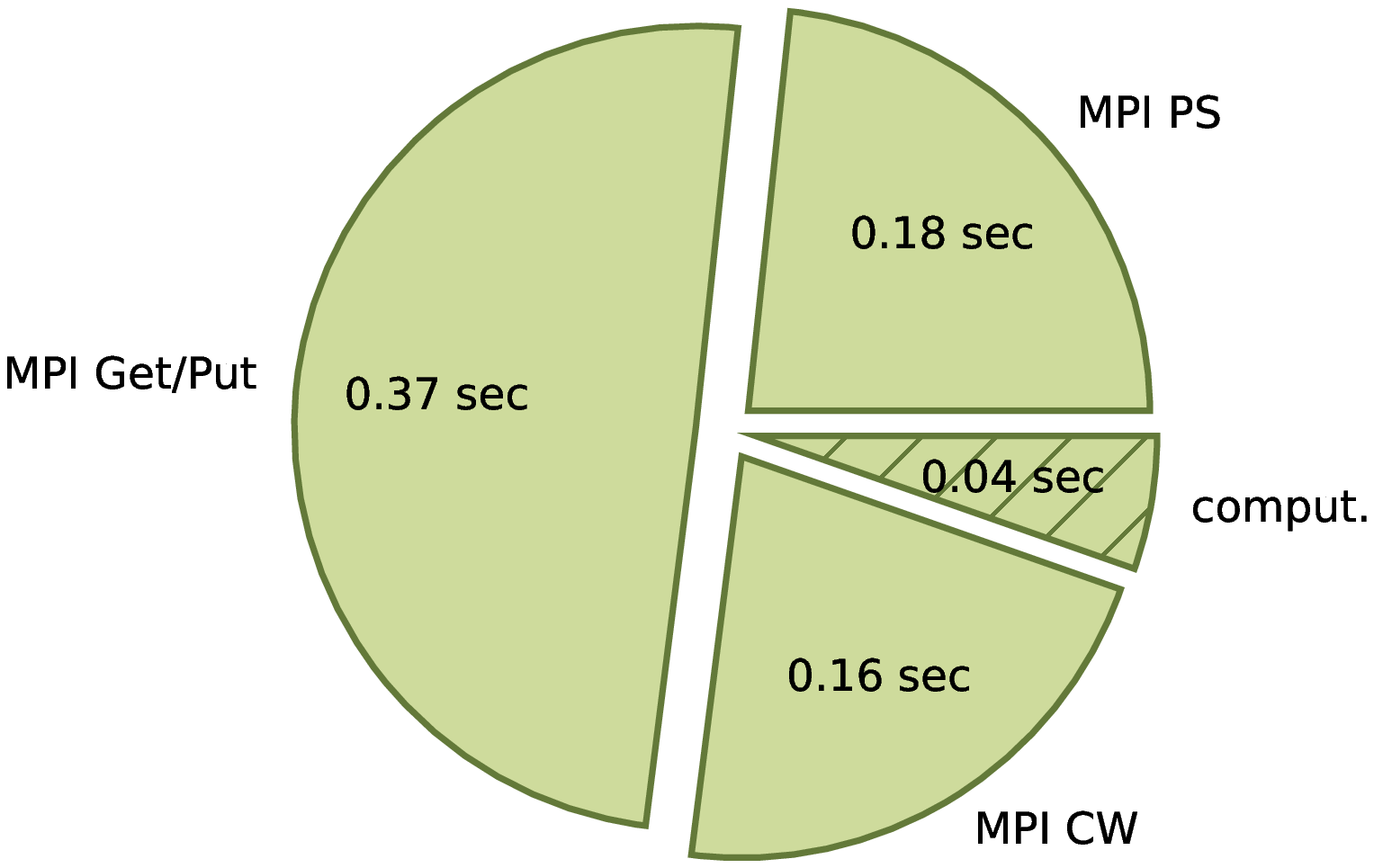}
	\fi
	\subcaption{\typo{$512$ node - $16,384$ ranks}}
\end{minipage}%
\hfill%
\caption{Breakdown of the time spent during ghost reconstruction in a weak scalability test.  The hatched region indicates operations only involving computations and no communications. The regions marked with ``\code{MPI PS}'' and ``\code{MPI CW}'' represent the calls to perform the active synchronization, ``\code{MPI Get/Put}`` contains the put and get accesses, and ``\code{comp.}'' is the time spent during the wavelet-based refinement and coarsening required to compute the ghost values.}
 \label{fig:weak_pie_ghost}
\end{figure}

\begin{figure}[ht!]
\begin{minipage}{0.32\textwidth}
	\centering
	\ifarxiv
	\includegraphics[width=\textwidth]{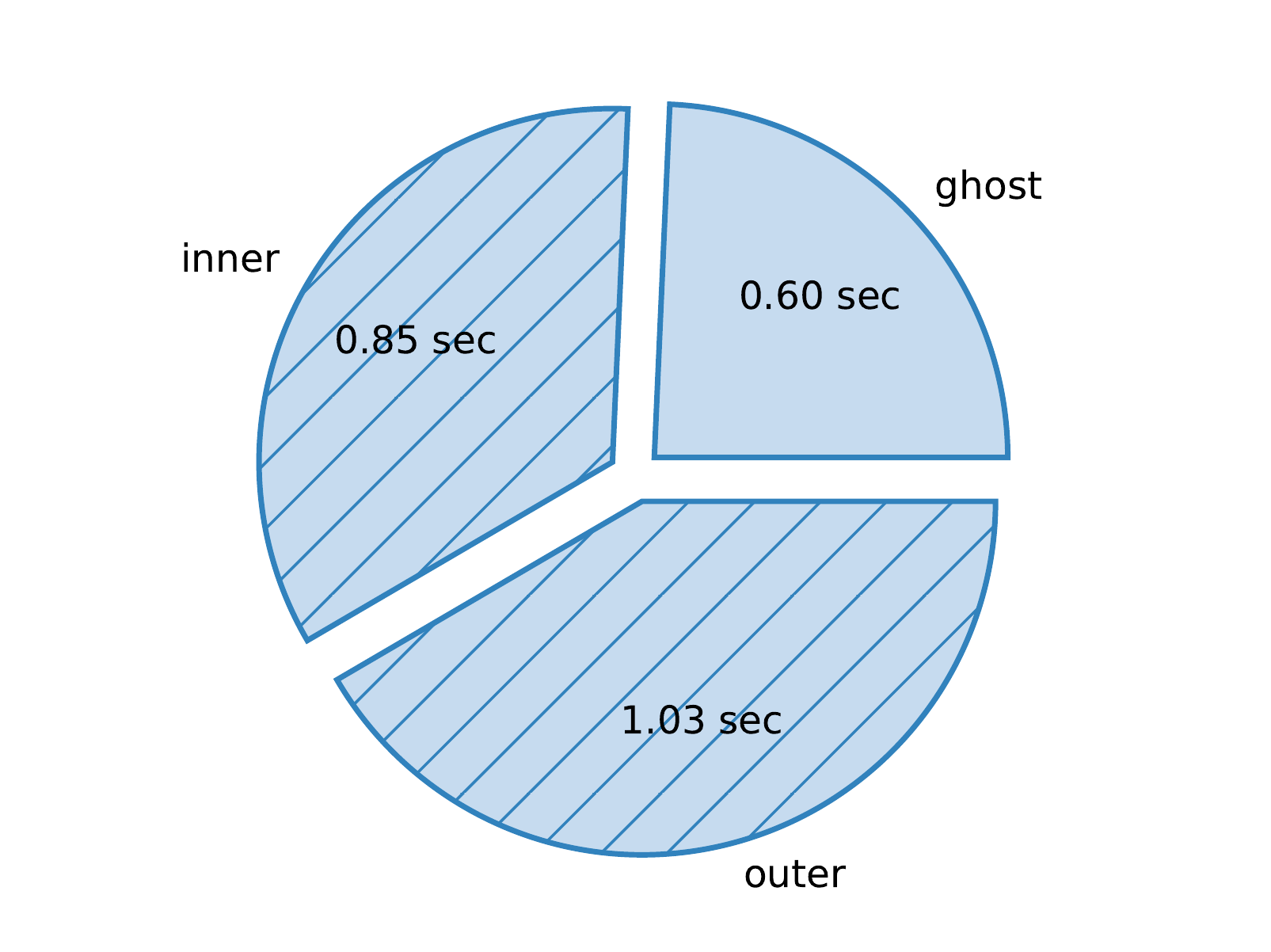}
	\else
	\includegraphics[width=\textwidth]{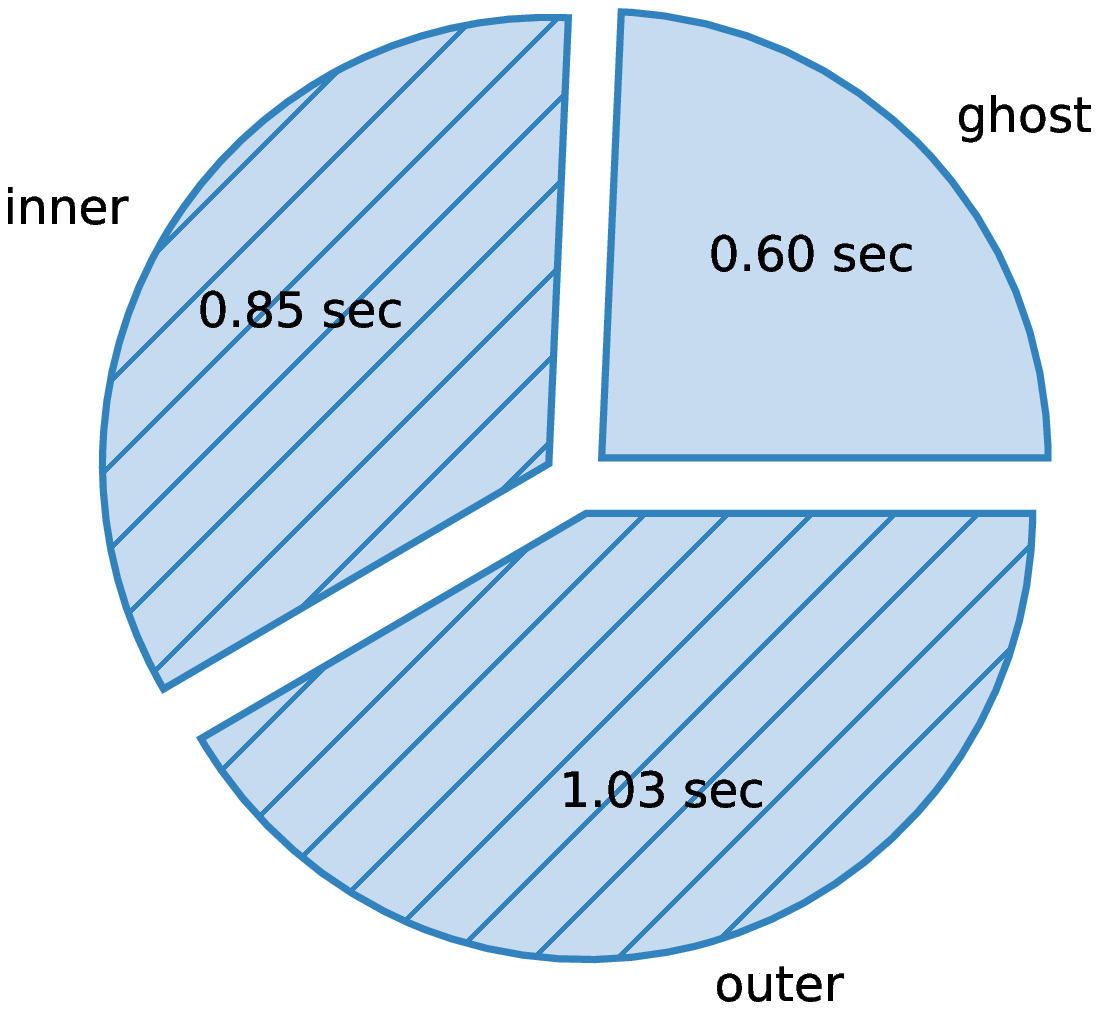}
	\fi
	\subcaption{$1$ node - $32$ ranks}
\end{minipage}%
\hfill%
\begin{minipage}{0.32\textwidth}
	\centering
	\ifarxiv
	\includegraphics[width=\textwidth]{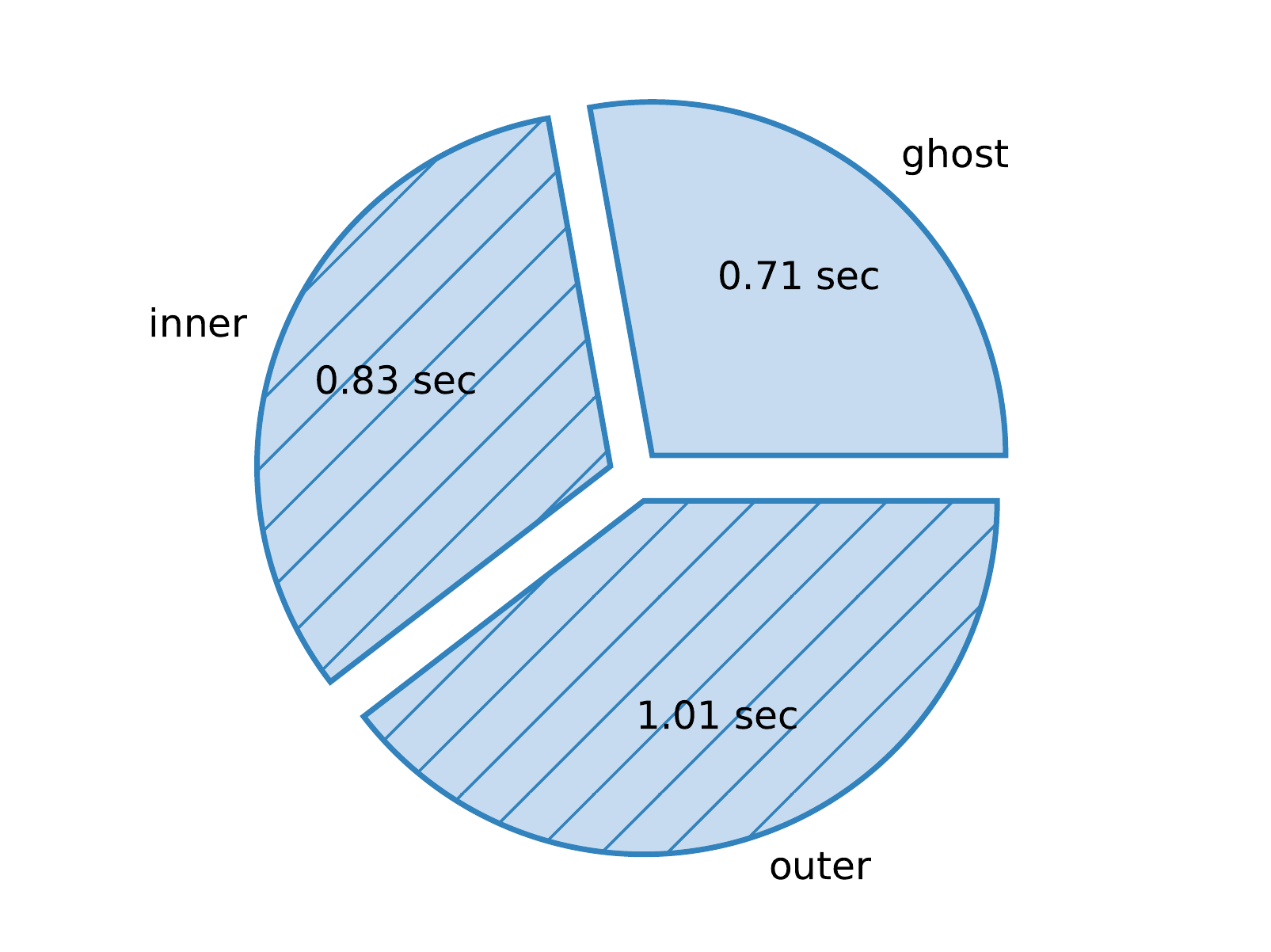}
	\else
	\includegraphics[width=\textwidth]{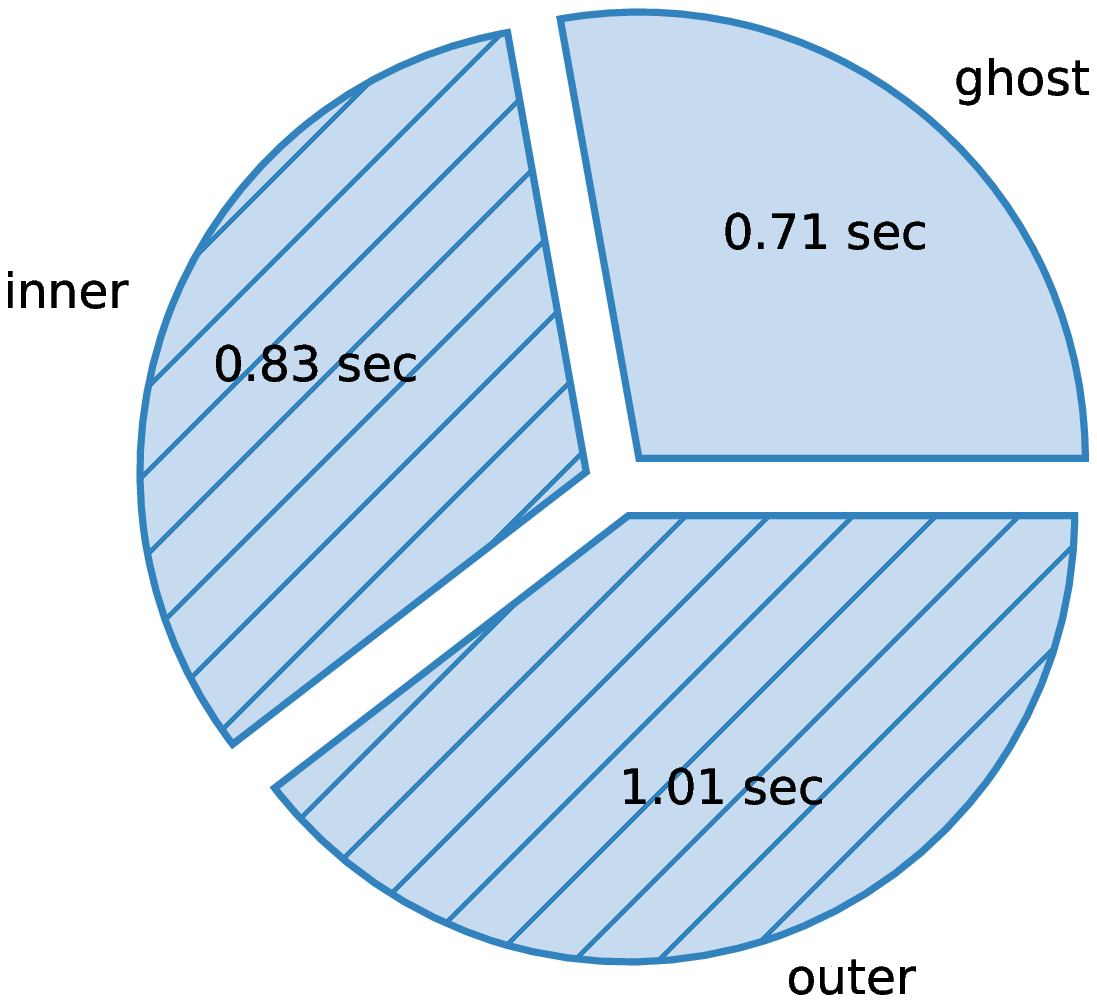}
	\fi
	\subcaption{$128$ node - $4,096$ ranks}
\end{minipage}%
\hfill%
\begin{minipage}{0.32\textwidth}
	\centering
	\ifarxiv
	\includegraphics[width=\textwidth]{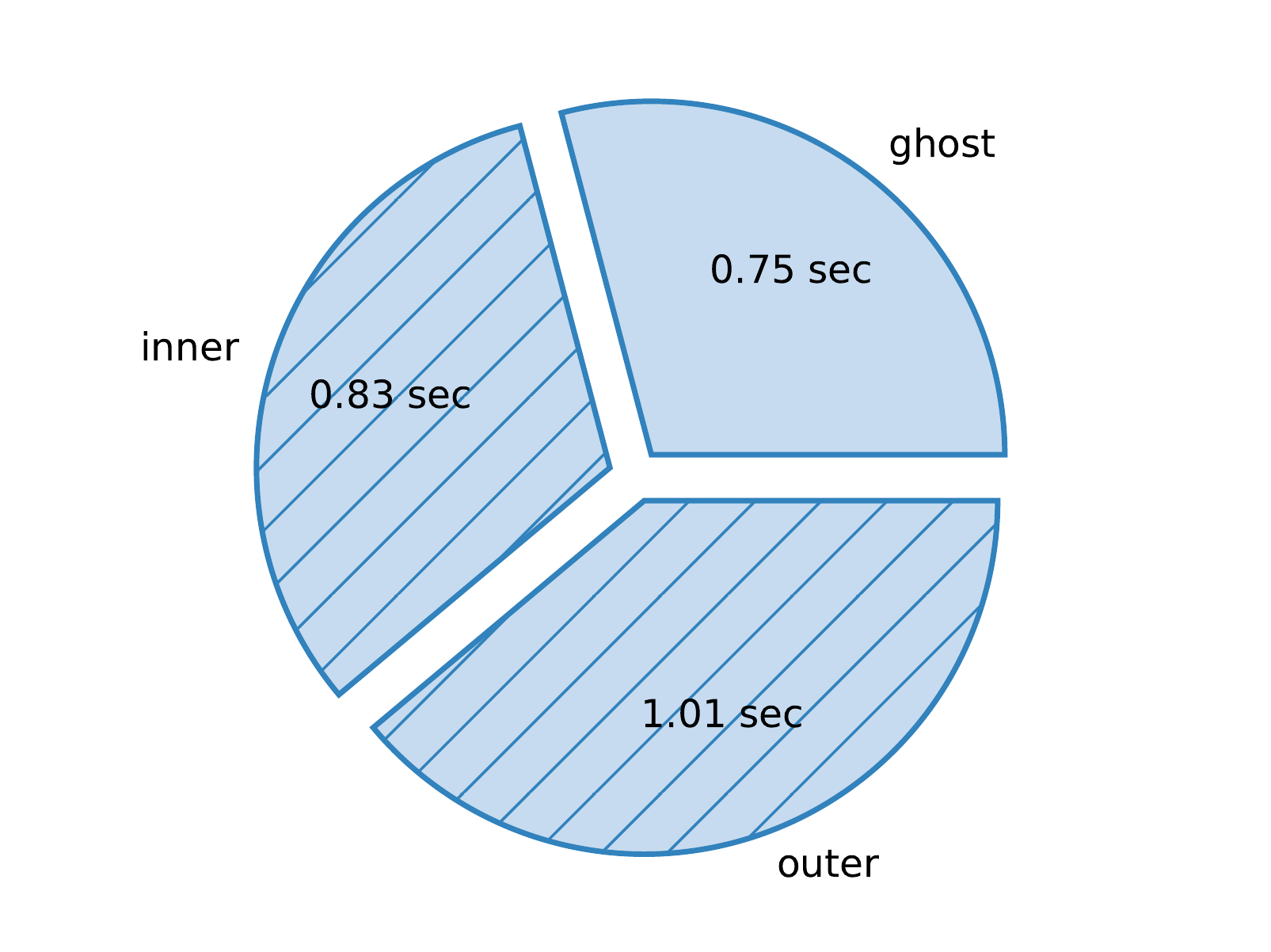}
	\else
	\includegraphics[width=\textwidth]{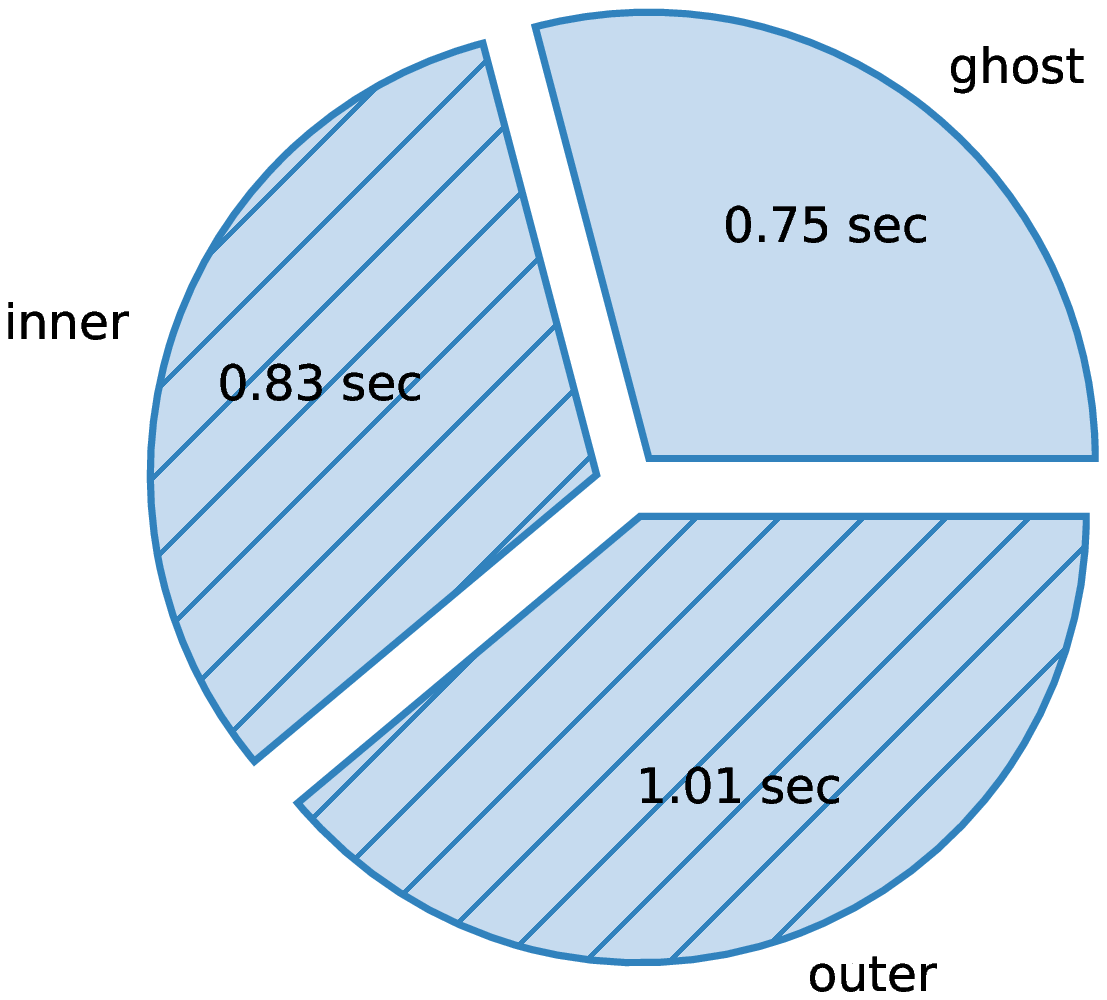}
	\fi
	\subcaption{$512$ node - $16,384$ ranks}
\end{minipage}%
\hfill%
\caption{Breakdown of the time spent during stencil computations in a weak scalability test. The hatched regions indicates operations only involving computations and no communications. The regions marked ``\code{inner}'' and ``\code{outer}'' are the computation of the stencil on the the inner and outer points in a block respectively, and ``\code{ghost}'' represents the computation of the ghost values overlapping with the ``\code{inner}'' computations \typo{as detailed in \fref{fig:weak_pie_ghost}}}
 \label{fig:weak_pie_stencil}
\end{figure}

\begin{figure}[ht!]
\begin{minipage}{0.32\textwidth}
	\centering
	\ifarxiv
	\includegraphics[width=\textwidth]{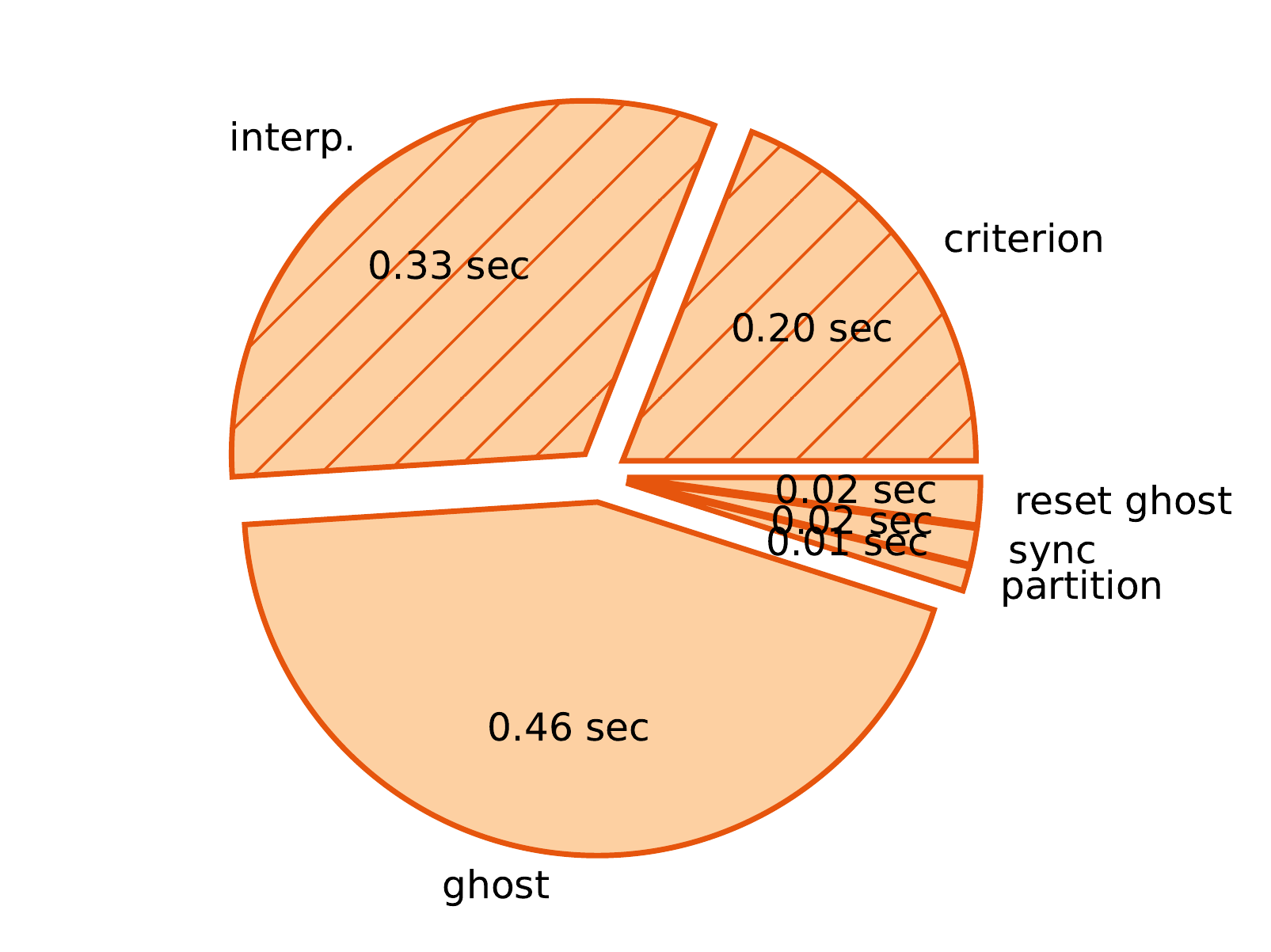}
	\else
	\includegraphics[width=\textwidth]{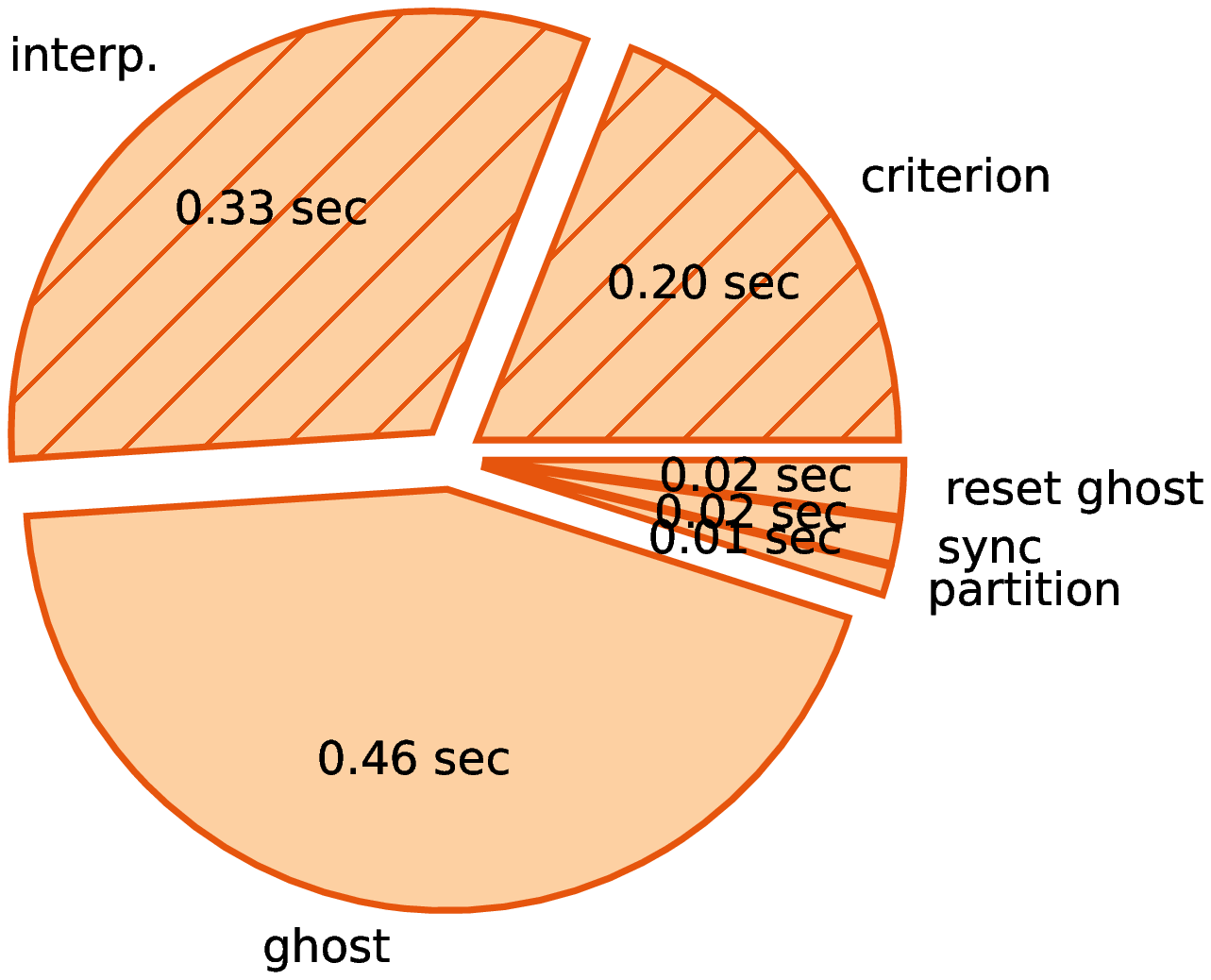}
	\fi
	\subcaption{$1$ node - $32$ ranks}
\end{minipage}%
\hfill%
\begin{minipage}{0.32\textwidth}
	\centering
	\ifarxiv
	\includegraphics[width=\textwidth]{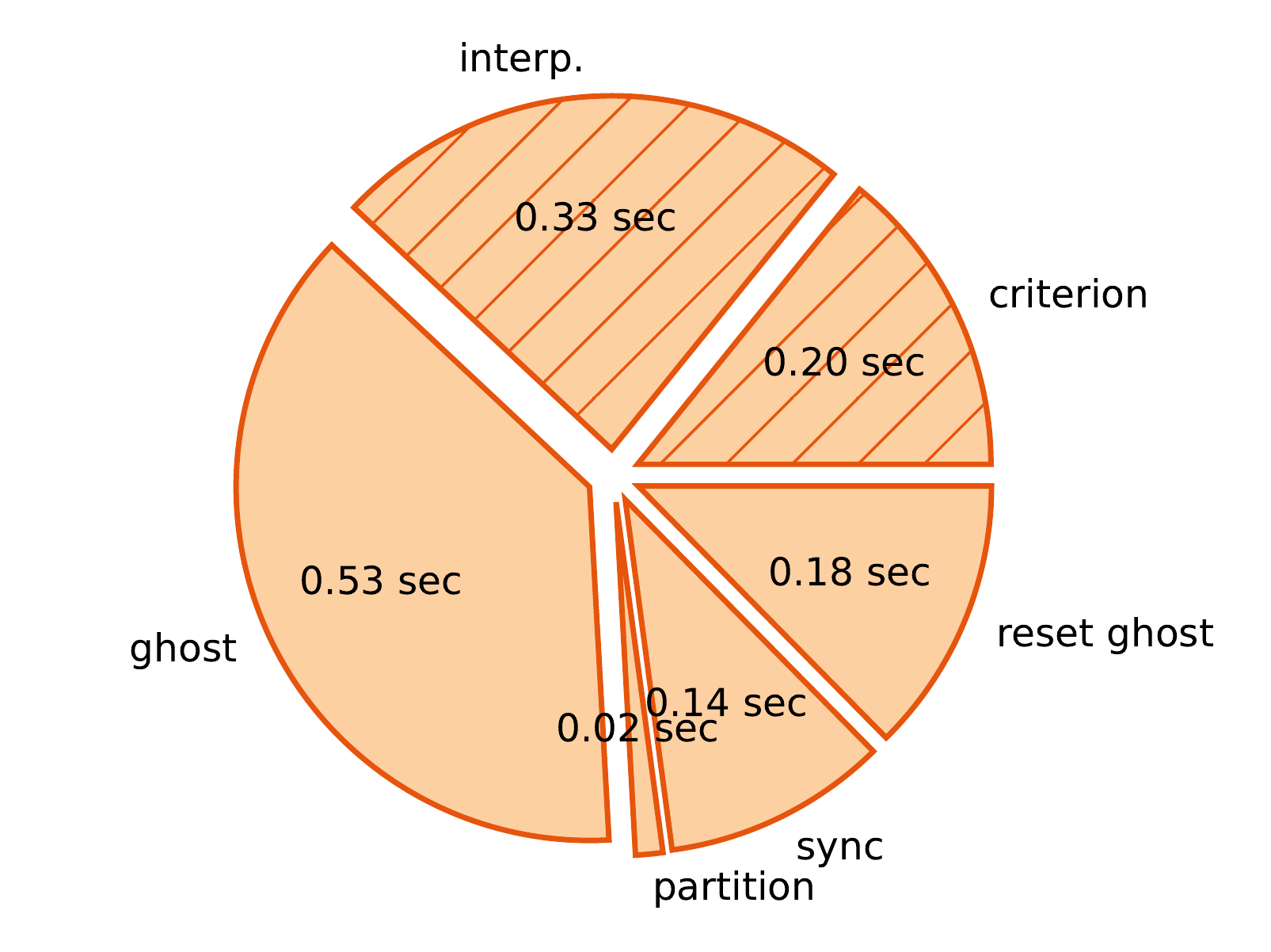}
	\else
	\includegraphics[width=\textwidth]{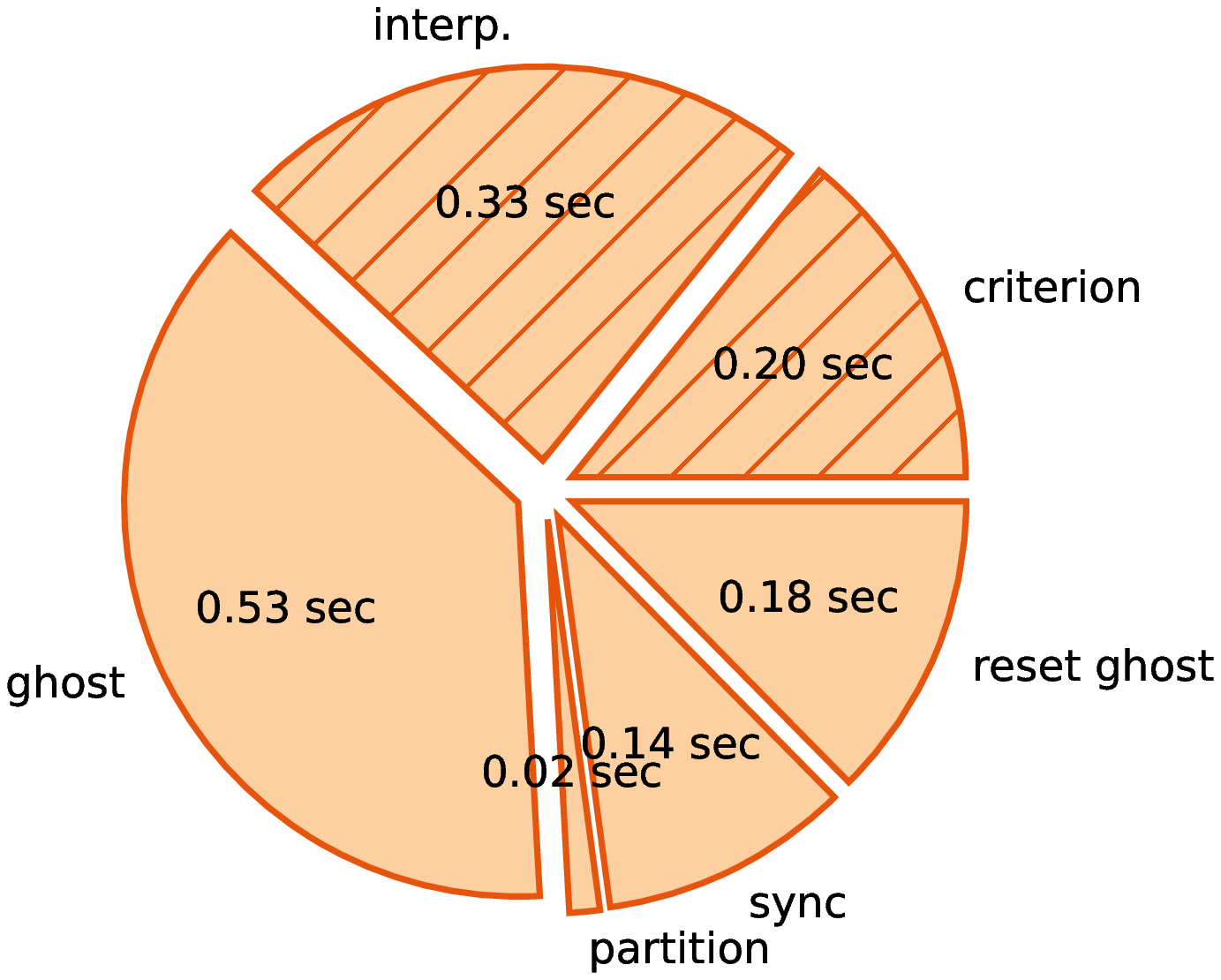}
	\fi
	\subcaption{$128$ node - $4,096$ ranks}
\end{minipage}%
\hfill%
\begin{minipage}{0.32\textwidth}
	\centering
	\ifarxiv
    	\includegraphics[width=\textwidth]{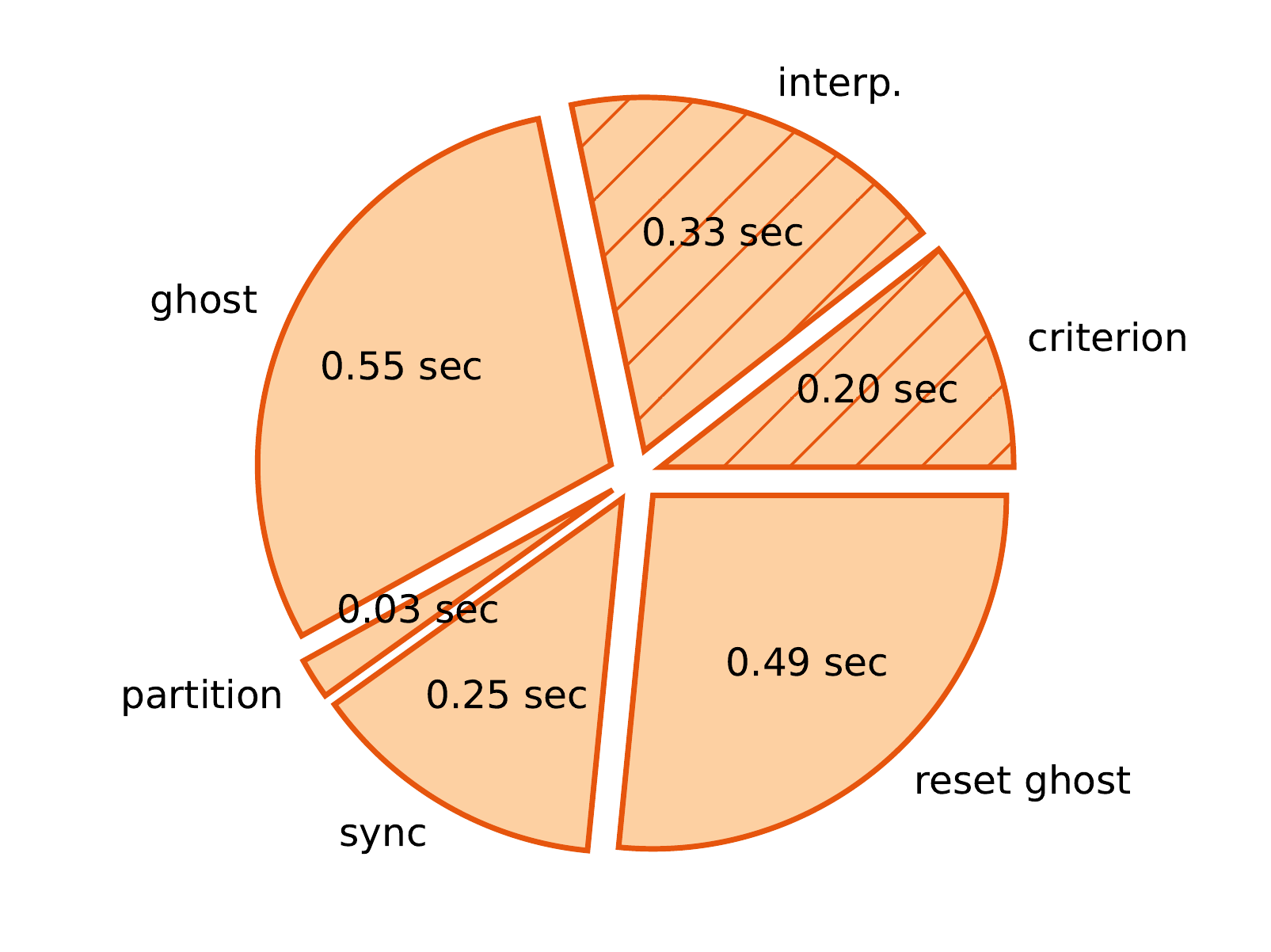}
	\else
    	\includegraphics[width=\textwidth]{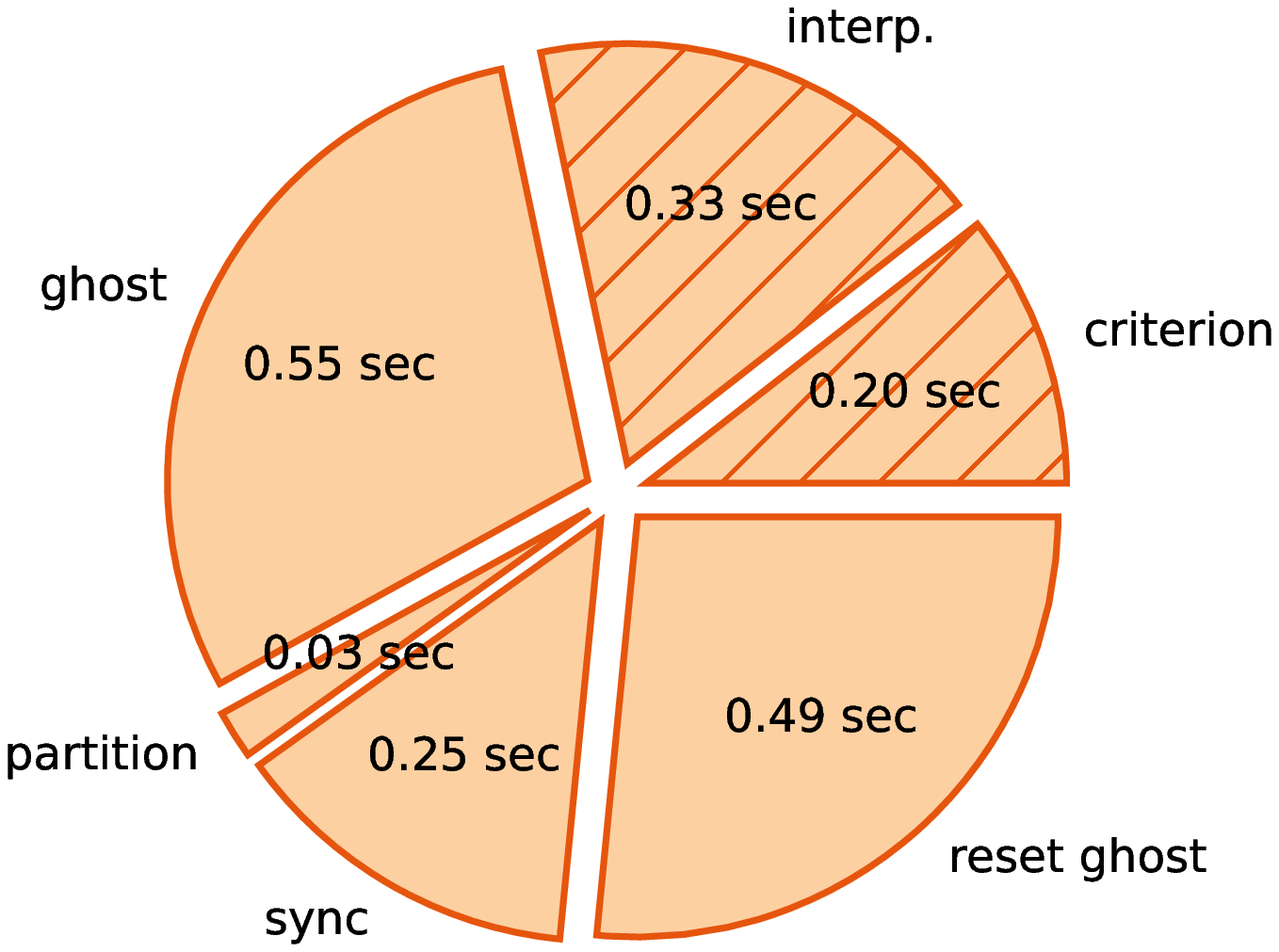}
	\fi
	\subcaption{$512$ node - $16,384$ ranks}
\end{minipage}%
\hfill%
\caption{Breakdown of the time spent during grid adaptation operations in a weak scalability test. The hatched regions indicates operations only involving computations and no communications.  The region marked ``\code{interp.}'' represents the coarsening/refinement of the blocks through the wavelets, ``\code{criterion}'' stands for the computation of the detail coefficients to evaluate the desired status of each block (refine, compress, or unaltered), ``\code{reset ghost}'' includes the reinitialization of the ghost meta-data including the \code{MPI\_Win\_create} and \code{MPI\_Win\_free} calls, ``\code{sync}'' contains the synchronizations and reductions on the block statuses to enforce the adaptation policy, ``\code{partition}'' is the load balancing of the grid and the re-distribution of the blocks, and ``\code{ghost}'' represents the computation of the ghost points needed for the coarsening/refinement of the blocks and the computation of the detail coefficients.}
\label{fig:weak_pie_adapt}
\end{figure}

In this section we provide more details on the scalability of the different sub-operations discussed in \sect{sec:scalability} of the main text. In \fref{fig:weak_pie_ghost} we show the breakdown of the timings during the ghost computation, where the hatched area represent purely computational operations that are expected to scale perfectly. In \fref{fig:weak_pie_stencil} we decompose similarly the operations involved in the stencil computation, and in \fref{fig:weak_pie_adapt} we do the same for the grid adaptation. 

